\documentclass[10pt]{article}
\usepackage{latexsym,amsmath,amssymb,graphics,amscd}
\makeatletter
\@addtoreset{figure}{section}
\def\thefigure{\thesection.\@arabic\c@figure}
\def\fps@figure{h,t}
\@addtoreset{table}{bsection}

\def\thetable{\thesection.\@arabic\c@table}
\def\fps@table{h, t}
\@addtoreset{equation}{section}

\makeatother

\begin{document}

\newtheorem{theorem}{Theorem}[section]
\newtheorem{definition}[theorem]{Definition}
\newtheorem{lemma}[theorem]{Lemma}
\newtheorem{remark}[theorem]{Remark}
\newtheorem{proposition}[theorem]{Proposition}
\newtheorem{corollary}[theorem]{Corollary}
\newtheorem{example}[theorem]{Example}
\newtheorem{examples}[theorem]{Examples}

\newcommand{\bfi}{\bfseries\itshape}

\newsavebox{\savepar}
\newenvironment{boxit}{\begin{lrbox}{\savepar}
\begin{minipage}[b]{15.8cm}}{\end{minipage}
\end{lrbox}\fbox{\usebox{\savepar}}}

\makeatletter
\title{{\bf Stochastic Hamiltonian dynamical systems}}
\author{Joan-Andreu L\'azaro-Cam\'{\i}$^1$ and Juan-Pablo Ortega$^2$}
\addtocounter{footnote}{1}
\footnotetext{Departamento de F\'{\i}sica Te\'orica. Universidad de Zaragoza. Pedro Cerbuna, 
12. E-50009 Zaragoza. Spain. {\texttt lazaro@unizar.es}}
\addtocounter{footnote}{1}
\footnotetext{Centre National de la Recherche Scientifique, D\'epartement de Math\'ematiques de
Besan\c con, Universit\'e de Franche-Comt\'e, UFR des Sciences et Techniques. 16, route de Gray.
F-25030 Besan\c con cedex. France. {\texttt Juan-Pablo.Ortega@univ-fcomte.fr }}

\date{}
\makeatother
\maketitle

\begin{abstract}  
We use the global stochastic analysis tools introduced by P. A. Meyer and L. Schwartz to
write down a stochastic generalization of the Hamilton equations on a Poisson
manifold that, for exact symplectic manifolds, are characterized by a natural critical
action principle similar to the one encountered in classical mechanics.
Several features and examples in relation with the solution
semimartingales of these equations are presented.
\end{abstract}

\medskip

\medskip

\medskip


\section{Introduction}

The generalization of classical mechanics to the context of stochastic dynamics has been an
active research subject ever since K. It\^o introduced  the theory of stochastic
differential equations in the 1950s (see for
instance~\cite{nelson 67, bismut 81, yasue 81, zambrini yasue 82, zheng meyer, thieullen zambrini
97, thieullen zambrini 97a, ludwig arnold, cresson darses 06, bou 1, bou 2}, and references therein). The motivations behind 
some pieces of work related to
this field lay in the hope that a suitable stochastic generalization of classical
mechanics should provide an explanation of the intrinsically random effects exhibited by
quantum mechanics within the context of the theory of diffusions . In other
instances the goal is establishing a
framework adapted to the handling of mechanical systems subjected to random
perturbations or whose parameters are not precisely determined and are hence
modeled as realizations of a random variable.

Most of the pieces of work in the first category use a class
of processes that have a stochastic derivative introduced in~\cite{nelson 67} and that has
been subsequently refined over the years. This derivative can be used to formulate a real
valued action and various associated variational principles whose extremals are the
processes of interest.

The approach followed in this paper is closer to the one introduced
in~\cite{bismut 81} in which the action has its image in the space of real valued
processes and the variations are taken in the space of processes with values in the
phase space of the system that we are modeling. This paper can be actually seen as a
generalization of some of the results in~\cite{bismut 81} in  the following directions:
\begin{description}
\item [(i)] We make extensive use of the global stochastic analysis tools introduced by
P. A. Meyer~\cite{meyer 81, meyer 82} and L. Schwartz~\cite{schwartz 82} to handle
non-Euclidean phase spaces. This feature not only widens the spectrum of systems that can
be handled but it is also of paramount importance at the time of reducing
them with respect to the symmetries that they may eventually have
(see~\cite{lo}); indeed, the orbit spaces obtained after reduction are generically
non-Euclidean, even if the original phase space is.
\item [(ii)] The stochastic dynamical components of the system are modeled by continuous
semimartingales and are not limited to Brownian motion.
\item [(iii)]  We handle stochastic Hamiltonian systems on Poisson manifolds and not only
on symplectic manifolds.
\item [(iv)] The variational principle that we propose in Theorem~\ref{Second Critical Action Principle} is not just satisfied by the stochastic Hamiltonian equations (as in~\cite{bismut 81}) but fully characterizes them.
\end{description}

There are various reasons that have lead us to consider these generalized Hamiltonian
systems. First, even though the laws that govern the dynamics of classical mechanical
systems  are, in principle, completely known,  the  finite precision of experimental
measurements yields impossible the estimation of the parameters of a particular given one
with total accuracy. Second, the modeling of complex physical systems involves most of the
time simplifying assumptions or idealizations of parts of the system, some
of which could be included in the description as a stochastic component; this
modeling philosophy has been extremely successful in the social
sciences~\cite{box jenkins 76}. Third, even if the model and the parameters
of the system are known with complete accuracy, the solutions of the
associated differential equations may be of great complexity and exhibit high
sensitivity to the initial conditions hence making the probabilistic
treatment and description of the solutions appropriate. Finally, we will see
(Section~\ref{The Langevin equation and viscous damping}) how stochastic
Hamiltonian modeling of microscopic systems can be used to model dissipation
and macroscopic damping.

The paper is structured as follows: in Section~\ref{The  stochastic Hamilton equations} 
we introduce the stochastic Hamilton equations with phase space a given
Poisson manifold and we study some of the fundamental properties of the solution
semimartingales like, for instance, the preservation of symplectic leaves or the
characterization of the conserved quantities. This section contains a discussion on two notions on non-linear stability, almost sure Lyapunov stability and stability in probability, that reduce in the deterministic setup to the standard definition of Lyapunov stability. We formulate criteria that generalize to the Hamiltonian stochastic context the standard energy methods to conclude the stability of a Hamiltonian equilibrium using existing conservation laws. More specifically, there are two different natural notions of conserved quantity in the stochastic context that, via a stochastic Dirichlet criterion (Theorem~\ref{Stochastic Dirichlet's Criterion}) allow one to conclude the different kinds of stability that we have mentioned above.
Section~\ref{examples section 1} contains
several examples: in the first one we show how the systems studied by Bismut
in~\cite{bismut 81}  fall in the category introduced in Section~\ref{The  stochastic
Hamilton equations}. We also see that a damped oscillator can be described as the average
motion of the solution semimartingale of a natural stochastic Hamiltonian
system, and that Brownian motion in a manifold is the projection onto the
base  space of very simple Hamiltonian stochastic semimartingale defined on
the cotangent bundle of the manifold or of its orthonormal frame bundle,
depending on the availability or not of a parallelization  for the manifold
in question.   Section~\ref{A critical action principle for the stochastic
Hamilton equations} is dedicated to showing that the stochastic Hamilton
equations are characterized by a critical action principle that generalizes the one
found in the treatment of deterministic systems. In order to make this part more readable, the proofs of most of the technical results needed to prove the  theorems in this section have been included separately at the end of the paper.

One of the goals of this paper is conveying to the geometric mechanics community the
plentitude of global tools available to handle mechanical problems that contain a stochastic
component and that do not seem to have been exploited to the full extent of their potential.
In order to facilitate the task of understanding the paper to non-probabilists we have
included an appendix that provides a self-contained presentation of some major facts in
stochastic calculus on manifolds needed for a first comprehension of our results. Those
pages are a very short and superficial presentation of a deep and technical field of
mathematics so the reader interested in a more complete account is encouraged to check with
the references quoted in the appendix and especially with the excellent
monograph~\cite{emery}. 

\medskip

\noindent {\bf Conventions:} All the manifolds in this paper are finite
dimensional, second-countable, locally compact, and Hausdorff (and hence
paracompact).

\section{The  stochastic Hamilton equations}
\label{The  stochastic Hamilton equations}

In this section we present a natural generalization of the standard Hamilton
equations in the stochastic context. Even though the arguments gathered in
the following paragraphs as motivation for these equations are of formal
nature, we will see later on that, as it was already the case for
the standard Hamilton equations, they satisfy a natural
variational principle.

We recall that a {\bfi symplectic manifold\/}
is a pair $(M,\,\omega)$, where $M$ is a manifold and $\omega\in\Omega^2(M)$ 
is a closed non-degenerate two-form on $M$, that is, $\mathbf{d}\omega = 0$ and, for every $m \in M$, the map $v \in T_mM \mapsto \omega(m)(v, \cdot) \in T^\ast_mM$ is a linear isomorphism  between the tangent space $T_mM$ to $M$ at $m$ and the cotangent space $ T^\ast_mM$. 
Using the nondegeneracy of the symplectic form
$\omega$,  one can associate each function $h \in  C^\infty(M) $ a vector field
$X_h\in\mathfrak{X}(M)$, defined by the equality
\begin{equation}
\label{hamilton intrinsic}
\mathbf{i}_{X_h}\omega=\mathbf{d} h.
\end{equation}
We will say that $X _h $ is the {\bfi  Hamiltonian vector field} associated to the {\bfi  Hamiltonian 
function } $h$. The expression~(\ref{hamilton intrinsic}) is referred to as the {\bfi 
Hamilton equations}.

A {\bfi Poisson manifold\/}
is a pair $(M,\,\{\cdot,\cdot\})$, where $M$ is a manifold 
and $\{\cdot,\cdot\}$ is a bilinear operation on $C^\infty(M)$
such that $(C^\infty(M),\,\{\cdot,\cdot\})$ is a Lie algebra
and $\{\cdot,\cdot\}$ is a derivation (that is, the Leibniz
identity holds) in each argument. 
The functions in the center $\mathcal{C}(M)$
of the Lie algebra $(C^\infty(M),\,\{\cdot ,\cdot \})$ are
called {\bfi Casimir functions\/}. From the natural isomorphism between derivations on 
$C^\infty(M)$
and  vector fields on $M$ it follows that
each $h\in C^\infty(M)$  induces a vector field on $M$ via the expression
$X_h = \{\cdot , h\}$,
called the {\bfi Hamiltonian vector field\/}
\index{Hamiltonian!vector field}
\index{vector!field!Hamiltonian}%
associated to  the {\bfi Hamiltonian function\/}
\index{Hamiltonian!function}
\index{function!Hamiltonian}%
$h$. Hamilton's equations $\dot z =  X_h(z)$ can be equivalently
written in Poisson bracket form as 
$\dot f = \{f, h\}$,
for any $f \in C^\infty(M)$. The derivation property of the Poisson bracket implies that for 
any two functions $f,\,g\in C^\infty(M)$, the value of the
bracket $\{f,\,g\}(z)$ at an arbitrary point $z\in M$ (and
therefore $X_f(z)$ as well), depends on $f$ only through $\mathbf{d}
f(z)$  which allows us to define  a contravariant
antisymmetric two--tensor $B\in\Lambda^2(M)$ by 
$B(z)(\alpha_z,\,\beta_z)=\{f,\,g\}(z)$,
where $\mathbf{d} f(z)=\alpha_z \in T^\ast _z M$ and $\mathbf{d} g(z)=\beta_z\in T^\ast _z
M$.  This tensor is called the {\bfi Poisson tensor\/} of $M$.
The vector bundle map $B^\sharp:T^\ast  M\rightarrow TM$ naturally associated to $B$ is
defined by
$
B(z)(\alpha_z,\,\beta_z)=
\langle\alpha_z,\,B^\sharp(\beta_z)\rangle$.

We start by rewriting the solutions of the standard Hamilton equations in a form that we will
be able to mimic in the stochastic differential equations context. All the necessary
prerequisites on stochastic calculus on manifolds can be found in a short review in the
appendix at the end of the paper.

\begin{proposition}
Let  $(M, \omega )$ be a symplectic manifold and $h \in  C^\infty(M) $. The
smooth curve 
$\gamma:[0,T]\rightarrow  M $ is an integral curve  of the Hamiltonian vector field $X _h $
if and only if for any $\alpha \in  \Omega (M) $ and for any $t \in  [0,T]$
\begin{equation}
\label{integral form of Hamiltonian solutions}
\int_{\gamma|_{[0,t]}} \alpha=-\int _0^t \mathbf{d} h(\omega ^{\sharp}(\alpha))\circ \gamma
(s) d s, 
\end{equation}
where $\omega^{\sharp}: T ^\ast M \rightarrow  TM $ is the vector bundle isomorphism induced
by $\omega$. More generally, if $M$ is a Poisson manifold with bracket $\{ \cdot , \cdot \}
$ then the same result holds with~(\ref{integral form of Hamiltonian solutions}) replaced by
\begin{equation}
\label{integral form of Hamiltonian solutions Poisson}
\int_{\gamma|_{[0,t]}} \alpha=-\int _0^t \mathbf{d} h(B ^{\sharp}(\alpha))\circ \gamma
(s) d s, 
\end{equation} 
\end{proposition}

\noindent\textbf{Proof.\ \ } Since in the symplectic case $\omega^{\sharp }=B ^{\sharp } $,
it suffices to prove~(\ref{integral form of Hamiltonian solutions Poisson}). 
As~(\ref{integral form of Hamiltonian solutions Poisson}) holds for
any $t \in [0,T]$, we can take derivatives with respect to $t$ on both sides and we obtain
the equivalent form 
\begin{equation}
\label{equivalent derivative}
\langle \alpha(\gamma (t)), \dot \gamma (t)\rangle=-\langle \mathbf{d} h (\gamma (t)),
B^\sharp(\gamma (t))(\alpha (\gamma (t)))\rangle.
\end{equation}
Let $f \in C^\infty(M) $ be such that
$\mathbf{d}f (\gamma (t))= \alpha (\gamma(t)) $. Then~(\ref{equivalent derivative})  can be
rewritten as 
\[
\langle\mathbf{d}f (\gamma (t)),\dot{\gamma}(t)\rangle=-\langle \mathbf{d} h (\gamma (t)),
B^\sharp(\gamma (t))(\mathbf{d}f (\gamma (t)))\rangle=\{f,h\}(\gamma(t)),
\]
which is equivalent to
$
\dot \gamma(t)= X _h (\gamma (t))$,
as required. \quad $\blacksquare$

\medskip

We will now introduce the stochastic Hamilton equations by mimicking in the context of
Stratonovich integration the integral expressions~(\ref{integral form of Hamiltonian
solutions}) and~(\ref{integral form of Hamiltonian
solutions Poisson}).  In the next definition we will use the following notation: let $f:M
\rightarrow  W $ be a differentiable function that takes values on the vector space $W $. We
define the {\bfi  differential} $ \mathbf{d} f: TM \rightarrow  W $ as the map given by
$\mathbf{d} f=p _2 \circ Tf $, where $Tf : TM \rightarrow  TW=W \times  W $ is the tangent
map of $f$ and  $p _2:W \times W \rightarrow  W $ is the projection onto the second factor.
If $W = \mathbb{R}$ this definition coincides with the usual differential. If $\{ e _1,
\ldots, e _n\}$ is a basis of  $W $ and $f=\sum _{i=1}^n f ^i e _i $ then $\mathbf{d}f =
\sum _{i=1}^n \mathbf{d}f ^i \otimes e _i $.

\begin{definition}
\label{hamiltonian semimartingale strato}
Let $(M, \{ \cdot , \cdot \} )$ be a Poisson manifold, $X: \mathbb{R}_+ \times  \Omega
\rightarrow   V
$ a semimartingale that takes values on the vector space $V$ with $X _0=0 $, and $h : M
\rightarrow  V ^\ast  $ a smooth function. Let $\{ \epsilon^1, \ldots, \epsilon^r\} $ be a
basis of $V ^\ast $ and  $h=\sum _{i=1}^r h _i \epsilon^i $. The {\bfi  Hamilton equations}
with {\bfi  stochastic component}
$X$, and {\bfi   Hamiltonian function}  $h $ are the Stratonovich stochastic differential
equation
\begin{equation}
\label{Hamiltonian equations differential form}
\delta\Gamma^h=H (X, \Gamma)\delta X,
\end{equation}
defined by the Stratonovich operator $H(v, z): T _v V   \rightarrow T _z M $ given by 
\begin{equation}
\label{stratonovich operator in full glory}
H(v, z)(u):=\sum _{j=1}^r\langle \epsilon^j,u\rangle X_{h _j}(z).
\end{equation}
\end{definition}

The dual Stratonovich operator $H ^\ast (v,z):T _z^\ast M \rightarrow T _v^\ast V  $ of
$H(v, z)$ is given by $H ^\ast (v,z)(\alpha_z)=- \mathbf{d}h (z) \cdot B
^{\sharp}(z)(\alpha_z)$. Hence, the 
results quoted in Appendix~\ref{Stochastic differential equations on manifolds} show
that for any $\mathcal{F} _0$ measurable random variable $\Gamma_0$,
there exists a unique semimartingale $\Gamma ^h $  such that  $\Gamma_0^h= \Gamma_0 $ and  a
maximal stopping time
$\zeta ^h
$  that solve~(\ref{Hamiltonian equations differential form}), that is, for any $\alpha\in
\Omega(M) $,
\begin{equation}
\label{solution Hamiltonian integral}
\int \langle \alpha, \delta \Gamma ^h\rangle=-\int  \langle\mathbf{d}h(
B^\sharp (\alpha))(\Gamma ^h), \delta X\rangle.
\end{equation}
We will refer to $\Gamma^h $ as the {\bfi  Hamiltonian semimartingale} associated to $h$
with initial condition $\Gamma _0$.

\begin{remark}
\normalfont
The stochastic component $X$ encodes the random behavior exhibited by the stochastic
Hamiltonian system that we are modeling and the Hamiltonian function $h$ specifies how it
embeds in its phase space. Unlike the situation encountered in the deterministic setup we
allow the Hamiltonian function to be vector valued in order to accommodate higher
dimensional stochastic dynamics. 
\end{remark}

\begin{remark}
\normalfont
The generalization of Hamilton's equations proposed in
Definition~\ref{hamiltonian semimartingale strato}  by using a Stratonovich
operator is inspired by one of the transfer principles presented
in~\cite{emery transfer} to provide stochastic versions of ordinary
differential equations.
This procedure can be also used to carry out a similar
generalization of the equations induced by a Leibniz bracket
(see~\cite{leibniz}).
\end{remark}

\begin{remark}
\normalfont
{\bf Stratonovich versus It\^o integration: }at the time of proposing the equations in
Definition~\ref{hamiltonian semimartingale strato} a choice has been made, namely, we have chosen Stratonovich integration instead of It\^o or other kinds of stochastic integration. The option that we took is motivated by the fact that by using Stratonovich integration, most of the geometric features underlying classical deterministic Hamiltonian mechanics are preserved in the stochastic context (see the next section). Additionally, from the mathematical point of view, this choice is the most economical one in the sense that the classical geometric ingredients of Hamiltonian mechanics plus a noise semimartingale suffice to construct the equations; had we used It\^o integration we would have had to provide a Schwartz operator (see Section~\ref{Stochastic differential equations on manifolds}) and the construction of such an object via a transfer principle like in~\cite{emery transfer} involves the choice of a connection.

The use of It\^o integration in the modeling of physical phenomena is sometimes preferred because the definition of this integral is not \emph{anticipative}, that is, it does not assume any knowledge about the behavior of the system in future times. Even though we have used Stratonovich integration to write down our equations, we also share this feature because the equations in  Definition~\ref{hamiltonian semimartingale strato} can be  naturally translated to the It\^o framework (see Proposition~\ref{equivalent in ito form}). This is a particular case of a more general fact since given any Stratonovich stochastic differential equation there always exists an equivalent It\^o stochastic differential equation, in the sense that both equations have the same solutions. Note that the converse is in general not true.
\end{remark}

\subsection{Elementary properties of the stochastic Hamilton's equations}

\begin{proposition}
\label{time evolution strato}
Let $(M, \{ \cdot , \cdot \} )$ be a Poisson manifold, $X: \mathbb{R}_+ \times  \Omega
\rightarrow   V
$ a semimartingale that takes values on the vector space $V$ with $X _0=0 $ and $h : M \rightarrow  V
^\ast  $ a smooth function. Let $\Gamma _0 $ be  a $\mathcal{F}_0 $ measurable random
variable and
$\Gamma ^h$ the Hamiltonian semimartingale associated to $h$
with initial condition $\Gamma _0$. Let $\zeta ^h $ be the corresponding maximal stopping
time. Then, for any stopping time $\tau< \zeta ^h  $, the
Hamiltonian semimartingale $\Gamma ^h  $ satisfies
\begin{equation}
\label{equation with Poisson brackets}
f (\Gamma^h_{\tau})-f (\Gamma^h_{0})=\sum _{j=1}^{r}\int _0 ^\tau\{f, h _j\}(\Gamma^h)\delta
X ^j,
\end{equation}
where $\{h _j\}_{j \in  \{1, \ldots, r\}}$ and $\{X^j\}_{j \in  \{1, \ldots, r\}}$ are the
components of $h$ and $X$ with respect to two given dual bases $\{e _1, \ldots, e _r\} $ and $\{
\epsilon ^1, \ldots, \epsilon^r \}$ of  $V $ and $V ^\ast  $, respectively.  
Expression~(\ref{equation with Poisson brackets})  can  be rewritten in differential
notation as
\begin{equation*}
\delta f (\Gamma^h)=\sum _{j=1}^{r}\{f, h
_j\}(\Gamma^h)\delta X ^j.
\end{equation*}
\end{proposition}

\noindent\textbf{Proof.\ \ } It suffices to take $\alpha = \mathbf{d} f  $ in~(\ref{solution
Hamiltonian integral}). Indeed, by~(\ref{properties strato})
\[
\int _0^\tau\langle \mathbf{d} f, \delta \Gamma ^h\rangle=f (\Gamma^h_{\tau})-f
(\Gamma^h_{0}).
\]  
At the same time
\begin{equation*}
-\int _0^\tau\langle\mathbf{d}h (
B^\sharp (\mathbf{d}  f))(\Gamma ^h), \delta X\rangle=-\sum_{j=1}^r\int _0 ^\tau
\langle  (\mathbf{d} h _j \otimes\epsilon ^j  ( B^{\sharp}(\mathbf{d} f)))(\Gamma ^h),
\delta X
\rangle=
\sum_{j=1}^r\int _0 ^\tau
\langle  \{ f, h _j \}(\Gamma ^h)\epsilon ^j,
\delta X
\rangle.
\end{equation*}
By the second statement in~(\ref{properties strato}) this equals 
$\sum_{j=1}^r\int _0 ^\tau
\{ f, h _j \}(\Gamma ^h)
\delta  \left(\int \langle \epsilon^j, \delta X\rangle \right)$. Given that $\int \langle
\epsilon^j, \delta X\rangle= X  ^j -X _0 ^j $, the equality follows. \quad $\blacksquare$

\begin{remark}
\label{relation classical equations 1}
\normalfont
Notice that if in Definition~\ref{hamiltonian semimartingale strato} we take
$V ^\ast  = \mathbb{R} $, $h \in C^\infty(M)$,  and $X:
\mathbb{R}_+ \times  \Omega \rightarrow  \mathbb{R} $ the deterministic process given by
$(t, \omega )\longmapsto t $, then the stochastic Hamilton equations~(\ref{solution
Hamiltonian integral}) reduce to 
\begin{equation}
\label{reduces to classical case}
\int \langle \alpha, \delta \Gamma^h\rangle=\int \langle \alpha, X _h
\rangle \left(\Gamma^h_t \right) d t.
\end{equation}
A straightforward application of~(\ref{equation with Poisson brackets})
shows that  $\Gamma^h_t(\omega) $ is necessarily a
differentiable curve, for any $\omega\in \Omega $, and hence the
Riemann-Stieltjes integral in the left hand side of~(\ref{reduces to
classical case})  reduces, when evaluated at a given $\omega \in \Omega $, to
a Riemann integral identical to the one in the left hand side
of~(\ref{integral form of Hamiltonian solutions Poisson}), hence proving
that~(\ref{reduces to classical case})  reduces to the standard Hamilton
equations.

Indeed, let $\Gamma^h _{t_0}(\omega) \in  M$ be an arbitrary point in the
curve
$\Gamma^h_t(\omega) $, let
$U$ be a coordinate patch around $\Gamma^h _{t_0}(\omega) $ with coordinates
$\{ x ^1,
\ldots, x ^n\} $,  and let $x (t)=(x ^1(t),
\ldots, x ^n(t)) $ be the expression of $\Gamma^h _t(\omega) $ in these
coordinates. Then by~(\ref{equation with Poisson brackets}), for $h \in 
\mathbb{R}$  sufficiently small, and $i \in \{1, \ldots, n\} $,
\[
x ^i(t _0+ h)- x ^i(t _0)=\int_{t _0}^{t _0+h}\{x ^i, h\}(x (t))d t.
\]
Hence, by the Fundamental Theorem of Calculus, $x ^i(t)$ is differentiable
at $t _0$, with derivative
\[
\dot x ^i (t _0)=\lim\limits_{h \rightarrow 0}\frac{1}{h} \left( x ^i(t _0+
h)- x ^i(t _0)\right)=\lim\limits_{h \rightarrow 0}\frac{1}{h} \left(\int_{t
_0}^{t _0+h}\{x ^i, h\}(x (t))d t\right)=\{x ^i, h\}(x (t_0)),
\]
as required.
\end{remark}

\noindent The following proposition provides an equivalent  expression of the Stochastic Hamilton
equations in the It\^o form (see Section~\ref{Stochastic differential equations on manifolds}).

\begin{proposition}
\label{equivalent in ito form}
The stochastic Hamilton's equations in  Definition~\ref{hamiltonian semimartingale
strato} admit an equivalent description using It\^o integration by using  the Schwartz operator ${\mathcal H}(v,m): \tau_v V  \rightarrow  \tau _m M $
naturally associated to the Hamiltonian Stratonovich operator
$H$ and  that can be described as follows. Let $L \in\tau_{v}M$ be a second order vector and $f\in
C^{\infty}\left(  M\right)$  arbitrary, then
\[
\mathcal{H}\left(  v,m\right)  \left(  L \right)  \left[  f\right]
=\left\langle \sum_{i,j=1}^r\{f, h _j\} (m) \epsilon ^j+\{\{f,h _j\}, h _i\}
(m)\epsilon^i \cdot \epsilon^j, L \right\rangle.
\]
Moreover, expression~(\ref{equation with Poisson brackets}) in the It\^o representation is
given by
\begin{equation}
\label{f in ito form increment}
f\left(  \Gamma_{\tau}^{h}\right)  -f\left(  \Gamma_{0}^{h}\right)
=\sum_{j=1}^{r}\int_{0}^{\tau}\left\{  f,h_{j}\right\}  \left(  \Gamma
^{h}\right)  dX^{j}+\frac{1}{2}\sum_{j,i=1}^{r}\int_{0}^{\tau}\left\{
\left\{  f,h_{j}\right\}  ,h_{i}\right\}  \left(  \Gamma^{h}\right)  d\left[
X^{j},X^{i}\right].
\end{equation}
We will refer to ${\mathcal H} $  as the {\bfi  Hamiltonian Schwartz operator} associated to $h $.

\end{proposition}

\noindent\textbf{Proof.\ \ } According to the remarks made in the Appendix~\ref{Stochastic
differential equations on manifolds},  the Schwartz operator ${\mathcal H}$ naturally
associated to $H $ is constructed as follows. For any second order vector $L_{\ddot{v}}
\in\tau_{v}M$ associated to the acceleration of a
curve $v\left(  t\right)  $ in $V$ such that
$v\left( 0\right)  =v$ we define
$\mathcal{H}\left(  v,m\right)  \left(
L_{\ddot{v}}\right):=L_{\ddot{m}\left(  0\right)  }\in\tau_{m}M$, where $m\left(  t\right) 
$ is a curve in $M$ such that
$m\left(  0\right)  =m$ and $\dot{m}\left(  t\right)  =H\left(  v\left(
t\right)  ,m\left(  t\right)  \right)  \dot{v}\left(  t\right)  $, for $t$ in a
neighborhood of $0.$ 
Consequently, 
\begin{align*}
\mathcal{H}\left(  v,m\right)  \left(  L_{\ddot{v}}\right)  \left[  f\right]
&  =\left.\frac{d^2}{dt^2}\right|_{t=0}f\left(  m\left(  t\right)
\right) 
= \left.\frac{d}{dt}\right|_{t=0} \langle\mathbf{d}f (m (t)) ,\dot{m}(t)\rangle
=\left.\frac{d}{dt}\right|_{t=0} \langle\mathbf{d}f (m (t)),H(v (t), m (t)) 
\dot{v }(t)\rangle \\ 
 &= \left.\frac{d}{dt}\right|_{t=0} \sum_{j=1}^{r}\langle\epsilon ^j, \dot{v}(t)
\rangle \langle\mathbf{d}f(m (t)),X_{h _j}(m (t))\rangle=\left.\frac{d}{dt}\right|_{t=0}
\sum_{j=1}^{r}\langle\epsilon ^j, \dot{v}(t)
\rangle  \{f, h _j\}(m (t))\\
&=\sum_{j=1}^{r}\langle\epsilon ^j, \ddot{v}(0)\rangle\{f,h _j\} (m)+\langle\epsilon ^j,
\dot{v}(0)\rangle\langle \mathbf{d}\{ f,h _j\} (m) , \dot{m}(0)\rangle
\\
&=\sum_{j=1}^{r}\langle\epsilon ^j, \ddot{v}(0)\rangle\{f,h _j\} (m)+\langle\epsilon ^j,
\dot{v}(0)\rangle \sum_{i=1}^r\langle\epsilon ^i,\dot{v}(0)\rangle\{\{f, h _j\},h _i\} (m)
\\
&=\left\langle\sum_{i,j=1}^{r}\{f,h _j\}
(m)\epsilon ^j+\{\{f, h _j\},h _i\}
(m) \epsilon^i \cdot  \epsilon^j,L_{\ddot v}\right\rangle.
\end{align*}

In order to establish~(\ref{f in ito form increment}) 
we need to calculate $\mathcal{H}^{\ast}\left(  v,m\right)  (d_{2}f(m))$ for
a second order form $d_{2}f(m)\in\tau_{m}^{\ast}M$ at $m\in M,$ $f\in C^{\infty
}\left(  M\right)  $. Since $\mathcal{H}^{\ast}\left(  v,m\right) 
(d_{2}f(m))$ is fully
characterized by its action on elements of the form
$L_{\ddot{v}}\in\tau_{v}V$ for some curve $v\left(  t\right)  $ in $V$ such
that $v\left(  0\right)  =v$, we have 
\begin{align*}
\left\langle \mathcal{H}^{\ast}\left(  v,m\right)  (d_{2}f(m)),L_{\ddot{v}%
}\right\rangle  &  =\left\langle d_{2}f(m),\mathcal{H}\left(  v,m\right)
(L_{\ddot{v}})\right\rangle =\mathcal{H}\left(  v,m\right)  (L_{\ddot{v}})\left[
f\right]  \\
&  =\left\langle\sum_{i,j=1}^{r}\{f,h _j\}
(m)\epsilon ^j+\{\{f, h _j\},h _i\}
(m) \epsilon^i \cdot  \epsilon^j,L_{\ddot v}\right\rangle.
\end{align*}
Consequently,
$
\mathcal{H}^{\ast}\left(  v,m\right)  (d_{2}f(m))=\sum_{i,j=1}^{r}\{f,h _j\}
(m)\epsilon ^j+\{\{f, h _j\},h _i\}
(m) \epsilon^i \cdot  \epsilon^j$.

Hence, if  $\Gamma _h$ is the Hamiltonian semimartingale associated to $h$
with initial condition $\Gamma _0$, $\tau < \zeta ^h  $ is any 
stopping time, and $f \in  C^\infty(M)$, we have by~(\ref{properties strato})
and~(\ref{integral with quadratic variation}) 
\begin{align*}
f\left(  \Gamma_{\tau}^{h}\right)  -f\left(  \Gamma_{0}^{h}\right)   &
=\int_{0}^{\tau}\left\langle d_{2}f,d\Gamma^{h}\right\rangle 
 =\int_{0}^{\tau}\left\langle
\mathcal{H}^{\ast}\left(  X,\Gamma^{h}\right) (d_{2}f),dX\right\rangle \\
&  =\sum_{j=1}^{r}\int_{0}^{\tau}\left\langle \left\{  f,h_{j}\right\}
\left(  \Gamma^{h}\right)  \epsilon^{j},dX\right\rangle
+\sum_{j,i=1}^{r}\int_{0}^{\tau}\left\langle \left\{  \left\{  f,h_{j}%
\right\}  ,h_{i}\right\}  \left(  \Gamma^{h}\right)  
\epsilon^{i}\cdot\epsilon^{j}  ,dX\right\rangle \\
&  =\sum_{j=1}^{r}\int_{0}^{\tau}\left\{  f,h_{j}\right\}  \left(  \Gamma
^{h}\right)  dX^{j}+\frac{1}{2}\sum_{j,i=1}^{r}\int_{0}^{\tau}\left\{
\left\{  f,h_{j}\right\}  ,h_{i}\right\}  \left(  \Gamma^{h}\right)  d\left[
X^{i},X^{j}\right].\quad \blacksquare
\end{align*}

\begin{proposition}[Preservation of the symplectic leaves by  Hamiltonian
semimartingales]
In the setup of Definition~\ref{hamiltonian semimartingale strato}, let
$\mathcal{L} $  be a symplectic leaf of 
$(M,\omega)$ and $\Gamma^h $ a  
Hamiltonian semimartingale with initial condition   $\Gamma _0 (\omega)= Z _0$, where $Z _0$ is
a random variable such that $Z _0(\omega)\in \mathcal{L}$ for all $\omega \in  \Omega$.  
Then, for any 
stopping time $\tau< \zeta^h $ we have that $\Gamma _\tau ^h \in  \mathcal{L} $.
\end{proposition}

\noindent\textbf{Proof.\ \ } Expression~(\ref{stratonovich operator in full glory}) shows
that for any $z \in  \mathcal{L}$, the Stratonovich operator $H(v,z)$ takes values in the
characteristic distribution associated to the Poisson structure $(M, \{\cdot , \cdot
\})$, that is, in the tangent space $T \mathcal{L} $ of $\mathcal{L}$. Consequently, $H $
induces another Stratonovich operator $H _{\mathcal L}(v,z): T _v V \rightarrow  T _z
\mathcal{L}$, $v \in  V $, $z \in  \mathcal{L}$, obtained from $H $ by restriction  of
its range. It is clear that if $i : \mathcal{L} \hookrightarrow M $ is the 
inclusion then
\begin{equation}
\label{one and other operator}
H _{\mathcal L}^\ast (v,z)\circ  T ^\ast  _z i=H ^\ast (v,z).
\end{equation}  
Let  $\Gamma _{\mathcal L} ^h $ be the semimartingale in $\mathcal{L}$ that is a solution
of the Stratonovich stochastic differential equation
\begin{equation}
\label{equation for l}
\delta \Gamma _{\mathcal L}^h=H _{\mathcal L}(X, \Gamma_{\mathcal L}^h)\delta X
\end{equation}
with initial condition $\Gamma_0$. We now show that $\overline{\Gamma}:=i \circ
\Gamma_{\mathcal L}^h $ is a solution of 
\[
\delta\overline{\Gamma}=H(X, \overline{\Gamma})\delta X.
\]
The uniqueness of the solution of a stochastic differential equation will guarantee in
that situation that $\Gamma^h $ necessarily coincides with $\overline{\Gamma} $, hence
proving the statement. Indeed, for any $\alpha\in \Omega (M) $,
\[
\int \langle \alpha, \delta\overline{\Gamma}\rangle=\int \langle \alpha, \delta(i
\circ \Gamma _{\mathcal L}^h)\rangle=\int \langle T ^\ast i \cdot \alpha, \delta
\Gamma_{\mathcal L} ^h\rangle.
\]
Since $\Gamma_{\mathcal L}^h$ satisfies~(\ref{equation for l}) and $T ^\ast i \cdot \alpha
\in  \Omega(\mathcal{L})$, by~(\ref{one and other operator}) this equals
\[
\int \langle H _{\mathcal L} ^\ast (X, \Gamma^h_{\mathcal L})(T ^\ast i \cdot
\alpha), \delta X\rangle=\int
\langle H ^\ast (X, i \circ \Gamma^h_{\mathcal L})(\alpha), \delta X
\rangle=\int
\langle H ^\ast (X, \overline{ \Gamma})(\alpha), \delta X
\rangle,
\] 
that is, $\delta\overline{\Gamma}=H(X, \overline{\Gamma})\delta X $, as required. \quad
$\blacksquare$

\begin{proposition}[The stochastic Hamilton equations in Darboux-Weinstein coordinates]
Let $(M, \{ \cdot , \cdot \})$ be a Poisson manifold and  $\Gamma ^h  $ be a solution of the
Hamilton equations~(\ref{Hamiltonian equations differential form}) with initial condition $x
_0 \in  M $. There exists an open neighborhood $U$ of $x _0 $ in $M$ and a stopping time
$\tau _U $ such that   $\Gamma_t ^h  (\omega) \in  U$, for any $\omega\in  \Omega
$ and any  $t\leq \tau_U(\omega)$. Moreover, $U$ admits local Darboux coordinates $(q ^1,
\ldots, q ^n, p _1, \ldots, p _n, z _1, \ldots, z _l )$ in which~(\ref{equation with Poisson
brackets}) takes the form
\begin{eqnarray*}
q ^i(\Gamma^h _\tau)-q ^i(\Gamma^h _0)&=&\sum _{j=1}^{r}\int _0 ^\tau \frac{\partial h
_j}{\partial p _i}\delta X ^j,\\
p _i(\Gamma^h _\tau)-p _i(\Gamma^h _0)&=&-\sum _{j=1}^{r}\int _0 ^\tau \frac{\partial h
_j}{\partial q^i}\delta X ^j,\\  
z _i(\Gamma^h _\tau)-z _i(\Gamma^h _0)&=&\sum _{j=1}^{r}\int
_0 ^\tau \{z _i, h _j\}_T\delta X ^j,
\end{eqnarray*} 
where $\{ \cdot , \cdot \} _T $ is the transverse Poisson structure of $(M, \{ \cdot , \cdot
\})$ at $x _0 $.
\end{proposition}

\noindent\textbf{Proof.\ \ } Let $U$ be an open neighborhood of $x _0  $ in $M$  for which
Darboux coordinates can be chosen. Define $\tau _U=\inf _{t \geq 0}\{ \Gamma_t ^h \in U ^c\}
$ ($\tau _U$ is the exit time of $U$). It is a standard fact in the theory of stochastic
processes that $\tau _U $ is a stopping time. The proposition follows by
writing~(\ref{equation with Poisson brackets}) for the Darboux-Weinstein coordinate functions
$(q ^1,
\ldots, q ^n, p _1, \ldots, p _n, z _1, \ldots, z _l )$. \quad $\blacksquare$

\medskip

Let $\zeta:M\times\Omega\rightarrow\left[  0,\infty\right]  $ be the map such
that, for any $z\in M,$ $\zeta\left(  z\right)  $ is the maximal stopping time
associated to the solution of the stochastic Hamilton equations~(\ref{Hamiltonian equations
differential form})  with initial condition $\Gamma_{0}=z$ a.s.. Let $F$ be the {\bfi  flow}
of~(\ref{Hamiltonian equations differential form}), that is, for any $z\in M$, $F\left( 
z\right)  :\left[  0,\zeta\left( z\right)  \right] 
\rightarrow M$ is the solution semimartingale of~(\ref{Hamiltonian equations differential form})
with initial condition $z$. The map $z \in M \longmapsto  F _t (z, \omega) \in  M$ is 
a local diffeomorphism  of  $M$, for each $t \geq 0  $ and almost all $\omega \in  \Omega $  in
which this map is defined (see~\cite{ikeda watanabe}). In the following result, we show that, in
the symplectic context, Hamiltonian flows preserve the symplectic form and hence the associated
volume form
$\theta=\omega\wedge\overset{n}{...}\wedge\omega$. This has already been shown for
Hamiltonian diffusions (see Example~\ref{bismut random mechanical}) by Bismut~\cite{bismut 81}.

\begin{theorem}[Stochastic Liouville's Theorem]
\label{Stochastic Liouville's Theorem}
Let $\left(  M,\omega\right)  $ be a symplectic
manifold, $X:\mathbb{R}_{+}\times\Omega\rightarrow V^{\ast}$ a semimartingale,
and $h:M\rightarrow V^{\ast}$ a Hamiltonian function. Let $F$ be the
associated Hamiltonian flow. Then, for any $z\in M$ and any $\left(
t,\eta\right)  \in\left[  0,\zeta\left(  z\right)  \right]  $,
\[
F_{t}^{\ast}\left(  z,\eta\right)  \omega=\omega.
\]
\end{theorem}

\medskip

\noindent\textbf{Proof.\ \ }
By \cite[Theorem 3.3]{Kunita} (see also \cite{Watanabe}), given an arbitrary
form $\alpha\in\Omega^{k}\left(  M\right)  $ and $z\in M$, the process
$F\left(  z\right)  ^{\ast}\alpha$ satisfies the following stochastic
differential equation:%
\[
F\left(  z\right)  ^{\ast}\alpha=\alpha\left(  z\right)  +\sum_{j=1}^{r}\int
F\left(  z\right)  ^{\ast}\left(  \pounds _{X_{h_{j}}}\alpha\right)  \delta
X^{j}.
\]
In particular, if $\alpha=\omega$ then $\pounds _{X_{h_{j}}}\omega=0$ for any
$j\in\left\{  1,...,r\right\}  $, and hence the result follows.
\quad $\blacksquare$

\subsection{Conserved quantities and stability}

Conservation laws in Hamiltonian mechanics are extremely important since they make easier the integration of the systems that have them and, in some instances, provide qualitative information about the dynamics. A particular case of this is their use in concluding the nonlinear stability of certain equilibrium solutions using Dirichlet type criteria that we will generalize to the stochastic setup using the following definitions.

\begin{definition}
A function $f\in C^{\infty}\left(  M\right)  $ is said to be a 
{\bfi strongly} (respectively, {\bfi  weakly}) {\bfi  conserved quantity} of the stochastic Hamiltonian system associated to $h: M
\rightarrow  V ^\ast $ if for any solution
$\Gamma^{h}$ of the stochastic Hamilton equations~(\ref{Hamiltonian equations
differential form}) we have that
$f\left( \Gamma^{h}\right)    =f\left(
\Gamma_{0}^{h}\right)  $ (respectively,  $E[f\left( \Gamma^{h}_{\tau}\right)]  =E[f\left(
\Gamma_{0}^{h}\right) ] $, for any stopping time $\tau $).
\end{definition}

\noindent Notice that strongly conserved quantities are obviously weakly conserved and that the two definitions  coincide for deterministic systems with the standard definition of conserved quantity.
The following result provides in  the stochastic setup an analogue of the
classical  characterization of the conserved quantities in terms of Poisson involution
properties.

\begin{proposition}
\label{prop 1} Let $(M,\{ \cdot , \cdot \})$ be a Poisson manifold, $X:\mathbb{R}%
_{+}\times\Omega\rightarrow V$ a semimartingale that takes values on the
vector space $V$ such that $X _0=0 $, and $h:M\rightarrow V^{\ast}$ and $f\in C^{\infty}(M)$
two smooth functions. If $\left\{  f,h_{j}\right\}  =0$ for every component
$h_{j}$ of $h$ then $f$ is a strongly conserved quantity of the stochastic Hamilton 
equations~(\ref{Hamiltonian equations differential form}).

Conversely, suppose that the semimartingale $X=\sum_{j=1}^{r}X^{j}%
\epsilon_{j}$ is such that $\left[  X^{i},X^{j}\right]  =0$ if $i\neq j$.
If $f$ is a strongly conserved quantity then $\left\{  f,h_{j}\right\}  =0$,
for any $j\in\left\{  1,...,r\right\}  $ such that $\left[  X^{j}
,X^{j}\right]  $ is an strictly increasing process at $0$. The last condition
means that there exists $A \in\mathcal{F}$ and $\delta>0 $ with $P (A)>0 $
such that for any $t< \delta$ and $\omega\in A $ we have $[X^{j}, X ^{j} ]_{t}
(\omega)>[X^{j}, X ^{j} ]_{0} (\omega)$, for all $j\in\left\{  1,...,r\right\}  $.
\end{proposition}

\noindent\textbf{Proof.\ \ } Let $\Gamma ^h$ be the Hamiltonian semimartingale associated to
$h$ with initial condition $\Gamma _0^h$. As we saw in
(\ref{f in ito form increment}),%
\begin{equation}
\label{decomposition to be used}
f\left(  \Gamma^{h}\right)  =f\left(  \Gamma^{h}_0\right)  +\sum_{j=1}^{r}
\int\left\{  f,h_{j}\right\}  \left(  \Gamma^{h}\right)  dX^{j}+\frac{1}
{2}\sum_{j,i=1}^{r}\int\left\{  \left\{  f,h_{j}\right\}  ,h_{i}\right\}
\left(  \Gamma^{h}\right)  d\left[  X^{i},X^{j}\right].
\end{equation}
If $\left\{  f,h_{j}\right\}  =0$ for every component
$h_{j}$ of $h$ then all the integrals in the previous expression vanish and therefore
$f\left( 
\Gamma^{h}\right)  =f\left(  \Gamma^{h}_0\right)$ which implies that $f$ is a strongly
 conserved quantity of the Hamiltonian stochastic equations associated to $h$.

Conversely, suppose now that $f$ is a  strongly conserved quantity. This implies that for
any initial condition $\Gamma ^h _0$, the semimartingale $f\left( 
\Gamma^{h}\right)  $ is actually time independent and hence of finite variation.
Equivalently, the (unique) decomposition of $f\left( 
\Gamma^{h}\right)  $ into two processes, one of finite variation plus a local martingale,
only has the first term. In order to isolate the local martingale term of $f\left( 
\Gamma^{h}\right)  $ recall first that the quadratic variations $\left[ 
X^{i},X^{j}\right] $ have finite variation and that the integral with respect to a finite
variation process has finite variation (see~\cite[Proposition 4.3]{legall integration}).
Consequently, the last summand in~(\ref{decomposition to be used}) has finite variation. As
to the second summand, let $M^{j}$ and $A ^j $, $j=1,\ldots, r $, local martingales and
finite variation processes, respectively, such that  $X^{j}=A^{j}+M^{j}$. Then,
\[
\int\left\{  f,h_{j}\right\}  \left(  \Gamma^{h}\right)  dX^{j}=\int\left\{ 
f,h_{j}\right\}  \left(  \Gamma^{h}\right)  dM^{j}+\int\left\{  f,h_{j}\right\}  \left( 
\Gamma^{h}\right)  dA^{j}.
\]
Given that for each $j$, $\int\left\{  f,h_{j}\right\}  \left( 
\Gamma^{h}\right)  dA^{j} $  is a finite variation process and $\int\left\{ 
f,h_{j}\right\}  \left(  \Gamma^{h}\right)  dM^{j} $ is a local martingale 
(see~\cite[Theorem 29, page 128]{protter})
we conclude that $Z:=\sum_{j=1}^{r}
\int\left\{  f,h_{j}\right\}  \left(  \Gamma^{h}\right)  dM^{j} $ is the local martingale
term of $f\left( 
\Gamma^{h}\right)  $ and hence equal to zero.

We notice now that any continuous local
martingale $Z:\mathbb{R}_{+}\times\Omega\rightarrow\mathbb{R}$ is also  a
local $L^{2}\left(  \Omega\right)  $-martingale. Indeed, consider the sequence of
stopping times
$\tau^{n}=\left\{  \inf t\geq0\mid\left\vert Z_{t}\right\vert =n\right\}$, $ n \in 
\mathbb{N} $. Then $E\left[  \left(  Z^{\tau^{n}
}\right)  _{t}^{2}\right]  \leq E\left[  n^{2}\right]  =n^{2}$, for all $t \in 
\mathbb{R}_+$. Hence, $Z^{\tau^{n}}\in L^{2}\left(  \Omega\right)  $ for any $n$. In
addition, $E\left[  \left(  Z^{\tau^{n}}\right)  _{t}^{2}\right]  =E\left[
\left[  Z^{\tau^{n}},Z^{\tau^{n}}\right]  _{t}\right]  $ (see~\cite[Corollary
3, page 73]{protter}). On the other
hand by Proposition~\ref{lemma restringir intervalo temporal}, 
\[
Z^{\tau^{n}}=\left(  \sum_{j=1}^{r}\int\left\{  f,h_{j}\right\}  \left(  \Gamma^{h}\right)
dM^{j}\right)  ^{\tau^{n}}=\sum_{j=1}^{r}\int\mathbf{1}_{\left[  0,\tau
^{n}\right]  }\left\{  f,h_{j}\right\}  \left(  \Gamma^{h}\right)  dM^{j}.
\]
Thus, by~\cite[Theorem 29, page 75]{protter} and the hypothesis $\left[  X^{i},X^{j}\right] 
=0$ if $i\neq j$,
\begin{align*}
E\left[  \left(  Z^{\tau^{n}}\right)  _{t}^{2}\right]   &  =E\left[  \left[
Z^{\tau^{n}},Z^{\tau^{n}}\right]  _{t}\right]  =\sum_{j,i=1}^{r}E\left[
\left[  \int\mathbf{1}_{\left[  0,\tau^{n}\right]  }\left\{  f,h_{j}\right\}
\left(  \Gamma^{h}\right)  dM^{j},\int\mathbf{1}_{\left[  0,\tau^{n}\right]
}\left\{  f,h_{i}\right\}  \left(  \Gamma^{h}\right)  dM^{i}\right]
_{t}\right]  \\
&  =\sum_{j,i=1}^{r}E\left[  \left(  \int\mathbf{1}_{\left[  0,\tau
^{n}\right]  }\left(  \left\{  f,h_{j}\right\}  \left\{  f,h_{i}\right\}
\right)  \left(  \Gamma^{h}\right)  d\left[  M^{j},M^{i}\right]  \right)
_{t}\right]  \\
&  =\sum_{j,i=1}^{r}E\left[  \left(  \int\mathbf{1}_{\left[  0,\tau
^{n}\right]  }\left(  \left\{  f,h_{j}\right\}  \left\{  f,h_{i}\right\}
\right)  \left(  \Gamma^{h}\right)  d\left[  X^{j},X^{i}\right]  \right)
_{t}\right]  \\
&  =\sum_{j=1}^{r}E\left[  \left(  \int\mathbf{1}_{\left[  0,\tau^{n}\right]
}\left\{  f,h_{j}\right\}  ^{2}\left(  \Gamma^{h}\right)  d\left[  X^{j}%
,X^{j}\right]  \right)  _{t}\right].
\end{align*}
Since $\left[  X^{j},X^{j}\right]  $ is an increasing process of
finite variation then $\int\mathbf{1}_{\left[  0,\tau^{n}\right]  }\left\{
f,h_{j}\right\}  ^{2}\left(  \Gamma^{h}\right)  d\left[  X^{j},X^{j}\right]  $
is a Riemann-Stieltjes integral and hence for any $\omega \in \Omega $
\[
\left(  \int\mathbf{1}_{\left[  0,\tau^{n}\right]  }\left\{  f,h_{j}\right\}
^{2}\left(  \Gamma^{h}\right)  d\left[  X^{j},X^{j}\right]  \right)  \left(
\omega\right)  =\int\mathbf{1}_{\left[  0,\tau^{n}\left(  \omega\right)
\right]  }\left\{  f,h_{j}\right\}  ^{2}\left(  \Gamma^{h}\left(
\omega\right)  \right)  d\left(  \left[  X^{j},X^{j}\right]  \left(
\omega\right)  \right).
\]
As $\left[  X^{j},X^{j}\right]  \left(  \omega\right)  $ is an increasing
function of $t\in\mathbb{R}_{+}$, then  for any $j \in \{1, \ldots , r\} $
\begin{equation}
\label{expected zero we will see}
E\left[  \int\mathbf{1}_{\left[  0,\tau^{n}\right]  }\left\{
f,h_{j}\right\}  ^{2}\left(  \Gamma^{h}\right)  d\left[  X^{j},X^{j}\right]
\right]  \geq 0.
\end{equation}
Additionally, since $E\left[  \left(  Z^{\tau^{n}}\right)  _{t}^{2}\right]=0 $, we
necessarily have that  the inequality in~(\ref{expected zero we will see}) is actually
an equality. 
Hence, 
\begin{equation}
\label{equal zero in integral}
\int_{0}^{t}\mathbf{1}_{\left[  0,\tau^{n}\right]  }\left\{  f,h_{j}\right\}
^{2}\left(  \Gamma^{h}\right)  d\left[  X^{j},X^{j}\right]  =0.
\end{equation}

Suppose now that $[X^{j},X^{j}]$ is strictly increasing at $0$ for a particular
$j.$ Hence, there exists $A\in\mathcal{F}$ with $P(A)>0$, and $\delta>0$ such
that $[X^{j},X^{j}]_{t}\left(  \omega\right)  >[X^{j},X^{j}]_{0}\left(
\omega\right)  $ for any $t<\delta.$ Take now a fixed $\omega\in A$. Since
$\tau^{n}\rightarrow\infty$ a.s., we can take $n$ large enough to ensure that
$\tau^{n}\left(  \omega\right)  >t,$ where $t\in\lbrack0,\delta)$. Thus, we
may suppose that $\mathbf{1}_{\left[  0,\tau^{n}\right]  }\left(
t,\omega\right)  =1.$ As $\left[  X^{j},X^{j}\right]  \left(  \omega\right)
$ is an strictly increasing process at zero  $\int_{0}^{t}\left\{
f,h_{j}\right\}  ^{2}\left(  \Gamma^{h}\left(  \omega\right)  \right)
d\left[  X^{j},X^{j}\right]  \left(  \omega\right)  >0$ unless $\left\{
f,h_{j}\right\}  ^{2}\left(  \Gamma^{h}\left(  \omega\right)  \right)  =0$ in
a neighborhood $[0,\widetilde{\delta}_{\omega})$ of $0$ contained in
$[0,\delta)$. In principle $\widetilde{\delta}_{\omega}>0$ might depend on
$\omega\in A,$ so the values of $t\in\lbrack0,\delta)$ for which $\left\{
f,h_{j}\right\}  ^{2}\left(  \Gamma_{t}^{h}\left(  \omega\right)  \right)  =0$
for any $\omega\in A$ are those verifying $0\leq t\leq\inf_{\omega\in
A}\widetilde{\delta}_{\omega}.$ In any case~(\ref{equal zero in integral}) allows us to
conclude that
$\left\{  f,h_{j}\right\} ^{2}\left(  \Gamma_{0}^{h}  \left(  \omega\right) \right)  =0$ for
any
$\omega\in A.$ Finally, consider any $\Gamma^{h}$ solution to the Stochastic
Hamilton equations with constant initial condition $\Gamma_{0}^h=m\in M$ an arbitrary
point. Then, for any $\omega\in A,$%
\[
0=\left\{  f,h_{j}\right\}  ^{2}\left(  \Gamma_{0}^{h}\left(  \omega\right)
\right)  =\left\{  f,h_{j}\right\}  ^{2}\left(  m\right)  .
\]
Since $m\in M$ is arbitrary we can conclude that  $\left\{  f,h_{j}\right\}  =0$.
\quad$\blacksquare$

\medskip

We now use the conserved  quantities of a system in order to formulate sufficient Dirichlet type stability criteria. Even though the statements that follow are enounced for processes that are not necessarily Hamiltonian, it is for these systems that the criteria are potentially most useful. We start by spelling out the kind of nonlinear stability that we are after.

\begin{definition}
\label{stability definition stochastic}
Let $M$ be a manifold and  let 
\begin{equation}
\label{Stratonovich to add}
\delta \Gamma= e(X,\Gamma) \delta X  
\end{equation}
 be a Stratonovich stochastic differential equation whose solutions $\Gamma: \mathbb{R} \times \Omega \rightarrow M $ take values on $M$. Given $x \in M  $ and $s \in \mathbb{R}$, denote by $\Gamma^{s,x} $ the unique solution of~(\ref{Stratonovich to add}) such that $\Gamma^{s,x}_s (\omega)=x $, for all $\omega \in \Omega $. Suppose that the point $z _0 \in M $ is an {\bfi  equilibrium} of ~(\ref{Stratonovich to add}), that is, the constant process $\Gamma_t(\omega):= z _0$, for all $t \in \mathbb{R}$ and $\omega \in \Omega $, is a solution of~(\ref{Stratonovich to add}). Then we say that the equilibrium $z _0 $ is
\begin{description}
\item [(i)] {\bfi Almost surely (Lyapunov) stable} when for any open neighborhood $U$ of $z _0 $ there exists another neighborhood $V \subset U $ of $z _0 $ such that for any $z \in V $ we have $\Gamma^{0,z} \subset U $, a.s.
\item [(ii)] {\bfi  Stable in probability.} For any $s\geq 0 $ and $\epsilon>0 $
\begin{equation*}
\lim_{x \rightarrow z _0 }P\left\{ \sup_{t>s} d\left(\Gamma^{s,x}_t, z _0 \right)> \epsilon\right\}=0,
\end{equation*}
where $d:M \times M \rightarrow \mathbb{R} $ is any distance function that generates the manifold topology of  $M$.
\end{description}
\end{definition}

\begin{theorem}[Stochastic Dirichlet's Criterion]
\label{Stochastic Dirichlet's Criterion}
Suppose that we are in the setup of the previous definition and assume that there exists a function $f \in C^\infty(M)  $ such that $\mathbf{d}f(z _0)=0 $  and that the quadratic form $\mathbf{d}^2 f (z _0)$ is (positive or negative) definite. If $f$ is a strongly (respectively, weakly) conserved quantity  for the solutions of~(\ref{Stratonovich to add}) then the equilibrium $z _0 $ is almost surely stable (respectively, stable in probability).
\end{theorem}

\noindent\textbf{Proof.\ \ } Since the stability of the equilibrium $z _0 $ is a local statement, we can work in a chart of $M$ around $z _0$ with coordinates $(x _1, \ldots, x _n )$ in which $z _0$ is modeled by the origin. Moreover, using the Morse lemma and the hypotheses on the function $f$, and assuming without loss of generality that $f (z _0)=0 $, we choose the coordinates $(x _1, \ldots, x _n )$ so that $f (x _1, \ldots, x _n )= x _1^2+ \cdots+ x _n^2$. Hence, in the definition of stability in probability, we can use the distance function $d(x, z _0)= f (x) $.

Suppose now that $f$ is a strongly conserved quantity and let $U $ be an open neighborhood of $z _0$. Let $r>0 $ be such that  $V:=f ^{-1}([0,r]) \subset U $. Let $z \in  V  $ with $f (z)=r' $. As $f$ is a strongly conserved quantity $f (\Gamma^{0,z})=r'\leq r $ and hence $\Gamma^{0,z} \subset U $, as required. 

In order to study the case in which  $f$ is a weakly conserved quantity, let $\epsilon>0 $ and let $U _\epsilon $ be the ball of radius $\epsilon $ around $z _0$. Then, for any $x \in U _\epsilon $ and $s \in \mathbb{R}_+ $, let $\tau_{U _\epsilon}$ be  the first exit time of  $\Gamma^{s,x}$ with respect to $U _\epsilon$. Notice first that if $\omega \in \Omega  $ belongs to the set $\{ \omega\in \Omega\mid \sup _{0\leq s<t} d\left(\Gamma^{s,x}_t, z _0 \right)> \epsilon \}=\{ \omega\in  \Omega\mid \sup _{0\leq s<t} f\left(\Gamma^{s,x}_t \right)> \epsilon^2 \}$, then $\tau_{U _\epsilon}(\omega)\leq t $ and hence the stopped process $\left(\Gamma^{s,x} \right)^{\tau_{U _\epsilon}} $ satisfies that 
\[
f \left( \left(\Gamma^{s,x} \right)^{\tau_{U _\epsilon}} _t (\omega)\right)= f \left(\Gamma^{s,x}_{\tau_{U _\epsilon}(\omega)}(\omega) \right)= \epsilon ^2,
\]
for those values of $\omega$. This ensures that
\begin{equation*}
\epsilon^2 {\bf 1}_{\{ \omega\in \Omega\mid \sup _{0\leq s<t} d\left(\Gamma^{s,x}_t, z _0 \right)> \epsilon \}}\leq f \left( \left(\Gamma^{s,x} \right)^{\tau_{U _\epsilon}} _t \right).
\end{equation*}
Taking expectations in both sides of this inequality we obtain
\begin{equation*}
P \left(\sup _{0\leq s<t} d\left(\Gamma^{s,x}_t, z _0 \right)> \epsilon\right)\leq \frac{E[f \left( \left(\Gamma^{s,x} \right)^{\tau_{U _\epsilon}} _t \right)]}{\epsilon^2}.
\end{equation*}
Since by hypothesis $f $ is a weakly conserved quantity, we can rewrite the right hand side of this inequality as
\begin{equation*}
\frac{E[f \left( \left(\Gamma^{s,x} \right)^{\tau_{U _\epsilon}} _t \right)]}{\epsilon^2}=\frac{E\left[f \left( \Gamma^{s,x}_{\tau_{U _\epsilon} \wedge t} \right)\right]}{\epsilon^2}= \frac{E[f(\Gamma^{s,x}_s)]}{\epsilon^2}=\frac{f (x)}{\epsilon^2},
\end{equation*}
and we can therefore conclude that 
\begin{equation}
\label{final inequality stability}
P \left(\sup _{0\leq s<t} d\left(\Gamma^{s,x}_t, z _0 \right)> \epsilon\right)\leq \frac{f (x)}{ \epsilon^2}.
\end{equation}
Taking the limit $x\rightarrow  z _0 $ in this expression and recalling that  $f (z _0)= 0 $, the result follows. \quad $\blacksquare$

\medskip

A careful inspection of the proof that we just carried out reveals that in order for~(\ref{final inequality stability}) to hold, it would suffice to have $E[f\left( \Gamma_{\tau}\right)]  \leq E[f\left(
\Gamma_{0}\right) ] $, for any stopping time $\tau $ and any solution $\Gamma  $, instead of the equality guaranteed by the weak conservation condition. This motivates the next definition.

\begin{definition}
Suppose that we are in the setup of Definition~\ref{stability definition stochastic}. Let $U$ be an open neighborhood of the equilibrium $z _0$ and let $V:U \rightarrow \mathbb{R}$ be a continuous function. We say that $V$ is a {\bfi  Lyapunov function} for the equilibrium $z _0$ if $V(z _0)= 0 $, $V(z)>0 $ for any $z \in U\setminus \{ z _0\} $, and 
\begin{equation}
\label{inequality lyapunov}
E[V\left( \Gamma _{\tau}\right)]  \leq E[V\left(
\Gamma_{0} \right) ], 
\end{equation}
for any stopping time $\tau $ and any solution $\Gamma $ of~(\ref{Stratonovich to add}).
\end{definition}

This definition generalizes to the stochastic context the standard notion of Lyapunov function that one encounters in dynamical systems theory. If~(\ref{Stratonovich to add}) is the stochastic differential equation associated to an It\^o diffusion and the Lyapunov function is twice differentiable, the inequality~(\ref{inequality lyapunov}) can be ensured by requiring that $A[V](z)\leq 0 $, for any $z \in U\setminus \{ z _0\} $, where $A$ is the infinitesimal generator of the diffusion, and by using Dynkin's formula. 

\begin{theorem}[Stochastic Lyapunov's Theorem]
Let $z _0 \in  M$ be an equilibrium solution of the stochastic differential equation~(\ref{Stratonovich to add}) and let $V:U \rightarrow \mathbb{R}$ be a continuous Lyapunov function for $z _0$. Then $z _0 $ is stable in probability.
\end{theorem}

\noindent\textbf{Proof.\ \ } Let $U _\epsilon $ be the ball of radius $\epsilon $ around $z _0$ and let $V _\epsilon:=\inf _{x \in U\setminus U _\epsilon}V (x) $. Using the same notation as in the previous theorem we denote, for any $x \in U _\epsilon $ and $s \in \mathbb{R}_+ $, $\tau_{U _\epsilon}$ as  the first exit time of  $\Gamma^{s,x}$ with respect to $U _\epsilon$. Using the same approach as above we notice that if $\omega \in \Omega  $ belongs to the set $\{ \omega\in \Omega\mid \sup _{0\leq s<t} d\left(\Gamma^{s,x}_t, z _0 \right)> \epsilon \}$, then $\tau_{U _\epsilon}(\omega)\leq t $ and hence the stopped process $\left(\Gamma^{s,x} \right)^{\tau_{U _\epsilon}} $ satisfies that 
\[
V \left( \left(\Gamma^{s,x} \right)^{\tau_{U _\epsilon}} _t (\omega)\right)= V\left(\Gamma^{s,x}_{\tau_{U _\epsilon}(\omega)}(\omega) \right)\geq V _\epsilon,
\]
for those values of $\omega$, since $\Gamma^{s,x}_{\tau_{U _\epsilon}(\omega)}(\omega) $ belongs to the boundary of $U _\epsilon $. This ensures that
\begin{equation*}
V_\epsilon {\bf 1}_{\{ \omega\in \Omega\mid \sup _{0\leq s<t} d\left(\Gamma^{s,x}_t, z _0 \right)> \epsilon \}}\leq V \left( \left(\Gamma^{s,x} \right)^{\tau_{U _\epsilon}} _t \right).
\end{equation*}
Taking expectations in both sides of this inequality we obtain
\begin{equation*}
P \left(\sup _{0\leq s<t} d\left(\Gamma^{s,x}_t, z _0 \right)> \epsilon\right)\leq \frac{E[V \left( \left(\Gamma^{s,x} \right)^{\tau_{U _\epsilon}} _t \right)]}{V_\epsilon}.
\end{equation*}
We now use that $V$ being a Lyapunov function satisfies~(\ref{inequality lyapunov}) and hence 
\begin{equation*}
\frac{E[V \left( \left(\Gamma^{s,x} \right)^{\tau_{U _\epsilon}} _t \right)]}{V_\epsilon}=\frac{E\left[V\left( \Gamma^{s,x}_{\tau_{U _\epsilon} \wedge t} \right)\right]}{V_\epsilon}\leq \frac{E[V(\Gamma^{s,x}_s)]}{V_\epsilon}=\frac{V (x)}{V_\epsilon}.
\end{equation*}
We can therefore conclude that 
\begin{equation*}
P \left(\sup _{0\leq s<t} d\left(\Gamma^{s,x}_t, z _0 \right)> \epsilon\right)\leq \frac{V (x)}{ V_\epsilon}.
\end{equation*}
Taking the limit $x\rightarrow  z _0 $ in this expression and recalling that  $V (z _0)= 0 $, the result follows. \quad $\blacksquare$

\begin{remark}
\normalfont
This theorem has been proved by Gihman~\cite{gihman} and Hasminskii~\cite{hasminskii} for  It\^o  diffusions.
\end{remark}

\section{Examples}
\label{examples section 1}

\subsection{Stochastic perturbation of a
Hamiltonian mechanical system and Bismut's Hamiltonian diffusions}
\label{bismut random mechanical}
\normalfont
Let $(M, \{ \cdot , \cdot \})$ be a Poisson
manifold and
$h_{j}\in C^{\infty }\left(  M\right)  $, $j=0,...,r$, smooth functions. Let
$h: M  \longrightarrow \mathbb{R}^{r+1}$ be the Hamiltonian function 
$m  \longmapsto  \left(  h_{0}\left(  m\right)  ,\ldots,h_{r}\left(  m\right)
\right)$,
and consider the semimartingale
$X:\mathbb{R}_{+}\times\Omega\rightarrow\mathbb{R}^{r+1}$ given by 
$\left(  t,\omega\right)   \longmapsto  \left(  t,B_{t}^{1}\left(
\omega\right)  ,\ldots ,B_{t}^{r}\left(  \omega\right)  \right)$,
where $B^{j},$ $j=1,...,r$, are $r$-independent Brownian motions. 
L\'evy's characterization of Brownian motion shows  (see for instance~\cite[Theorem
40, page 87]{protter}) that   
$\left[  B^{j},B^{i}\right]_t  =t\delta^{ji}$. In this setup, the equation~(\ref{equation
with Poisson brackets}) reads
\begin{equation}
\label{increment for random}
f\left(  \Gamma_{\tau}^{h}\right)  -f\left(  \Gamma_{0}^{h}\right)  =\int
_{0}^{\tau}\left\{  f,h_{0}\right\}  \left(  \Gamma^{h}\right)  dt+\sum
_{j=1}^{r}\int_{0}^{\tau}\left\{  f,h_{j}\right\}  \left(  \Gamma^{h}\right)
\delta B^{j}
\end{equation}
for any $f\in C^{\infty}\left(  M\right)$. According to~(\ref{f in ito form increment}), the
equivalent It\^o version of this equation is 
\[
f\left(  \Gamma_{\tau}^{h}\right)  -f\left(  \Gamma_{0}^{h}\right)  =\int
_{0}^{\tau}\left\{  f,h_{0}\right\}  \left(  \Gamma^{h}\right)  dt+\sum
_{j=1}^{r}\int_{0}^{\tau}\left\{  f,h_{j}\right\}  \left(  \Gamma^{h}\right)
d B^{j}+\int_{0}^{\tau}\left\{\{  f,h_{j}\right\} , h _j\} \left(  \Gamma^{h}\right)
d t.
\]
Equation~(\ref{increment for random}) may be interpreted as a stochastic
perturbation of the classical Hamilton equations associated to $h _0$, that is,
\[
\frac{d(f \circ \gamma)}{dt}(t)=\left\{
f,h_{0}\right\}  \left(  \gamma\left(  t\right)  \right)  .
\]
by  the $r$
Brownian motions
$B^{j}$. These equations have been studied by Bismut in~\cite{bismut 81} in the particular
case in which the Poisson manifold
$(M, \{ \cdot , \cdot \})$  is just the symplectic Euclidean space
$\mathbb{R} ^{2n}$ with the canonical symplectic form. He refers to these particular processes
as {\bfi  Hamiltonian diffusions}.

If we apply Proposition~\ref{prop 1} to the stochastic Hamiltonian system~(\ref{bismut random
mechanical}), we obtain a generalization to Poisson manifolds of a result originally formulated
by Bismut (see~\cite[Th\'eor\`emes 4.1 and 4.2, page 231]{bismut 81}) for Hamiltonian diffusions.
See also~\cite{misawa}.

\begin{proposition}
Consider the stochastic Hamiltonian system introduced in~(\ref{bismut random
mechanical}). Then $f\in C^{\infty}\left(  M\right)  $ is a 
conserved quantity if and only if
\begin{equation}
\label{involution random}
\left\{  f,h_{0}\right\}  =\left\{  f,h_{1}\right\}  =\ldots=\left\{
f,h_{r}\right\}  =0.
\end{equation}
\end{proposition}

\noindent\textbf{Proof.\ \ } If~(\ref{involution random}) holds then $f$ is clearly a
 conserved quantity by Proposition~\ref{prop 1}. 
Conversely, notice that as $[B ^i, B
^j]=t \delta^{ij}$, $i,j \in  \{ 1, \ldots,r\} $, and $X ^0(t, \omega)=t $ is a finite
variation process then $[X ^i, X ^j]=0 $ for any $i,j \in  \{0, 1, \ldots,r\} $ such that 
$i\neq j $. Consequently, by Proposition~\ref{prop 1}, if $f$ is a  conserved
quantity then 
\begin{equation}
\label{to be contradict}
\left\{  f,h_{1}\right\}  =\ldots=\left\{  f,h_{r}\right\}  =0.
\end{equation}
Moreover,~(\ref{increment for random}) reduces to 
\[
\int_{0}^{\tau}\left\{  f,h_{0}\right\}  \left(  \Gamma^{h}\right)  dt=0,
\]
for any Hamiltonian semimartingale $\Gamma^h $ and any stopping time
$\tau\leq\zeta^h$. Suppose that
$\left\{  f,h_{0}\right\}  \left(  m_{0}\right)  >0$ for some $m_{0}\in M$. By continuity 
there exists a compact neighborhood $U$ of $ m _0$ such that $\left.  \left\{
f,h_{0}\right\}  \right\vert _{U}>0.$ Take $\Gamma^{h}$ the Hamiltonian semimartingale
with initial condition $\Gamma^{h}_0=m_{0},$ and let $\xi$ be the first exit time of
$U$ for
$\Gamma^{h}$. Then, defining
$\tau:=\xi\wedge\zeta,$
\[
\int_{0}^{\tau}\left\{  f,h_{0}\right\}  \left(  \Gamma^{h}\right)  dt\geq
\int_{0}^{\tau}\min\left\{  \left\{  f,h_{0}\right\}(m)\mid m \in  U  \right\} 
dt>0,
\]
which contradicts~(\ref{to be contradict}). Therefore, 
$\left\{  f,h_{0}\right\}  =0$ also, as required.
\quad $\blacksquare$

\begin{remark}
\normalfont
Notice that, unlike what happens for standard deterministic Hamiltonian systems, the energy $h
_0$ of a Hamiltonian diffusion does not need to be conserved if the other components of the
Hamiltonian are not involution with $h _0 $. This is a general fact about stochastic
Hamiltonian systems that makes them useful in the modeling of dissipative phenomena. We see
more of this in the next example.
\end{remark}

\subsection{Integrable stochastic Hamiltonian dynamical systems.}

Let $\left(  M,\omega\right)  $ be a $2n$-dimensional manifold, $X:\mathbb{R}
_{+}\times\Omega\rightarrow V$ a semimartingale, and $h:M\rightarrow V^{\ast}$
such that $h=\sum_{i=1}^{r}h_{i}\epsilon^{i}$, with $\left\{
\epsilon^{1},...,\epsilon^{r}\right\}  $ a basis of $V^{\ast}$. Let $H$
be the associated Stratonovich operator in~(\ref{stratonovich operator in full glory}).

Suppose that there exists a family of functions $\left\{  f_{r+1}%
,...,f_{n}\right\}  \subset C^{\infty}\left(  M\right)  $ such that the
$n$-functions
$
\left\{  f_{1}:=h_{1},...,f_{r}:=h_{r},f_{r+1},...,f_{n}\right\}  \subset
C^{\infty}\left(  M\right)
$
are in Poisson involution, that is, $\left\{  f_{i},f_{j}\right\}  =0$, for any
$i,j\in\left\{  1,...,n\right\}  $. Moreover, assume that $F:=\left(
f_{1},...,f_{n}\right)  $ satisfies the hypotheses of the Liouville-Arnold
Theorem~\cite{arnold}: $F$ has compact and connected fibers and its
components are independent. In this setup, we will say that the stochastic
Hamiltonian dynamical system associated to $H$ is {\bfi  integrable}.

As it was already the case for standard (Liouville-Arnold) integrable systems,
there is a symplectomorphism that takes $\left(  M,\omega\right)  $ to
$\left(  \mathbb{T}^{n}\times\mathbb{R}^{n},\sum_{i=1}^{n}\mathbf{d}\theta
^{i}\wedge\mathbf{d}I_{i}\right)  $ and for which $F\equiv F\left(
I_{1},...,I_{n}\right)  $. In particular, in the action-angle coordinates
$\left(I_{1},...,I_{n},  \theta^{1},...,\theta^{n}\right)  $, $h_{j}\equiv
h_{j}\left(  I_{1},...,I_{n}\right)  $ with $j\in\left\{  1,...,r\right\}  $.
In other words, the components of the Hamiltonian function depend only on the
actions ${\bf I}:=\left(  I_{1},...,I_{n}\right)  $. Therefore, for any random variable
$\Gamma_{0}$ and any $i\in\left\{ 1,...,n\right\}  $
\begin{subequations}
\label{action-angle}%
\begin{align}
I_{i}\left(  \Gamma\right)  -I_{i}\left(  \Gamma_{0}\right)   &  =\sum
_{j=1}^{r}\int\left\{  I_{i},h_{j}\left(  \mathbf{I}\right)  \right\}  \left(
\Gamma\right)  \delta X^{j}=0\label{action evolution}\\
\theta^{i}\left(  \Gamma\right)  -\theta^i\left(  \Gamma_{0}\right)   &
=\sum_{j=1}^{r}\int\left\{  \theta^{i},h_{j}\left(  \mathbf{I}\right)
\right\}  \left(  \Gamma\right)  \delta X^{j}=\sum_{j=1}^{r}\int\frac{\partial
h_{j}}{\partial I_{i}}\left(  \Gamma\right)  \delta X^{j}.
\label{angle evolution}%
\end{align}
Consequently, the tori determined by fixing $\mathbf{I}%
=\operatorname*{constant}$ are left invariant by the stochastic flow
associated to (\ref{action-angle}). In particular, as the paths of the
solutions are contained in compact sets, the stochastic flow is defined for
any time and the flow is complete. Moreover, the restriction of this
stochastic differential equation to the torus given by say, $\mathbf{I}_{0}$,
yields the solution%
\end{subequations}
\begin{equation}
\theta^{i}\left(  \Gamma\right)  -\theta^i\left(  \Gamma_{0}\right)  =\sum
_{j=1}^{r}\omega_{j}\left(  \mathbf{I}_{0}\right)  X^{j},
\label{angulo integrado}%
\end{equation}
where $\omega_{j}\left(  \mathbf{I}_{0}\right)  :=\frac{\partial h_{j}
}{\partial I _i}\left(  {\bf I} _0\right)  $ and where we have assumed that
$X_{0}=0$. Expression~(\ref{angulo integrado}) clearly resembles the integration that
can be carried out for deterministic integrable systems.

Additionally, the Haar measure $\mathbf{d}\theta^{1}\wedge...\wedge
\mathbf{d}\theta^{n}$ on each invariant torus is left invariant by the
stochastic flow (see Theorem~\ref{Stochastic Liouville's Theorem}  and \cite{Xue-Mei}). 
Therefore, if we can
ensure that there exists a unique invariant measure $\mu$ (for instance, if
(\ref{angulo integrado}) defines a non-degenerate diffusion on the torus
$\mathbb{T}^{n}$, the invariant measure is unique up to a multiplicative constant by the
compactness of
$\mathbb{T}^{n}$ (see~\cite[Proposition 4.5]{ikeda watanabe})) then $\mu$ coincides necessarily
with the Haar measure.

\subsection{The Langevin equation and viscous damping}
\label{The Langevin equation and viscous damping}

Hamiltonian stochastic differential equations can be used to model dissipation phenomena. The
simplest example in this context is the damping force experienced by a particle in motion in a
viscous fluid. This dissipative phenomenon is usually modeled using a force in Newton's second
law  that depends linearly on the velocity of the particle (see for instance~\cite[\S
25]{landau}). The standard microscopic description of this motion is carried out using the
{\bfi  Langevin} stochastic differential equation (also called the {\bfi  Orstein-Uhlenbeck}
equation) that says that the velocity $\dot q (t) $ of the particle with mass $m$ is a
stochastic process that solves the stochastic differential equation
\begin{equation}
\label{langevin equation}
m\, d \dot q(t)=- \lambda \dot q (t) d t+ b d B _t,
\end{equation}
where $\lambda >0$ is the damping coefficient, b is a constant, and $B _t $ is a Brownian
motion. A common physical interpretation for this equation (see~\cite{chorin hald}) is that the
Brownian motion models random instantaneous bursts of momentum that are added to the particle
by collision with lighter particles, while the mean effect of the collisions is the slowing
down of the particle. This fact is mathematically described by saying that the expected value
$q _e:= \mbox{\rm E}[q] $ of the process $q $ determined by~(\ref{langevin equation}) satisfies
the  ordinary differential equation $\ddot q _e=- \lambda \dot q _e $. Even though this
description is accurate it is not fully satisfactory given that it {\it does not provide any
information about the mechanism that links the presence of the Brownian perturbation to the
emergence of damping} in the equation. In order for the physical explanation to be complete, a
relation between the coefficients $b$ and $\lambda $ should be provided in such a way that {\it
the damping vanishes when the Brownian collisions disappear}, that is, $\lambda=0 $ when $b = 0
$.

We now show that the motion of a particle of mass $m$ in one dimension subjected to viscous
damping with coefficient $\lambda  $ and to a harmonic potential with Hooke constant $k$ is a
Hamiltonian stochastic differential equation. More explicitly, we will give a stochastic
Hamiltonian system such that the expected value $q _e$ of its solution
semimartingales satisfies the ordinary differential equation of the damped harmonic oscillator,
that is,
\[
m\,\ddot q_e (t)=- \lambda \dot q_e(t)-k q_e (t).
\]
This description provides a mathematical mechanism by which the stochastic perturbations in the
system generate an average damping.

Consider $\mathbb{R}^{2}$ with its canonical symplectic form and let
$X:\mathbb{R}_{+}\times\Omega\rightarrow\mathbb{R}$ be the real semimartingale
given by $X_{t}\left(  \omega\right)  =\left(  t+\nu B_{t}\left(
\omega\right)  \right)  $ with $\nu\in\mathbb{R}$ and $B_{t}$ a Brownian
motion. Let now $h:\mathbb{R}^{2}\rightarrow\mathbb{R}$ be the energy of a
harmonic oscillator, that is, $h\left(  q,p\right)  :=\frac{1}{2m}p^{2}%
+\frac{1}{2}\rho q^{2}$. By (2.10), the solution semimartingales $\Gamma^{h}$
of the Hamiltonian stochastic equations associated to $h$ and $X$ satisfy%
\begin{align}
q\left(  \Gamma^{h}\right)  -q\left(  \Gamma_{0}^{h}\right)   &  =\frac{1}%
{2m}\int\left(  2p\left(  \Gamma_{t}^{h}\right)  -\nu^{2}\rho q\left(
\Gamma_{t}^{h}\right)  \right)  dt+\frac{\nu}{m}\int p\left(  \Gamma_{t}%
^{h}\right)  dB_{t},\label{damped 1}\\
p\left(  \Gamma^{h}\right)  -p\left(  \Gamma_{0}^{h}\right)   &  =-\frac{\rho
}{2m}\int\left(  \nu^{2}p\left(  \Gamma_{t}^{h}\right)  +2mq\left(  \Gamma
_{t}^{h}\right)  \right)  dt-\nu\rho\int q\left(  \Gamma_{t}^{h}\right)
dB_{t}. \label{damped 2}%
\end{align}
Given that $E\left[  \int p\left(  \Gamma_{t}^{h}\right)  dB_{t}\right]
=E\left[  \int q\left(  \Gamma_{t}^{h}\right)  dB_{t}\right]  =0$, if we
denote%
\[
q_{e}\left(  t\right)  :=E\left[  q\left(  \Gamma_{t}^{h}\right)  \right]
\text{, \ \ \ \ \ \ \ \ }p_{e}\left(  t\right)  :=E\left[  p\left(  \Gamma
_{t}^{h}\right)  \right]  ,
\]
Fubini's Theorem guarantees that%
\begin{equation}
\dot{q}_{e}\left(  t\right)  =\frac{1}{m}p_{e}\left(  t\right)  -\frac{\nu
^{2}\rho}{2m}q_{e}\left(  t\right)  \text{ \ \ \ and \ \ \ }\dot{p}_{e}\left(
t\right)  =-\frac{\nu^{2}\rho}{2m}p_{e}\left(  t\right)  -\rho q_{e}\left(
t\right)  . \label{eq 7}%
\end{equation}
From the first of these equations we obtain that%
\[
p_{e}\left(  t\right)  =m\dot{q}_{e}+\frac{\nu^{2}\rho}{2}q_{e}%
\]
whose time derivative is%
\[
\dot{p}_{e}\left(  t\right)  =m\ddot{q}_{e}+\frac{\nu^{2}\rho}{2}\dot{q}_{e}.
\]
These two equations substituted in the second equation of (\ref{eq 7}) yield%
\begin{equation}
m\ddot{q}_{e}\left(  t\right)  =-\nu^{2}\rho\dot{q}_{e}\left(  t\right)
-\rho\left(  \frac{\nu^{4}\rho}{4m}+1\right)  q_{e}\left(  t\right)  ,
\label{eq 8}%
\end{equation}
that is, the expected value of the position of the Hamiltonian semimartingale
$\Gamma^{h}$ associated to $h$ and $X$ satisfies the differential equation of
a damped harmonic oscillator (\ref{eq 4}) with constants%
\[
\lambda=\nu^{2}\rho\text{ \ \ \ and \ \ \ }k=\rho\left(  \frac{\nu^{4}\rho
}{4m}+1\right)  .
\]
Notice that the dependence of the damping and elastic constants on the
coefficients of the system is physically reasonable. For instance, we see that
the more intense the stochastic perturbation is, that is, the higher $\nu$ is,
the stronger the damping becomes ($\lambda=\nu^{2}\rho$ increases). In
particular, if there is no stochastic perturbation, that is, if $\nu=0$, then
the damping vanishes, $k=\rho$ and (\ref{eq 8}) becomes the differential
equation of a free harmonic oscillator of mass $m$ and elastic constant $\rho$.

\medskip

\noindent {\bf The stability of the resting solution.} It is easy to see that the constant process $\Gamma _t( \omega)=(0,0) $, for all $t \in \mathbb{R} $ and $\omega \in \Omega $  is an equilibrium solution of~(\ref{damped 1}) and~(\ref{damped 2}). One can show using the stochastic Dirichlet's criterion (Theorem~\ref{Stochastic Dirichlet's Criterion}) that this equilibrium is almost surely Lyapunov stable since the Hamiltonian function $h$ is a strongly conserved quantity (by~(\ref{equation with Poisson brackets}))  that exhibits a critical point at the origin with definite Hessian.

\medskip

\noindent {\bf The Langevin equation.} In the previous paragraphs we succeeded in providing a microscopic Hamiltonian description of the harmonic oscillator subjected to Brownian perturbations whose macroscopic counterpart via expectations yields the equations of the damped harmonic oscillator. In view of this, is  such a stochastic Hamiltonian description available for the pure Langevin equation~(\ref{langevin equation})? The answer is no. More specifically, it can be easily shown (proceed by contradiction) that ~(\ref{langevin equation}) cannot be written as a stochastic Hamiltonian differential equation on $\mathbb{R}^2 $ with its canonical symplectic form with a noise semimartingale of the form $X _t(\omega)=(f _0(t, B _t),f _1(t, B _t)) $ and a Hamiltonian function $h(q,p)=(h_0(q,p), h_1(q,p))$, $f _0, f _1, h _0, h _1 \in C^{\infty}(\mathbb{R})$. Nevertheless, if we put aside for a moment the stochastic Hamiltonian category and we use It\^o integration, the Langevin equation can still be written in phase space, that is, 
\begin{equation}
\label{eq langevin 1}
dq_{t}   =v_{t}dt,\quad
dv_{t}   =-\lambda v_{t}dt+bdB_{t},
\end{equation}
as a stochastic perturbation of a deterministic system, namely, a free
particle whose evolution is given by the differential equations%
\begin{equation}
dq_{t}=v_{t}dt\text{ \ \ and \ \ }dv_{t}=0.\label{unperturbed}%
\end{equation}
Let $\left\{  u^{1},u^{2}\right\}  $ be global coordinates on
$\mathbb{R}^{2}$ associated to the canonical basis $\left\{  e_{1}%
,e_{2}\right\}  $ and consider the global basis $\left\{  d_{2}u^{i}%
,d_{2}u^{i}\cdot d_{2}u^{j}\right\}  _{i,j=1,2}$ of $\tau^{\ast}\mathbb{R}%
^{2}$. Define a dual Schwartz operator $\mathcal{S}^{\ast}\left(  x,\left(
q,v\right)  \right)  :\tau_{\left(  q,v\right)  }^{\ast}\mathbb{R}%
^{2}\longrightarrow\tau_{x}^{\ast}\mathbb{R}^2$  characterized by the  relations
\begin{equation*}
d_{2}q    \longmapsto v d_{2}u^{1},\quad
d_{2}v \longmapsto bd_{2}u^{2}-\lambda v \left(  d_{2}u^{2}\cdot
d_{2}u^{2}\right),
\end{equation*}
where $\left(  q,v\right)  \in\mathbb{R}^{2}$ is an arbitrary point in phase
space and $x\in\mathbb{R}^{2}$. If  $X:\mathbb{R}_{+}\times
\Omega\rightarrow\mathbb{R}^{2}$ is such that $X\left(  t,\omega\right)
=\left(  t,bB_{t}\left(  \omega\right)  \right)  $, for any $\left(
t,\omega\right)  \in\mathbb{R}_{+}\times\Omega$, it is immediate to see that
the It\^{o} equations associated to $\mathcal{S}^{\ast}$ and $X$ are
(\ref{eq langevin 1}). Moreover, if we set $b=0 $, that is, we switch off  the Brownian perturbation then we recover (\ref{unperturbed}), as required.

\subsection{Brownian motions on manifolds}

The mathematical formulation of Brownian motions (or Wiener processes) on
manifolds has been the subject of much research and it is a central topic in
the study of stochastic processes on manifolds (see~\cite[Chapter 5]{ikeda
watanabe},~\cite[Chapter V]{emery}, and references therein for a good general review
of this subject).

In the following paragraphs we show that  Brownian motions can be defined in a particularly
simple way using the stochastic Hamilton equations introduced in
Definition~\ref{hamiltonian semimartingale strato}. More specifically we
will show that Brownian motions on manifolds can be obtained as the
projections onto the base  space of very simple Hamiltonian stochastic
semimartingales defined on the cotangent bundle of the manifold or of its
orthonormal frame bundle, depending on the availability or not of a
parallelization  for the manifold in question. 

We will first present  the case in which the manifold in question is
parallelizable or, equivalently, when the coframe bundle on the manifold
admits a global section, for the construction is particularly simple in
this situation. The parallelizability hypothesis is verified by many important
examples. For instance, any Lie group is parallelizable; the spheres $S^{1}$,
$S^{3}$, and
$S^{7}$ are parallelizable too. At the end of the section we describe the
general case.

The notion of manifold valued Brownian motion that we will use is the
following. A $M$-valued process $\Gamma $ is called a Brownian motion on
$(M,g)$, with
$g$ a Riemannian metric on $M$, whenever $\Gamma $ is continuous and adapted
and for every $f \in C^\infty(M)$
\[
f(\Gamma)-f(\Gamma_0)-\frac{1}{2}\int \Delta_Mf(\Gamma) dt
\]
is a local martingale. We recall that the Laplacian $\Delta_{M}\left( 
f\right) 
$ is defined  as $\Delta_{M}\left(
f\right)  =\operatorname*{Tr}\left(  \operatorname*{Hess}  f 
\right)  $, for any
$f\in C^{\infty}\left(  M\right)  $, where $\operatorname*{Hess}f:=\nabla(\nabla
f)$, with $\nabla:\mathfrak{X}(M)\times\mathfrak{X}(M)\rightarrow
\mathfrak{X}(M)$, the Levi-Civita connection of $g$. $\operatorname*{Hess}f$
is a symmetric $(0,2)$-tensor such that for any $X,Y\in\mathfrak{X}(M)$,
\begin{equation}
\operatorname*{Hess}f(X,Y)=X\left[  g(\operatorname*{grad}f,Y)\right]
-g(\operatorname*{grad}f,\nabla_{X}Y). \label{property hessian}%
\end{equation}

\medskip

\noindent {\bf Brownian motions on parallelizable manifolds.} Suppose that
the
$n$-dimensional manifold
$\left(  M,g\right) 
$ is parallelizable and  let
$\left\{ Y_{1},...,Y_{n}\right\}  $ be a family of vector fields such that for each
$m\in  M$, $\left\{ Y_{1}(m),...,Y_{n}(m)\right\}  $ forms a basis of $T_{m}M$ (a
parallelization). Applying the Gram-Schmidt orthonormalization procedure if necessary, we
may suppose that this parallelization is orthonormal, that is,
$g\left( Y_{i},Y_{j}\right)  =\delta_{ij}$,  
for any $i,j=1,...,n$.

Using this structure we are going to construct a stochastic Hamiltonian system on the
cotangent bundle $T ^\ast M $ of $M$, endowed with its canonical symplectic structure, and
we will show that the projection of the solution semimartingales of this system onto
$M$ are
$M$-valued Brownian motions in the sense specified above. Let $X: \mathbb{R}_+ \times \Omega
\rightarrow 
\mathbb{R}^{n+1}$  be the semimartingale given by
$X(t, \omega):=(t,B _t ^1 (\omega ), \ldots, B _t ^n (\omega )) $,
where $B^{j},$ $j=1,...,n$, are $n$-independent  Brownian
motions and let $h=(h _0, h _1, \ldots, h _n): T ^\ast M \rightarrow \mathbb{R}^{n+1}$ be
the function whose components are given by
\begin{equation}
\label{hamiltonian functions parallelizable}
\begin{array}{cccc}
h_{0}:&T^{\ast}M & \longrightarrow & \mathbb{R}\\
	&\alpha _m & \longmapsto & -\frac{1}{2}\sum_{j=1}^{n}\langle
\alpha_m, \left(\nabla_{Y_{j}}Y_{j}\right)(m)\rangle
\end{array}
\quad \mbox{ and } \quad
\begin{array}
{cccc}
h_{j}:&T^{\ast}M & \longrightarrow & \mathbb{R}\\
	&\alpha_m & \longmapsto &\langle\alpha_m, Y_{j}(m)\rangle.
\end{array}
\end{equation}
We will now study the projection onto $M$ of Hamiltonian semimartingales $\Gamma^h$ that
have $X$ as stochastic component and $h$ as Hamiltonian function and will prove that they
are $M$-valued Brownian motions. In order to do so we will be particularly interested in the
{\bfi  projectable functions}
$f$ of $T ^\ast M$, that is, the functions $f \in C^\infty(T ^\ast M)$ that can be written
as
$f =
\overline{f} \circ  \pi $ with $ \overline{f} \in C^\infty(M) $  and $\pi: T ^\ast M
\rightarrow  M $ the canonical projection.

We start by proving that for any projectable function $f= \overline{f} \circ
\pi
\in C^\infty(T ^\ast M) $
\begin{equation}
\label{brackets with projectable}
\{f, h _0\}=g\left(\mbox{\rm grad}\,\overline{f}, 
-\frac{1}{2}\sum_{j=1}^{n}\nabla_{Y_{j}}Y_{j}\right)\quad \mbox{ and } \quad
\{f, h _j\}=g\left(\mbox{\rm grad}\,\overline{f}, 
Y_{j}\right),
\end{equation}
 and where $\{ \cdot , \cdot \} $ is the Poisson bracket
associated to the canonical symplectic form on $T ^\ast M $. Indeed, let $U$ a Darboux
patch for $T ^\ast M $ with associated coordinates
$\left(  q^{1}, \ldots, q ^n, p _1, \ldots,p_{n}\right)$ such that $ \{ 
q^{i},p_{j} \}  =\delta_{j}^{i}$. There exists functions $f ^k_j \in C^{\infty}(\pi(U))$,
with $k,j \in  \{ 1 , \ldots, n\} $ such that the
vector fields
may be locally written as
$Y_{j}=\sum_{k=1}^{n}f_{j}^{k}%
\frac{\partial}{\partial q^{k}}$. Moreover, $h_{j}\left(
q,p\right)  =\sum_{k=1}^{n}f_{j}^{k}\left(  q\right)  p_{k}$ and
\begin{align*}
\left\{  f,h_{j}\right\}  
&=\left\{
\overline{f}\circ \pi,\sum_{k=1}^{n}f_{j}^{k}p_{k}\right\}  
=\sum_{k=1}^{n}f_{j}^{k}\left\{
\overline{f}\circ \pi,p_{k}\right\} =\sum_{k,i=1}^{n}f_{j}^{k}\frac{\partial
(\overline{f}\circ \pi)}{\partial q^{i}}\left\{ q^{i},p_{k}\right\} 
=\sum_{k,i=1}^{n}f_{j}^{k}\delta_{k}^{i}\frac{\partial \overline{f}}{\partial
q^{i}}\\
	&=Y_{j}[\overline{f}]\circ \pi=g\left( 
\operatorname*{grad}\overline{f},Y_{j}\right)\circ  \pi,
\end{align*}
as required. The first equality in~(\ref{brackets with projectable}) is proved analogously.
Notice that the formula that we just proved shows that if $f$ is projectable then so is
$\{f, h _j\} $, with $j \in  \{ 1 , \ldots, n\} $. Hence, 
using~(\ref{brackets with projectable}) again and~(\ref{property hessian}) we obtain that
\begin{equation}
\label{bracket with projectable 2}
\left\{  \left\{  f,h_{j}\right\}  ,h_{j}\right\}  =Y_{j}\left[  g\left( 
\operatorname*{grad} \overline{f},Y_{j}\right)  \right] \circ \pi
=\operatorname*{Hess}\overline{f}\left(  Y_{j},Y_{j}\right) \circ 
\pi+g\left( 
\operatorname*{grad}\overline{f},\nabla_{Y_{j}}Y_{j}\right)\circ  \pi,
\end{equation}
for $j \in  \{ 1 , \ldots, n\} $. Now, using~(\ref{brackets with projectable})
and~(\ref{bracket with projectable 2}) in~(\ref{f in ito form increment}) we have shown
that for any projectable function $f= \overline{f} \circ\pi$, the
Hamiltonian semimartingale $\Gamma^h $ satisfies that 
\begin{equation}
\overline{f}\circ \pi\left(  \Gamma^{h} \right)-\overline{f}\circ
\pi\left( 
\Gamma_{0}^{h}\right)  =\sum_{j=1} ^{n}\int g\left( 
\operatorname*{grad}\overline{f},Y_{j}\right)  \left(\pi\circ  
\Gamma ^{h}\right) 
dB_{s}^{j}+\frac{1}{2}\sum_{j=1}^{n}\int\operatorname*{Hess}
\overline{f}\left(  Y_{j},Y_{j}\right)  \left( \pi\circ  \Gamma^{h}\right)
dt,\label{increment brownian}
\end{equation}
or equivalently
\begin{equation}
\overline{f}\circ \pi\left(  \Gamma^{h} \right)-\overline{f}\circ
\pi\left( 
\Gamma_{0}^{h}\right)  -\frac{1}{2}\int\Delta_M(\overline{f})
  \left( \pi\circ  \Gamma^{h}\right)
dt=\sum_{j=1} ^{n}\int g\left( 
\operatorname*{grad}\overline{f},Y_{j}\right)  \left(\pi\circ  
\Gamma ^{h}\right) 
dB_{s}^{j}.\label{increment brownian second version}
\end{equation}
Since $\sum_{i=1}^{n}\int
g\left(  \operatorname*{grad}\overline{f},Y_{j}\right)  (\overline{\Gamma}^{h})dB^{i}$ is a
local martingale (see \cite[Theorem 20, page 63]{protter}), $ \pi(\Gamma ^h)$ is a Brownian
motion.

\medskip

\noindent\textbf{Brownian motions on Lie groups. }Let now $G$ be a (finite
dimensional) Lie group with Lie algebra $\mathfrak{g}$ and assume that $G$
admits a bi-invariant metric $g$, for example when $G$ is Abelian or compact.
This metric induces a pairing in $\mathfrak{g}$ invariant with respect to the
adjoint representation of $G$ on $\mathfrak{g}$. Let $\{\xi_{1},\ldots,\xi
_{n}\}$ be an orthonormal basis of $\mathfrak{g}$ with respect to this
invariant pairing and let $\{\nu_{1},\ldots,\nu_{n}\}$ be the corresponding
dual basis of $\mathfrak{g}^{\ast}$. The infinitesimal generator vector fields
$\{\xi_{1G},\ldots,\xi_{nG}\}$ defined by $\xi_{iG}(h)=T_{e}L_{h}\cdot\xi$,
with $L_{h}:G\rightarrow G$ the left translation map, $h\in G$, $i\in
\{1,\ldots n\}$, are obviously an orthonormal parallelization of $G$, that is
$g(\xi_{iG},\xi_{jG}):=\delta_{ij}$. Since $g$ is bi-invariant then $\nabla
_{X}Y=\frac{1}{2}[X,Y]$, for any $X,Y\in\mathfrak{X}(G)$ (see~\cite[Proposition 9,
page 304]{neill 83}), and hence $\nabla_{\xi_{iG}}\xi_{iG}=0$. Therefore, in this
particular case the first component
$h_{0}$ of the Hamiltonian function introduced in~(\ref{hamiltonian functions
parallelizable}) is zero and we can hence take  $h_{G}=\left(
h_{1},...,h_{n}\right)  $ and $X_{G}=\left(  B_{t}^{1},...,B_{t}^{n}\right)  $
when we consider the Hamilton equations that define the Brownian motion with
respect to $g$. 

As a special case of the previous construction that serves as a particularly
simple illustration, we are going to explicitly build the \textbf{Brownian
motion on a circle}. Let $S^{1}=\{e^{i\theta}\mid\theta\in\mathbb{R}\}$ be the
unit circle. The stochastic Hamiltonian differential equation for the
semimartingale $\Gamma^{h}$ associated to $X:\mathbb{R}_{+}\times
\Omega\rightarrow\mathbb{R}$, given by $X_{t}(\omega):=B_{t}(\omega)$, and the
Hamiltonian function $h:TS^{1}\simeq S^{1}\times\mathbb{R}\rightarrow
\mathbb{R}$ given by $h(e^{i\theta},\lambda):=\lambda$, is simply obtained by
writing (\ref{increment brownian}) down for the functions $f_{1}(e^{i\theta
}):=\cos\theta$ and $f_{2}(e^{i\theta}):=\sin\theta$ which provide us with the
equations for the projections $X^{h}$ and $Y^{h}$ of $\Gamma^{h}$ onto the
$OX$ and $OY$ axes, respectively. A straightforward computation yields
\begin{equation}
dX^{h}=-Y^{h}dB-\frac{1}{2}X^{h}dt\quad\mbox{ and  }\quad dY^{h}=X^{h}%
dB-\frac{1}{2}Y^{h}dt,\label{brownian on the circle equations}%
\end{equation}
which, incidentally, coincides with the equations proposed in expression
(5.1.13) of \cite{oksendal}. A solution of
(\ref{brownian on the circle equations}) is $(X_{t}^{h},Y_{t}^{h})=(\cos
B_{t},\sin B_{t})$, that is, $\Gamma_{t}^{h}=e^{iB_{t}}$.

\medskip

\noindent {\bf Brownian motions on arbitrary manifolds.}
Let $\left(  M,g\right)  $ be a not necessarily parallelizable Riemannian
manifold. In this case we will reproduce the same strategy as in the
previous paragraphs but replacing the cotangent bundle of the manifold by
the cotangent bundle of its orthonormal frame bundle.

Let $\mathcal{O}_{x}\left(  M\right)  $ be the set of orthonormal frames
for the tangent space $T_{x}M$. The orthonormal frame bundle $\mathcal{O}
\left(  M\right)  =\bigcup_{x\in M}\mathcal{O}_{x}\left(  M\right)  $ has a
natural smooth manifold structure of dimension $\left. 
n\left( n+1\right)  \right/  2$. We denote by $\pi:\mathcal{O}\left(  M\right)
\rightarrow M$ the canonical projection. We recall that a curve
$\gamma:\left( -\varepsilon,\varepsilon\right) 
\subset\mathbb{R}\rightarrow\mathcal{O}
\left(  M\right)  $ is called horizontal if $\gamma_{t}$ is the parallel
transport of
$\gamma_{0}$ along the projection $\pi\left(  \gamma_{t}\right)  $. The set of
tangent vectors of horizontal curves that contain a point
$u\in\mathcal{O}\left(  M\right)  $ defines the horizontal subspace
$H_{u}\mathcal{O}\left(  M\right)  \subset T_{u}\mathcal{O}\left(  M\right)
,$ with dimension $n$. The
projection $\pi:\mathcal{O}\left(  M\right)  \rightarrow M$ induces an
isomorphism $T_{u}\pi:H_{u}\mathcal{O}\left(  M\right)  \rightarrow
T_{\pi\left(  u\right)  }M$. On the orthonormal frame bundle, we have $n$
horizontal vector fields $Y_{i},$ $i=1,...,n$, defined as follows. For each
$u\in\mathcal{O}\left(  M\right)$, let $Y_{i}\left(  u\right)  $ be the
unique horizontal vector in $H_{u}\mathcal{O}\left(  M\right)  $ such that
$T_{u}%
\pi\left(  Y_{i}\right)  =u_{i},$ where $u_{i}$ is the $i$th unit vector of
the orthonormal frame $u$.
Now, given a smooth function $F\in C^{\infty}\left(  \mathcal{O}\left( 
M\right)
\right)  $, the operator%
\[
\Delta_{\mathcal{O}\left(  M\right)  }\left(  F\right)  =\sum_{i=1}^{n}%
Y_{i}\left[  Y_{i}\left[  F\right]  \right]
\]
is called Bochner's horizontal Laplacian on $\mathcal{O}\left(  M\right)  $.
At the same time, we recall that the Laplacian $\Delta_{M}\left(  f\right) 
$, for any
$f\in C^{\infty}\left(  M\right)  $, is defined  as $\Delta_{M}\left(
f\right)  =\operatorname*{Tr}\left(  \operatorname*{Hess}  f 
\right)  $. These two Laplacians are related by the relation
\begin{equation}
\Delta_{\mathcal{O}\left(  M\right)  }\left(  \pi^{\ast}f\right)  =\Delta
_{M}\left(  f\right),  \label{relation laplacians}%
\end{equation}
for any $f\in C^{\infty}\left(  M\right)  $ (see \cite{Hsu-Brownian motion}).

The Eells-Elworthy-Malliavin construction of Brownian motion can be summarized
as follows. Consider the following stochastic differential equation on
$\mathcal{O}\left(  M\right)  $ (see \cite{ikeda watanabe}):%
\begin{equation}
\delta U_{t}=\sum_{i=1}^{n}Y_{i}\left(  U_{t}\right)  \delta B_{t}%
^{i}\label{Brownian equation orthonormal}%
\end{equation}
where $B^{j},$ $j=1,...,n$, are $n$-independent  Brownian
motions.
Using the conventions introduced in the appendix~\ref{Stochastic
differential equations on manifolds} the expression~(\ref{Brownian equation
orthonormal})  is the Stratonovich stochastic differential equation
associated to the Stratonovich operator:
\[%
\begin{array}
[c]{cccc}%
e\left(  v,u\right)  : & T_{v}\mathbb{R}^{n} & \longrightarrow &
T_{u}\mathcal{O}\left(  M\right)  \\
& v=\sum_{i=1}^{n}v^{i}e_{i} & \longmapsto & \sum_{i=1}^{n}v^{i}Y_{i}\left(
u\right),
\end{array}
\]
where $\{e _1, \ldots, e _n\} $ is a fixed basis for $\mathbb{R}^n $. 
A solution of the stochastic differential equation~(\ref{Brownian equation
orthonormal})  is called a horizontal Brownian motion on
$\mathcal{O}\left(  M\right)  $ since, by the It\^{o} formula,
\[
F\left(  U\right)  -F\left(  U_{0}\right)  =\sum_{i=1}^{n}\int Y_{i}\left[
F\right]  \left(  U_{s}\right)  \delta B_{s}^{i}=\sum_{i=1}^{n}\int
Y_{i}\left[  F\right]  \left(  U_{s}\right)  dB_{s}^{i}+\frac{1}{2}\int
\Delta_{\mathcal{O}\left(  M\right)  }(F)\left(  U_{s}\right)  ds,
\]
for any $F\in C^{\infty}\left(  \mathcal{O}\left(  M\right)  \right)  $. In
particular, if $F=\pi^{\ast}\left(  f\right)  $ for some $f\in C^{\infty
}\left(  M\right)  ,$ by (\ref{relation laplacians})%
\[
f\left(  X\right)  -f\left(  X_{0}\right)  =\sum_{i=1}^{n}\int Y_{i}\left[
\pi^{\ast}\left(  f\right)  \right]  \left(  U_{s}\right)  dB_{s}^{i}+\frac
{1}{2}\int\Delta_{M}f\left(  X_{s}\right)  ds,
\]
where $X_{t}=\pi\left(  U_{t}\right)  $, which implies precisely that $X_{t}$
is a Brownian motion on $M$.

In order to generate (\ref{Brownian equation orthonormal}) as a Hamilton
equation, we introduce the functions $h_{i}:T^{\ast}   \mathcal{O}\left(
M\right)     \rightarrow\mathbb{R}$, $i=1,...,n$, given by $h_{i}\left(
\alpha\right)  =\langle\alpha, Y_{i}\rangle$. Recall that $T^{\ast} 
\mathcal{O}
\left(  M\right) $ being a cotangent bundle it has a canonical
symplectic structure. Mimicking the computations carried out in the
parallelizable case it can be seen that  the
Hamiltonian vector field
$X_{h_{i}}$ coincides with
$Y_{i}$ when acting on functions of the form
$F\circ\pi_{T^{\ast}\mathcal{O}\left(  M\right)  },$ where $F\in
C^{\infty}\left(  \mathcal{O}\left(  M\right) 
\right)  $ and
$\pi_{T^{\ast }\mathcal{O}\left(  M\right)  }$ is the canonical projection
$\pi_{T^{\ast }\mathcal{O}\left(  M\right)  }:T^{\ast}\mathcal{O}\left( 
M\right)
\rightarrow\mathcal{O}\left(  M\right)  $. By~(\ref{equation with Poisson
brackets}), the Hamiltonian semimartingale 
$\Gamma^{h}$ associated to $h=\left(  h_{1},...,h_{n}\right)$ and
to the stochastic Hamiltonian equations on $T ^\ast \mathcal{O}(M) $ with
stochastic component
$X=\left( B_{t}^{1},...,B_{t}^{n}\right)  $ is such that
\begin{multline*}
F\circ\pi_{T^{\ast}\mathcal{O}\left(  M\right)  }\left(  \Gamma^{h}\right) 
-F\circ\pi_{T^{\ast}\mathcal{O}\left(  M\right)  }\left( 
\Gamma_{0}^{h}\right)  \\
=\sum_{i=1}^{n}\int\left\{
F\circ\pi_{T^{\ast}\mathcal{O}\left(  M\right)  },h_{i}\right\} 
\left( 
\Gamma_{s}^{h}\right)  \delta B_{s}^{i}=\sum_{i=1}^{n}\int Y_{i}\left[ 
F\right]  \left(\pi_{T^{\ast}\mathcal{O}\left(  M\right)  }\left( 
\Gamma_{s}^{h}\right)\right)
\delta B_{s}^{i}
\end{multline*}
for any $F\in C^{\infty}\left(  \mathcal{O}\left(  M\right)  \right)  $. This
expression obviously implies that $U^{h}=\pi_{T^{\ast}\mathcal{O}\left(  M\right)  }\left(  \Gamma^{h}\right)
$ is a
solution of~(\ref{Brownian equation orthonormal}) and consequently
$X^{h}=\pi\left( U^{h}\right)  $ is a Brownian motion on $M$.

\subsection{The inverted pendulum with stochastically vibrating suspension point}

The equation of motion for small angles of a damped inverted
unit mass pendulum of length $l$ with a vertically vibrating suspension point is 
\begin{equation}
\label{inverted pendulum}
\ddot{\phi}=\left(  \frac{\ddot{y}}{l}+\frac{g}{l}\right)  \phi-\lambda
\dot{\phi}, 
\end{equation}
where $\phi $ is the angle that measures the separation of the pendulum from the vertical upright position, $y=y\left(  t\right)  $ is the height of the suspension point (externally controlled), $\lambda$
is the friction coefficient, and $g$ is the gravity constant. By construction, the point
$ (  \phi,\dot{\phi})  = (  0,0 )  $ corresponds to the
upright equilibrium position. It can be shown that if the function $ y(t)$ is of the form $y (t)= az\left(  \omega t\right) $, with $z $ periodic,   the   amplitude $a$ is sufficiently small, and  the frequency $\omega $ is sufficiently high, then this equilibrium becomes nonlinearly stable.

We now consider the case in which the external forcing of the suspension point is given by a continuous stochastic process $\dot{z}:\mathbb{R}_{+}\times\Omega\rightarrow\mathbb{R}$ such that 
$\dot{z}^{2}$ is  continuous and stationary. Under this assumptions, the equation~(\ref{inverted pendulum}) becomes the stochastic differential equation
\begin{equation}
\label{inverted pendulum stochastic}
d\phi  =\dot{\phi}dt, \quad
d\dot{\phi}   =\left(  \frac{g}{l}\phi-\lambda\dot{\phi}\right)
dt+\varepsilon^{2}\omega^{2}\phi d\dot{z}_{t}, 
\end{equation}
where $\varepsilon:=\sqrt{a/l}$. Observe that this equation is not 
Hamiltonian unless the friction term $-\lambda\dot{\phi}$ vanishes ($\lambda=0 $), in which case one obtains a Hamiltonian stochastic system with Hamiltonian function $h(\phi,\dot \phi)=(\frac{1}{2}(l ^2\dot \phi^2- l  \phi^2),\frac{1}{4}(\varepsilon^2 \omega ^2 \phi l )^2, - \frac{1}{2}( \varepsilon  \omega  \phi  l) ^2)$ and noise semimartingale $X _t=(t, [\dot z,\dot z], \dot z)$ (the symplectic form is obviously $ l ^2d \phi\wedge d \dot \phi $).

The stability of the upright position of the stochastically forced pendulum has been studied in~\cite{inverted pendulum, imkeller}, and references therein.  In~\cite{inverted pendulum} it is assumed that   the noise has the fairly strong mixing property. We recall that a continuous, adapted, stationary process  $\Gamma:\mathbb{R}_{+}\times\Omega\rightarrow\mathbb{R}$  has the fairly
strong mixing property if $E\left[  \Gamma_{t}^{2}\right]  <\infty$, there
exists a real function $c$ such that $\int_{0}^{\infty}c\left(  s\right)
ds<\infty$, and for any $t>s$%
\[
\left\Vert E\left[  \left.  \Gamma_{t}-E\left[  \Gamma_{t}\right]  \right\vert
\mathcal{F}_{s}\right]  \right\Vert _{L^{2}}\leq c\left(  t-s\right)
\left\Vert \Gamma_{s}-E\left[  \Gamma_{s}\right]  \right\Vert _{L^{2}},
\]
where $\left\Vert \cdot\right\Vert _{L^{2}}$ stands for the $L^{2}$ norm.
For example, if $x$ is the unique stationary solution with zero mean of the
It\^{o} equations
\begin{equation*}
dx_{t} =y_{t}dt, \quad
dy_{t} =-\left(  x_{t}+y_{t}\right)  dt+dB_{t},
\end{equation*}
where $B_{t}$ is a standard Brownian motion, then $\dot{x}_{t}^{2}-\frac{1}%
{2}=y_{t}^{2}-\frac{1}{2}$ has the fairly strong mixing property. Using this hypothesis, it can be shown~\cite[Theorem 1]{inverted pendulum} that if
$z:\mathbb{R}_{+}\times
\Omega\rightarrow\mathbb{R}$ is a continuously differentiable and stationary
process such that, for any $t\in\mathbb{R}_{+}$, $E\left[  z_{t}\right]  =0$,
$E\left[  \exp\left(  \varepsilon\left\vert z_{t}\right\vert \right)  \right]
<\infty$ if $\varepsilon=\sqrt{a/l}$ is sufficiently small, and the process
$\dot{z}^{2}$ has the fairly strong mixing property, then the solution
$ (  \phi,\dot{\phi} )  = (  0,0 )  $ of~(\ref{inverted pendulum stochastic}) is
exponentially stable in probability, if $\varepsilon$ is sufficiently small and
$\frac{g}{l\varepsilon^{4}}<E\left[  \dot{z}^{2}\right]  $.
Moreover, Ovseyevich shows in \cite[Section 4]{inverted pendulum} that if we
put $\lambda=0$ in (\ref{inverted pendulum stochastic}) and we consider hence the inverted pendulum
as a Hamiltonian system, then the equilibrium point $ (  \phi,\dot{\phi
} )  = (  0,0 )  $ is unstable.

\section{Critical action principles for the stochastic Hamilton equations}
\label{A critical action principle for the stochastic Hamilton equations}

Our goal in this section is showing that
the stochastic Hamilton equations can be characterized by  a variational principle that
generalizes the one used in the classical deterministic situation. In the
following pages we shall consider  an exact
symplectic manifold
$(M,
\omega)$, that is, there exist a one-form $\theta\in
\Omega\left(  M\right)  $ such that $\omega=-\mathbf{d}\theta$. The archetypical example of
an exact symplectic manifold is the cotangent bundle $T^{\ast}Q$ of any manifold $Q$, with
$\theta $ the Liouville one-form.   

In the following pages we will proceed in two stages. In the first subsection we will construct a critical action principle based on using variations of the solution semimartingale using the flow of a vector field on the manifold. Even though this approach is extremely natural and mathematically very tractable it yields a variational principle (Theorem~\ref{Critical Action Principle theorem}) that does not fully characterize the stochastic Hamilton's equations. In order to obtain such a characterization one needs to use more general variations associated to the flows of vector fields defined on the solution semimartingale, that is, they depend on $\Omega $. This complicates considerably the formulation and will be treated separately in the second subsection.

\begin{definition}
\label{action definition}
Let $(M, \omega=- \mathbf{d}\theta)$ be an exact symplectic manifold, $X: \mathbb{R} _+
\times  \Omega \rightarrow V $ a semimartingale taking values on the vector space $V $, and
$h :M \rightarrow  V ^\ast   $  a Hamiltonian function. We denote by $\mathcal{S}\left( 
M\right) 
$ and
$\mathcal{S}\left(  \mathbb{R}\right)  $  the sets of $M$ and real-valued
semimartingales, respectively.
We define the {\bfi  stochastic action} associated to $h$ as the
map $S: \mathcal{S}(M) \rightarrow \mathcal{S}(\mathbb{R}) $ given by
\[
S\left(  \Gamma\right)  =\int\left\langle \theta,\delta\Gamma\right\rangle
-\int\left\langle \widehat{h}\left(  \Gamma\right)  ,\delta X\right\rangle ,
\]
where in the previous expression, $\widehat{h}\left(  \Gamma\right): \mathbb{R}_+ \times 
\Omega \rightarrow V  \times  V^\ast $ is given by $\widehat{h}\left( 
\Gamma\right)(t,\omega):=(X _t (\omega), h(\Gamma_t(\omega)))$.
\end{definition}

\subsection{Variations involving vector fields on the phase space}

\begin{definition}
\label{def differentiable action}
Let $M$ be a manifold,
$F:\mathcal{S}\left(  M\right)
\rightarrow\mathcal{S}\left(  \mathbb{R}\right)  $ a map, and $\Gamma\in\mathcal{S}\left( 
M\right)  $. 
A local one-parameter group of diffeomorphisms $\varphi: \mathcal{D} \subset \mathbb{R}\times M \rightarrow  M $ is said to be {\bfi  complete with respect} to $\Gamma $ if there exists $\epsilon> 0$ such that $\varphi _s(\Gamma  )$ is a well-defined process for any $s \in (- \epsilon, \epsilon)$.
We  say that $F$ is {\bfi  differentiable} at $\Gamma $ in the direction of a
local one parameter group of diffeomorphisms
$\varphi$ complete with respect to $\Gamma$, if for any
sequence
$\left\{  s_{n}\right\}  _{n\in\mathbb{N}}\subset \mathbb{R}$ such that $s_{n}\underset
{n\rightarrow\infty}{\longrightarrow}0$, the family
\[
X_{n}=\frac{1}{s_{n}}\left(  F\left(  \varphi_{s_{n}}\left(  \Gamma\right)
\right)  -F\left(  \Gamma\right)  \right)
\]
converges uniformly on compacts in probability  (ucp) to a process  that we will
denote by $\left.  \frac{d}{ds}\right\vert _{s=0}F\left( 
\varphi_{s}\left( 
\Gamma\right) 
\right)  $ and that is referred to as the {\bfi  directional derivative} of $F$ at $\Gamma $
in the direction of $\varphi _s$.
\end{definition}

\begin{remark}
\label{remark with stopping time}
\normalfont
Note that global one-parameter groups of diffeomorphisms (for instance, flows of complete vector fields) are complete with respect to any semimartingale. 

Let
$\Gamma:
\mathbb{R}_+
\times
\Omega
\rightarrow M
$ be a
$M$-valued continuous and adapted stochastic process and $A\subset M$ a  set. We will denote by
$\tau_{A}=\inf\left\{ 
 t>0\mid\Gamma_{t}\left(  \omega\right)  \notin A\right\}  $ the {\bfi 
first exit time} of $\Gamma$ with respect to $A$. We recall that $\tau _A$
is a stopping time if $A$ is a Borel set.
Additionally, let  $\Gamma$ be a semimartingale and  $K$ a compact set such that $ \Gamma _0\subset K $. Then, any local one-parameter group of diffeomorphisms $\varphi $ is complete with respect to the stopped process $\Gamma^{\tau _K}$. Note that this conclusion could also hold for certain non-compact sets.
\end{remark}

\noindent The proof of the following proposition can be found in Section~\ref{proof of proposition crucial}. 

\begin{proposition}
\label{prop differentiation processes}
Let $M$ be a manifold, $\alpha\in\Omega\left(  M\right)  $
a one-form, and $F:\mathcal{S}\left(  M\right)
\rightarrow\mathcal{S}\left(  \mathbb{R}\right)  $ the map defined by 
$F\left(\Gamma\right)  :=\int\left\langle \alpha
,\delta  \Gamma   \right\rangle $. Then $F$ is differentiable in all directions. Moreover,
if $\Gamma:\mathbb{R}_{+}\times\Omega\rightarrow M$ is a continuous
semimartingale, $\varphi$ is an arbitrary local one-parameter
group of diffeomorphisms complete with respect to $\Gamma $, and $Y \in \mathfrak{X} (M) $ is the vector field associated to
$\varphi$, then 
\begin{equation}
\label{derivative of functional}
\left.  \frac{d}{ds}\right\vert _{s=0}F\left( 
\varphi_{s}\left( 
\Gamma\right) 
\right)=\left.  \frac{d}{ds}\right\vert _{s=0}\int\left\langle
\alpha,\delta\left(\varphi_{s} \circ \Gamma \right)\right\rangle =\left. 
\frac{d}{ds}\right\vert _{s=0}\int\left\langle \varphi_{s}^{\ast
}\alpha,\delta\Gamma\right\rangle =\int\left\langle \pounds_{Y}\alpha,\delta
\Gamma\right\rangle.
\end{equation}
The symbol $\pounds_{Y}\alpha $ denotes the Lie derivative of $\alpha$ in the direction
given by $Y$.
\end{proposition}

\begin{corollary}
\label{corollary differentiability of S}
In the setup of Definition~\ref{action definition} let
$\alpha=\omega^{\flat}\left(  Y\right)\in \Omega (M)
$, with $\omega^\flat $ the inverse of the vector bundle
isomorphism $\omega^{\sharp}: T ^\ast M
\rightarrow  TM
$  induced by $\omega$.   
Let
$\Gamma:\mathbb{R}_{+}
\times\Omega\rightarrow M$ be a continuous adapted semimartingale.
  $\varphi$ an arbitrary local
one-parameter group of diffeomorphisms complete with respect to $\Gamma $,  and $Y \in  \mathfrak{X} (M) $ the associated vector
field. Then, the action $S$ is
differentiable at $\Gamma $ 
in the direction of $\varphi$ and the directional derivative is given by
\begin{equation}
\label{derivative formula for s}
\left.  \frac{d}{ds}\right\vert _{s=0}S\left(  \varphi_{s}\left(
\Gamma\right)  \right)  =-\int\left\langle \alpha,\delta\Gamma\right\rangle
-\int\left\langle  \mathbf{d}h\left( \omega^{\#}\left(  \alpha\right)
\right)  \left(  \Gamma\right)  ,\delta X\right\rangle +{\bf i}_{Y}\theta\left(
\Gamma\right)  -{\bf i}_{Y}\theta\left(  \Gamma_{0}\right).
\end{equation}
\end{corollary}

\noindent\textbf{Proof.\ \ }
It is clear from Proposition~\ref{prop differentiation processes} that
\[
\frac{1}{s}\left[  \int\left\langle \varphi_{s}^{\ast}\theta-\theta
,\delta\Gamma\right\rangle \right]  \overset{s\rightarrow 0}{\longrightarrow
}\int\left\langle \pounds_{Y}\theta,\delta\Gamma\right\rangle
\]
in $ucp$. The proof of that result can be easily adapted to show that ucp
\[
\frac{1}{s}\left[  \int\left\langle \left(  \varphi_{s}^{\ast}\widehat{h}-\widehat{h}\right)
\left(  \Gamma\right)  ,\delta X\right\rangle \right]  \overset
{s\rightarrow0}{\longrightarrow}\int\left\langle \left(  \pounds_{Y}\widehat{h}\right)
\left(  \Gamma\right)  ,\delta X\right\rangle .
\]
Thus, using~(\ref{properties strato}) and $\alpha=\omega^{\flat}\left(  Y\right)\in \Omega
(M)$,
\begin{align*}
\left.  \frac{d}{ds}\right\vert _{s=0}S\left(  \varphi_{s}\left(
\Gamma\right)  \right)   & =\int\left\langle \pounds_{Y}\theta,\delta\Gamma
\right\rangle -\int\left\langle \left(  \pounds_{Y}\widehat{h}\right)  \left(  \Gamma\right)
,\delta X\right\rangle =\int\left\langle {\bf i}_{Y}\mathbf{d}\theta+
\mathbf{d}\left( 
{\bf i}_{Y}\theta\right)  ,\delta
\Gamma\right\rangle -\int\left\langle   \mathbf{d}h\left(   Y 
\right)  \left(  \Gamma\right)  ,\delta X\right\rangle \\
  &
=-\int\left\langle \alpha,\delta\Gamma\right\rangle +\int\left\langle \mathbf{d}\left( 
{\bf i}_{Y}\theta\right) ,\delta\Gamma\right\rangle -\int\left\langle
\mathbf{d}h\left( \omega^{\#}\left(  \alpha\right)  \right)  \left(
\Gamma\right)  ,\delta X\right\rangle \\
& =-\int\left\langle \alpha,\delta
\Gamma\right\rangle -\int\left\langle  \mathbf{d}h\left( \omega
^{\#}\left(  \alpha\right)  \right)  \left(  \Gamma\right)  ,\delta
X\right\rangle  +\left( {\bf i}_{Y}\theta\right)  \left(  \Gamma\right)  -\left( {\bf
i}_{Y}\theta\right)  \left(  \Gamma_{0}\right).\quad \blacksquare
\end{align*}

\begin{corollary}[Noether's theorem]
In the setup of Definition~\ref{action definition}, let $\varphi:\mathbb{R}
\times M\rightarrow M$ be a one parameter group of diffeomorphisms and
$Y \in  \mathfrak{X} (M)$ the associated vector field. If the action
$S:\mathcal{S}\left(  M\right)  \rightarrow\mathcal{S}\left(  \mathbb{R}\right)  $
is invariant by $\varphi$, that is,  $S\left(
\varphi_{s}\left(  \Gamma\right)  \right)  =S\left(  \Gamma\right) $, for any $s \in \mathbb{R}$, then
the function ${\bf i}_{Y}\theta$ is a  conserved quantity of the stochastic
Hamiltonian system associated to $h:M\rightarrow V^{\ast}$.
\end{corollary}

\noindent\textbf{Proof.\ \ } Let $ \Gamma^h $ be the Hamiltonian semimartingale associated
to $h $  with initial condition $\Gamma_0$. Since $\varphi_s  $ leaves invariant the action
we have that 
\[
\left.  \frac{d}{ds}\right\vert _{s=0}S\left(  \varphi_{s}\left(
\Gamma^h\right)  \right)  =0
\]
and hence by~(\ref{derivative formula for s}) we have that
\[
0 =-\int\left\langle \alpha,\delta\Gamma^h\right\rangle
-\int\left\langle \mathbf{d}h\left(  \omega^{\#}\left(  \alpha\right)
\right)  \left(  \Gamma^h\right)  ,\delta X\right\rangle +{\bf i}_{Y}\theta\left(
\Gamma^h\right)  -{\bf i}_{Y}\theta\left(  \Gamma_{0}\right).
\]
As $ \Gamma^h $ is the Hamiltonian semimartingale associated
to $h $ we have that 
\[-\int\left\langle \alpha,\delta\Gamma^h\right\rangle
=\int\left\langle  \mathbf{d}h\left( \omega^{\#}\left(  \alpha\right)
\right)  \left(  \Gamma^h\right)  ,\delta X\right\rangle\] 
and hence ${\bf i}_{Y}\theta\left(
\Gamma^h\right)  ={\bf i}_{Y}\theta\left(  \Gamma_{0}\right)$, as required. \quad
$\blacksquare$

\begin{remark}
\normalfont
The hypotheses of the previous corollary can be modified by requiring,
instead of the invariance of the action by $\varphi _s $, the existence of a
function $F \in C^\infty(M) $ such that 
\[
\left.  \frac{d}{ds}\right\vert _{s=0}S\left(  \varphi_{s}\left(
\Gamma^h\right)  \right)  =F(\Gamma)-F(\Gamma_0).
\]
In that situation, the conserved quantity is  ${\bf i}_{Y}\theta+F$.
\end{remark}

Before we state the Critical Action Principle for the stochastic Hamilton
equations  we need one more definition.

\begin{definition}
Let $M$ be a manifold and $A$ a set. We will say that a local one parameter group of diffeomorphisms $\varphi
:\mathcal{D} \times M\rightarrow M$  {\bfi 
fixes} $A$ if $\varphi_{s}\left(  y\right)  =y$ for any
$y\in A$ and any $s \in \mathbb{R}$ such that  $(s,y)\in \mathcal{D}$. The
corresponding vector field $Y \in \mathfrak{X} (M) $ given by $Y (m)=
\left.\frac{d}{ds}\right|_{s=0} \varphi _s (m) $ satisfies that $Y|_{A}=0$. 
\end{definition}

\begin{theorem}[First Critical Action Principle]
\label{Critical Action Principle theorem}Let $\left(  M,\omega=-d\theta
\right)  $ be an exact symplectic manifold, $X:\mathbb{R}_{+}\times
\Omega\rightarrow V$ a semimartingale taking values on the vector space $V$
such that $X _0=0 $, and $h:M\rightarrow V^{\ast}$ a Hamiltonian function.
Let
$m_{0}\in M$ be a point in $M$ and
$\Gamma:\mathbb{R}_{+}\times\Omega\rightarrow M$ a continuous
semimartingale such that $\Gamma_{0}=m_{0}$. Let $K $ be a compact set that contains the point $m_{0}$. If the
semimartingale $\Gamma$ satisfies the stochastic Hamilton
equations~(\ref{solution Hamiltonian integral})  (with initial condition
$\Gamma_{0}=m_{0}$) up to time $\tau _K $ then  for any
local one-parameter group of diffeomorphisms
$\varphi$ that fixes the set
$\{ m_{0}\}\cup\partial K$ we have
\begin{equation}
{\bf 1}_{\{\tau_K< \infty\}}\left[  \left.  \frac{d}{ds}\right\vert
_{s=0}S\left( 
\varphi_{s}\left(
\Gamma^{\tau_K}\right)  \right)  \right]  _{\tau_{K}}=0\text{ \ a.s.}.\label{critical action 6}%
\end{equation}

\end{theorem}

\noindent\textbf{Proof.\ \ } We start by emphasizing that when we write that $\Gamma$ satisfies the
stochastic Hamiltonian equations~(\ref{solution Hamiltonian integral}) up to time $\tau _K $ we mean that
\[
\left(  \int\left\langle \beta,\delta\Gamma\right\rangle +\int\left\langle
  \mathbf{d}h\left(\omega^{\#}\left(  \beta\right)  \right)  \left(
\Gamma\right)  ,\delta X\right\rangle \right)  ^{\tau_{K}}=0.
\]
For the sake of simplicity in our notation we define the linear operator
$\mbox{\rm Ham}: \Omega(M) \rightarrow \mathcal{S}(\mathbb{R})$ given by
\[
\operatorname*{Ham}\left(  \beta\right)  :=\left(  \int\left\langle
\beta,\delta\Gamma\right\rangle +\int\left\langle   \mathbf{d}%
h\left(\omega^{\#}\left(  \beta\right)  \right)  \left(  \Gamma\right) 
,\delta X\right\rangle \right), \qquad \beta\in\Omega\left(  M\right).  
\]
Suppose now that the semimartingale $\Gamma$ satisfies the
stochastic Hamilton equations up to time $\tau _K $.
Let 
$\varphi$ be
a local one-parameter group of diffeomorphisms that fixes $\{
m_{0}\}\cup\partial K$,  and let $Y \in \mathfrak{X}(M) $ be the associated
vector field. Then, taking
$\alpha=\omega^{\flat}\left(  Y\right) $,  we have by
Corollary~\ref{corollary differentiability of S}, 
\begin{equation}
\label{before imposing condition}
\left.  \frac{d}{ds}\right\vert _{s=0}S\left(  \varphi_{s}\left(
\Gamma^{\tau_{K}}\right)  \right)  =-\int\left\langle \alpha,\delta\Gamma^{\tau_{K}}\right\rangle
-\int\left\langle   \mathbf{d}h\left(\omega^{\#}\left(  \alpha\right)
\right)  \left(  \Gamma^{\tau_{K}}\right)  ,\delta X\right\rangle +i_{Y}\theta\left(
\Gamma^{\tau_{K}}\right),
\end{equation}
since $Y\left(  m_{0}\right)  =0$ and hence ${\bf i}_{Y}\theta\left(  \Gamma_{0}\right)=0 $.
Additionally, since $\Gamma $ is continuous, ${\bf 1}_{\{\tau_K<
\infty\}}\Gamma_{\tau_K}
\in
\partial K $ and $Y|_{\partial K}=0
$. Hence,
\begin{equation*}
{\bf 1}_{\{\tau_K< \infty\}}\left[\left.  \frac{d}{ds}\right\vert _{s=0}S\left(  \varphi_{s}\left(
\Gamma^{\tau_{K}}\right)  \right)\right]_{\tau _K}  =-{\bf 1}_{\{\tau_K< \infty\}}\left[\int\left\langle \alpha,\delta\Gamma^{\tau_{K}}\right\rangle
+\int\left\langle   \mathbf{d}h\left(\omega^{\#}\left(  \alpha\right)
\right)  \left(  \Gamma^{\tau_{K}}\right)  ,\delta X\right\rangle \right]_{\tau_K}.
\end{equation*}
Now,  Proposition~\ref{lemma restringir intervalo temporal} and the hypothesis on $\Gamma$ satisfying Hamilton's equation guarantee that the previous expression equals
\begin{eqnarray*}
{\bf 1}_{\{\tau_K< \infty\}}\left[  \left.  \frac{d}{ds}\right\vert
_{s=0}S\left(  \varphi_{s}\left(
\Gamma^{\tau_K}\right)  \right)  \right]  _{\tau_{K}}&=&
-{\bf 1}_{\{\tau_K< \infty\}}\left[\left[\int\left\langle \alpha,\delta\Gamma^{\tau_{K}}\right\rangle
+\int\left\langle   \mathbf{d}h\left(\omega^{\#}\left(  \alpha\right)
\right)  \left(  \Gamma^{\tau_{K}}\right)  ,\delta X\right\rangle\right]^{\tau_K} \right]_{\tau_K}\\
	&= &-{\bf 1}_{\{\tau_K< \infty\}}\left[\left[\int\left\langle \alpha,\delta\Gamma\right\rangle
+\int\left\langle   \mathbf{d}h\left(\omega^{\#}\left(  \alpha\right)
\right)  \left(  \Gamma\right)  ,\delta X\right\rangle\right]^{\tau_K} \right]_{\tau_K}\\
	&=&-{\bf 1}_{\{\tau_K<
\infty\}}\left[
\mbox{\rm Ham} (\alpha)_{\tau_{K}}\right]=0\text{ \ a.s.},
\end{eqnarray*}
as required. \quad $\blacksquare$

\begin{remark}
\normalfont
The relation between the Critical Action Principle stated in
Theorem~\ref{Critical Action Principle theorem} and the classical one for
Hamiltonian mechanics is not straightforward since the categories in which
both are formulated are very much different; more specifically, the
differentiability hypothesis imposed on the solutions of the
deterministic principle is not a reasonable assumption in the
stochastic context and this has serious consequences. For example, unlike
the situation encountered in classical mechanics, Theorem~\ref{Critical
Action Principle theorem} does not admit a converse within the set of
hypotheses in which it is formulated.

In order to elaborate a little bit more on this question let $\left(  M,\omega=-d\theta
\right)  $ be an exact symplectic manifold, take the
Hamiltonian function $h \in C^\infty(M)  $, and consider the stochastic
Hamilton equations with trivial stochastic component $X: \mathbb{R}_+
\times  \Omega \rightarrow  \mathbb{R}$ given by $X _t(\omega)=t $. As we
saw in Remark~\ref{relation classical equations 1} the paths
of the semimartingales that solve these stochastic Hamilton equations are the
smooth curves that integrate the standard Hamilton equations.
In this situation the action reads
\[
S\left(  \Gamma\right)  =\int\langle\theta, \delta \Gamma\rangle-\int
h\left( 
\Gamma_s\right)  ds.
\]
If the path $\Gamma _t(\omega)$ is differentiable then the integral
$\left(\int\langle\theta, \delta \Gamma\rangle\right)(\omega) $
reduces to the Riemann integral $\int_{\Gamma_t(\omega)}\theta $ and
$S (\Gamma)(\omega)$ coincides with the classical action. In particular, if
$\Gamma $ is a solution of the stochastic Hamilton equations then the paths
$\Gamma_t(\omega)$ are necessarily differentiable (see Remark~\ref{relation
classical equations 1}), they satisfy the standard Hamilton equations, and
hence make the action critical. The following elementary example shows that
the converse is not necessarily true, that is one may have semimartingales
that satisfy~(\ref{critical action 6}) and that do not solve the Hamilton equations up to time $\tau _K $. 

We will consider a deterministic example.
Let
$m_{0},$
$m_{1}\in M$ be two points. Suppose there exists an integral curve
$\gamma:\left[  t_{0},t_{1}\right]  \rightarrow M$ of the Hamiltonian vector
field $X_{h}$ defined on some time interval $\left[  t_{0},t_{1}\right]  $
such that $\gamma\left(  t_{0}\right)  =m_{0}$ and $\gamma\left(
t_{1}\right)  =m_{1}$. Define the continuous and piecewise smooth curve $\sigma:\left[
0,t_{1}\right] 
\rightarrow M$ as follows:%
\[
\sigma\left(  t\right)  =\left\{
\begin{array}{ll}
m_{0}&\text{ \ if }t\in\left[  0,t_{0}\right]  \\
\gamma\left(  t\right)  &\text{ \ if }t\in\left[  t_{0},t_{1}\right]  .
\end{array}
\right.
\]
Let $\varphi$ be
a local one-parameter group of diffeomorphisms that fixes
$\{m_{0},m_{1}\}$. Then by~(\ref{derivative formula for s}) 
\[
\left[\left.  \frac{d}{ds}\right\vert _{s=0}S\left(  \varphi_{s}\left(
\sigma\right)  \right) \right]_t
=-\int_{\sigma|_{[0,t]}}\alpha+\int_{0}^{t}\langle \alpha, 
X_{h}\rangle\left( 
\sigma(t)\right) dt+\langle \theta(\sigma(t)), Y(\sigma(t))\rangle-\langle
\theta(m _0), Y(m _0)\rangle,
\]
where  $Y\left(  m\right)  =\left.
\frac{d}{ds}\right\vert _{s=0} \varphi_{s}\left(  m\right) $, for any $m\in
M$ and $\alpha=\omega^{\flat}\left(  Y\right)  $. Using that $\sigma$
satisfies the Hamilton equations on $\left[  t_{0} ,t_{1}\right]$ and
$\alpha\left(  m_{0}\right)  =0$, it is easy to see
that
\[
\left[  \left.  \frac{d}{ds}\right\vert _{s=0}S\left(  \varphi_{s}\left(
\sigma\right)  \right)  \right]  _{t_{1}}=0,
\]
that is, $\sigma $ makes the action critical. However, it
does not satisfy the Hamilton equations on the interval $\left[ 
0,t_{1}\right]  ,$ because they do not hold on
$\left(  0,t_{0}\right)  $. This shows that the converse of the statement in
Theorem~\ref{Critical Action Principle theorem} is not necessarily true. In the following subsection we will obtain such a converse by generalizing the set of variations allowed in the variational principle.
\end{remark}

\subsection{Variations involving vector fields on the solution semimartingale}

We start by spelling out the variations that we will use in order to obtain a converse to Theorem~\ref{Critical Action Principle theorem}.

\begin{definition}
\label{def pathwise variation}
Let $M$ be a manifold and $\Gamma$ a $M$-valued semimartingale. Let $s_{0}>0$; we say that the map $\Sigma:\left(  -s_{0},s_{0}\right)  \times
\mathbb{R}_{+}\times\Omega \rightarrow M $ is a {\bfi  pathwise variation} of  $\Gamma$ whenever $\Sigma_{t}^{0}=\Gamma_{t}$
for any $t\in\mathbb{R}_{+}$ a.s..

We say that the pathwise variation $\Sigma$ of $\Gamma $ {\bfi  converges uniformly} to $\Gamma$
whenever the
following properties are satisfied:
\begin{description}
\item [(i)] For any $f\in C^{\infty}\left(  M\right)  $, $f\left(  \Sigma
^{s}\right)  \rightarrow f\left(  \Gamma\right)  $ in $ucp$ as $s\rightarrow0$.
\item [(ii)] There exists a process $Y:\mathbb{R}_{+}\times\Omega\rightarrow TM$ over
$\Gamma$  such that, for any $f\in C^{\infty}\left(  M\right)  $, the
Stratonovich integral $\int Y\left[  f\right]  \delta X$ exists for any
continuous real semimartingale $X$ (this is for instance guaranteed if $Y$ is a semimartingale)
and, additionally, the increments $\left.  \left(  f\left(  \Sigma^{s}\right)
-f\left(  \Gamma\right)  \right)  \right/  s$ converge in $ucp$ to $Y\left[
f\right]  $ as $s\rightarrow0$. We will call such a $Y$ the
\textbf{infinitesimal generator} of $\Sigma$.
\end{description}
We will say that $\Sigma $ (respectively, $Y $) is {\bfi  bounded} when its image lies in a compact set of  $M$ (respectively, $TM $).
\end{definition}

The next proposition shows that, roughly speaking, there exist bounded pathwise variations that converge uniformly to a given semimartingale with prescribed bounded infinitesimal generator.

\begin{proposition}
\label{existence pathwise variation}
Let $\Gamma$ be a  continuous $M$-valued
semimartingale $\Gamma$, $K\subseteq M$ a compact set, and $\tau_{K}$ the
first exit time of $\Gamma$ from $K$. Let $Y:\mathbb{R}_{+}\times
\Omega\rightarrow TM$ be a bounded process over $\Gamma^{\tau_{K}}$ (that is, the image of $Y$ lies in a compact subset of $TM $) such that $\int
Y\left[  f\right]  \delta X$ exists for any continuous real semimartingale $X$
and for any $f\in C^{\infty}\left(  M\right)  $. Then, there exists a bounded pathwise
variation $\Sigma$ that converges  uniformly to $\Gamma^{\tau_{K}}$ whose
infinitesimal generator is $Y$.
\end{proposition}

\noindent\textbf{Proof.\ \ }
Let $\left\{  (V_{k}, \varphi_k)\right\}  _{k\in\mathbb{N}}$ be a countable open covering
of $M$ by coordinate patches such that any $V_{k}$ is contained in a compact
set. This covering is always available by the second countability of the manifold and Lindel\"of's Lemma.
Let $\left\{  U_{k}\right\}  _{k\in\mathbb{N}}$ be an open subcovering
such that, if $U_{k}\subseteq V_{i}$ for some $k$, $i\in\mathbb{N}$, then
$U_{k}\subsetneq V_{i}$. Let $\left\{  \tau_{m}\right\}  _{m\in\mathbb{N}}$ be
a sequence of stopping times (available by Lemma 3.5 in~\cite{emery}) such that, a.s., $\tau_{0}=0,$ $\tau_{m}
\leq\tau_{m+1},$ $\sup_{m}\tau_{m}=\infty,$ and that, on each of the sets
$\left[  \tau_{m},\tau_{m+1}\right]  \cap\left\{  \tau_{m+1}>\tau_{m}\right\}
$ the semimartingale $\Gamma$ takes  values in the open set $U_{k\left(  m\right)
}$, for some $k\left(  m\right)  \in\mathbb{N}$. Since $K$ is compact, it can
be covered by a finite number of these open sets, i.e. $K\subseteq\cup_{j\in
J}U_{k_{j}}$, where $|J|<\infty$.

Let $x_{k_j}\equiv (x_{k_{j}}^{1}, \ldots, x_{k_{j}}^{n})$, $n=\dim\left(  M\right) $ be a set of coordinate
functions on $U_{k_{j}\left(  m\right)  }$ and $(x_{k_j},v_{k_j})\equiv (x_{k_{j}}^{1}, \ldots, x_{k_{j}}^{n}
,v_{k_{j}}^{1}, \ldots, v_{k_{j}}^{n})  $ the corresponding adapted coordinates for $TM$ on $\pi_{TM}^{-1}\left(
U_{k_{j}\left(  m\right)  }\right)  $. Since  $Y$ is bounded and covers $\Gamma^{\tau_K} $, and on $\left[  \tau_{m},\tau_{m+1}\right]  \cap\left\{  \tau_{m+1}>\tau_{m}\right\}
$ the semimartingale $\Gamma$ takes  values in the open set $U_{k_j\left(  m\right)
}$,  there exist a $s_{k_{j}}>0$
such that, on $\left[  \tau_{m},\tau_{m+1}\right]  \cap\left\{  \tau
_{m+1}>\tau_{m}\right\}  $, the points
$ (x_{k_{j}}^{1}\left(  \Gamma\right)  +sv_{k_{j}}^{1}\left(  Y\right), \ldots,x_{k_{j}}^{n}\left(  \Gamma\right)  +sv_{k_{j}}^{n}\left(  Y\right))
$
lie in the image of some coordinate patch $V_{k_{j}}$ containing
$U_{k_{j}\left(  m\right)  }$ for all $s\in\left(  -s_{k_{j}},s_{k_{j}
}\right)  $. Let $s_{0}=\min_{j\in J}\left\{  s_{k_{j}}\right\}  $. Now, since the sets of the form
$I _m:=\left[  \tau_{m},\tau_{m+1}\right)  \cap\left\{  \tau_{m+1}>\tau_{m}\right\}\subset \mathbb{R}_+ \times \Omega
$, $m \in  \mathbb{N} $  form a disjoint partition of $\mathbb{R}_+ \times  \Omega $ we define $\Sigma$ as the map that for any $m \in  \mathbb{N} $ satisfies
\[
\begin{array}
[c]{rrl}
\Sigma|_{I_m}:\left(  -s_{0},s_{0}\right)  \times\left[  \tau_{m},\tau_{m+1}\right)
\cap\left\{  \tau_{m+1}>\tau_{m}\right\}  & \longrightarrow & V_{k_{j}}\\
\left(  s,t,\omega\right)  & \longmapsto & \varphi_k ^{-1}\left(  x_{k_{j}} \left(
\Gamma_{t}\left(  \omega\right)  \right)  +sv_{k_{j}} \left(  Y_{t}\left(
\omega\right)  \right)  \right).
\end{array}
\]
Observe that by construction  the image of $\Sigma$ is covered by a finite number of
coordinated patches and therefore, by hypothesis, contained in a compact set. $\Sigma$ is hence bounded. More specifically 
\begin{equation}
\label{obs covering finite patches}
\left\{  \Sigma_{t}^{s}\left(  \omega\right)  \mid (s,t, \omega)\in \left(  -s_{0}
,s_{0}\right)  \times\mathbb{R}\times\Omega\right\}  \subseteq\bigcup_{j\in
J}V_{k_{j}}. 
\end{equation}
It is immediate to see that $\Sigma$ is a pathwise variation which
converges uniformly to $\Gamma^{\tau_{K}}$. Indeed, if $f\in C^{\infty}\left(
M\right)  $ has compact support within one of the elements in the family $\left\{  U_{k_{j}}\right\}
_{j\in J}$, it can be easily checked that
\begin{equation}
f\left(  \Sigma^{s}\right)  \underset{s\rightarrow0}{\underset{ucp}%
{\longrightarrow}}f\left(  \Gamma\right)  \text{ \ \ \ and\ \ \ \ }%
\frac{f\left(  \Sigma^{s}\right)  -f\left(  \Gamma\right)  }{s}\underset
{s\rightarrow0}{\underset{ucp}{\longrightarrow}}Y\left[  f\right]  .
\label{eq 5 original}
\end{equation}
If, more generally, $f\in C^{\infty}\left(  M\right)  $ has not compact support contained in
one of the $\left\{  U_{k_{j}}\right\}  _{j\in J}$, observe that, by
(\ref{obs covering finite patches}), we only need to consider the restriction of $f$ to
$\bigcup_{j\in J}V_{k_{j}}$. Take now a partition of the unity $\left\{
\phi_{k}\right\}  _{k\in\mathbb{N}}$ subordinated to the covering $\left\{
U_{k}\right\}  _{k\in\mathbb{N}}$. Since $\left\{  \operatorname*{supp}\left(
\phi_{k}\right)  \right\}  _{k\in\mathbb{N}}$ is a locally finite family and
$\bigcup_{j\in J}V_{k_{j}}$ is contained in a compact set because, by hypothesis, so is each
$V_{k_{j}}$ for any $j\in J$, then among all the $\left\{  \phi
_{k}\right\}  _{k\in\mathbb{N}}$ only a finite number of them have their
supports  in $\left\{  U_{k_{j}}\right\}  _{j\in J}$, 
say $\left\{  \phi_{k_{i}}\right\}  _{i\in I}$ with $\left\vert I\right\vert
<\infty$. Thus,
\[
f| _{\cup_{j\in J}V_{k_{j}}}=\sum_{i=1}^{\left\vert
I\right\vert }\phi_{k_{i}}f
\]
and since each $\phi_{k_{i}}f$ is a function similar to those  considered in~(\ref{eq 5 original}) it is
straightforward to see that those implications also hold for $f$.
\quad $\blacksquare$

\medskip

The following result generalizes Proposition~\ref{prop differentiation processes} to pathwise variations of a semimartingale. The proof can be found in Section~\ref{appendix proposition derivada pathwise}

\begin{proposition}
\label{derivada pathwise variation}
Let $\Gamma$ be a $M$-valued continuous
semimartingale $\Gamma$, $K\subseteq M$ a compact set, and $\tau_{K}$ the
first exit time of $\Gamma$ from $K$. Let $\Sigma$ be a bounded pathwise variation
that converges uniformly to $\Gamma^{\tau_{K}}$  and $Y:\mathbb{R}_{+}\times
\Omega\rightarrow TM$ the infinitesimal generator of $\Sigma$ that we will also assume to be bounded. Then, for any
$\alpha\in\Omega\left(  M\right)  $,%
\[
\underset{\underset{s\rightarrow0}{ucp}}{\lim}~\frac{1}{s}\left[
\int\left\langle \alpha,\delta\Sigma^{s}\right\rangle -\int\left\langle
\alpha,\delta\Gamma^{\tau_{K}}\right\rangle \right]  =\int\left\langle
i_{Y}\mathbf{d}\alpha,\delta\Gamma^{\tau_{K}}\right\rangle +\left\langle
\alpha\left(  \Gamma^{\tau_{K}}\right)  ,Y\right\rangle -\left\langle
\alpha\left(  \Gamma^{\tau_{K}}\right)  ,Y\right\rangle _{t=0}.
\]
\end{proposition}

The next theorem shows that the generalization of the Critical Action Principle in Theorem~\ref{Critical Action Principle theorem} to pathwise variations fully characterizes the stochastic Hamilton's equations.

\begin{theorem}[Second Critical Action Principle]
\label{Second Critical Action Principle}
Let $\left(  M,\omega=-d\theta\right)  $ be an exact symplectic manifold,
$X:\mathbb{R}_{+}\times\Omega\rightarrow V$ a semimartingale that takes values in
the vector space $V,$ and $h:M\rightarrow V^{\ast}$ a Hamiltonian function.
Let $m_{0}$ be a point in $M$ and $\Gamma:\mathbb{R}_{+}\times\Omega
\rightarrow M$ a continuous adapted semimartingale defined on $\left[  0,\zeta_{\Gamma}\right)  $ such that $\Gamma
_{0}=m_{0}$. Let $K\subseteq M$
be a compact set that contains $m _0 $ and $\tau_{K}$ the first exit time of $\Gamma$ from $K$.
Suppose that $\tau_{K}<\infty$ a.s.. Then,
\begin{description}
\item [(i)] For any bounded pathwise variation $\Sigma$ with bounded infinitesimal generator $Y$  which converges uniformly to $\Gamma^{\tau_{K}}$
uniformly, the action has a directional derivative that equals
\begin{eqnarray}
\!\!\!\!\!\!\!\!\!\!\!\!\!\!\!\!\!\!\!\!\!
\left.  \frac{d}{ds}\right\vert _{s=0}S\left(  \Sigma^{s}\right)
	&:=&\underset{s\rightarrow0}{\lim}\frac{1}{s}\left[  S\left(  \Sigma^{s}\right)
-S\left(  \Gamma^{\tau_{K}}\right)  \right]\notag \\
	&=&\int\left\langle {\bf  i}_{Y}
d\theta,\delta\Gamma^{\tau_{K}}\right\rangle -\int\left\langle \widehat
{Y\left[  h\right]  }(\Gamma^{\tau_K}),\delta X\right\rangle +\left\langle \theta\left(
\Gamma^{\tau_{K}}\right)  ,Y\right\rangle -\left\langle \theta\left(
\Gamma^{\tau_{K}}\right)  ,Y\right\rangle _{t=0}, \label{eq 6 action 2}
\end{eqnarray}
where the symbol $\widehat{Y\left[  h\right]  }(\Gamma^{\tau_K}) $ is consistent with the notation introduced in Definition~\ref{action definition}

\item [(ii)] The semimartingale $\Gamma$ satisfies the stochastic Hamiltonian
equations~(\ref{solution Hamiltonian integral}) with initial condition $\Gamma_{0}=m_{0}$ up to time $\tau
_{K}$ if and only if, for any bounded pathwise variation $\Sigma:\left(  -s_{0}%
,s_{0}\right)  \times\mathbb{R}_{+}\times\Omega\rightarrow M$  with bounded infinitesimal generator  which converges uniformly to
$\Gamma^{\tau_{K}}$ and such that $\Sigma_{0}^{s}=m_{0}$ and
$\Sigma_{\tau_{K}}^{s}=\Gamma_{\tau_{K}}$ a.s. for any $s\in\left(
-s_{0},s_{0}\right)  $,%
\[
\left[  \left.  \frac{d}{ds}\right\vert _{s=0}S\left(  \Sigma^{s}\right)
\right]  _{\tau_{K}}=0\text{ \ a.s..}
\]
\end{description}
\end{theorem}

\noindent\textbf{Proof.\ \ }
We first show that the limit (\ref{eq 6 action 2}) exist. Let
$\Sigma$ be an arbitrary bounded pathwise variation converging to $\Gamma$ uniformly and
$Y:\mathbb{R}_{+}\times\Omega\rightarrow TM$ its infinitesimal generator, that we also assume to be bounded. We
have%
\begin{equation}
\frac{1}{s}\left[  S\left(  \Sigma^{s}\right)  -S\left(  \Gamma^{\tau_{K}%
}\right)  \right]  =\frac{1}{s}\left[  \int\left\langle \theta,\delta
\Sigma^{s}\right\rangle -\int\left\langle \theta,\delta\Gamma^{\tau_{K}%
}\right\rangle \right]  -\frac{1}{s}\left[  \int\left\langle \widehat
{h}\left(  \Sigma^{s}\right)  -\widehat{h}\left(  \Gamma^{\tau_{K}}\right)
,\delta X\right\rangle \right]  . \label{eq 7}%
\end{equation}
By Proposition~\ref{derivada pathwise variation}, the first summand in the right hand side of~(\ref{eq 7}) converges $ucp $    to
\[
\int\left\langle {\bf i}_{Y}d\theta,\delta\Gamma^{\tau_{K}}\right\rangle
+\left\langle \theta\left(  \Gamma^{\tau_{K}}\right)  ,Y\right\rangle
-\left\langle \theta\left(  \Gamma^{\tau_{K}}\right)  ,Y\right\rangle _{t=0}.
\]
as $s\rightarrow0$. 
An argument similar to the one leading to Proposition~\ref{derivada pathwise variation} shows that the second summand converges to
$
\int\left\langle \widehat{Y\left[  h\right]  }(\Gamma^{\tau_K}),\delta X\right\rangle
$.
Hence,
\[
\underset{s\rightarrow0}{\lim}\frac{1}{s}\left[  S\left(  \Sigma^{s}\right)
-S\left(  \Gamma^{\tau_{K}}\right)  \right]  =\int\left\langle {\bf  i}_{Y}%
d\theta,\delta\Gamma^{\tau_{K}}\right\rangle -\int\left\langle \widehat
{Y\left[  h\right]  }(\Gamma^{\tau_K}),\delta X\right\rangle +\left\langle \theta\left(
\Gamma^{\tau_{K}}\right)  ,Y\right\rangle -\left\langle \theta\left(
\Gamma^{\tau_{K}}\right)  ,Y\right\rangle _{t=0}.
\]
If we denote by $\eta:=-{\bf i}_{Y}d\theta={\bf i}_{Y}\omega$ the one-form over $\Gamma
^{\tau_{K}}$ built using the vector field $Y$ over $\Gamma^{\tau_{K}}$, the
previous relation may be rewritten as
\begin{equation}
\left[  \left.  \frac{d}{ds}\right\vert _{s=0}S\left(  \Sigma^{s}\right)
\right]  =-\int\left\langle \eta,\delta\Gamma^{\tau_{K}}\right\rangle
-\int\left\langle \mathbf{d}h\left(  \Gamma^{\tau_{K}}\right)  \left(
\omega^{\#}\left(  \eta\right)  \right)  ,\delta X\right\rangle +\left\langle
\theta\left(  \Gamma^{\tau_{K}}\right)  ,Y\right\rangle -\left\langle
\theta\left(  \Gamma^{\tau_{K}}\right)  ,Y\right\rangle _{t=0}. \label{eq 9}%
\end{equation}

We are now going to prove  the assertion in part {\bf (ii)}. Recall that the hypothesis that $\Gamma$ satisfies
the stochastic Hamilton equations up to time $\tau_{K}$ means that
\begin{equation}
\left(  \int\left\langle \beta,\delta\Gamma\right\rangle +\int\left\langle
\left(  \mathbf{d}h\cdot\omega^{\#}\left(  \beta\right)  \right)  \left(
\Gamma\right)  ,\delta X\right\rangle \right)  ^{\tau_{K}}=0 \label{eq 8},
\end{equation}
for any $\beta\in\Omega\left(  M\right)  $. We now show that this expression is also true if we
replace $\beta$ with any process $\eta:\mathbb{R}_{+}\times\Omega\rightarrow
T^{\ast}M$ over $\Gamma$ such that the two Stratonovich integrals involved in
(\ref{eq 8}) are well-defined (for instance if $\beta$ is a semimartingale).
Indeed, invoking (\cite[7.7]{emery}) and Whitney's embedding theorem,
there exist an integer $p\in\mathbb{N}$ such that the manifold $M$ can be seen as an  embedded submanifold of
$\mathbb{R}^{p}$. In this embedded picture, there exists a family of functions $\left\{
f^{1},...,f^{p}\right\}  \subset C^{\infty}\left(  \mathbb{R}^{p}\right)  $
such that the one-form $\eta$ may be written as%
\[
\eta=\sum_{j=1}^{p}Z_{j}\mathbf{d}f^{j},
\]
where the $Z_{j}:\mathbb{R}_{+}\times\Omega\rightarrow\mathbb{R}$,  $j\in\left\{  1,...,p\right\}  $, are real
processes. Moreover, using the properties of the Stratonovich integral (see
\cite[Proposition 7.4]{emery}),%
\begin{align*}
\Big(  \int\left\langle \eta,\delta\Gamma\right\rangle &+\int\left\langle
\left(  \mathbf{d}h\cdot\omega^{\#}\left(  \eta\right)  \right)  \left(
\Gamma\right) ,\delta X\right\rangle \Big)  ^{\tau_{K}} \\
&=\left(  \sum_{j=1}^{p}\int Z_{j}\delta\left(  \int\left\langle
\left\langle \mathbf{d}f^{j},\delta\Gamma\right\rangle +\int\left\langle
\left(  \mathbf{d}h\cdot\omega^{\#}\left(  \mathbf{d}f^{j}\right)  \right)
\left(  \Gamma\right)  ,\delta X\right\rangle \right\rangle \right)  \right)
^{\tau_{K}}\\
& = \sum_{j=1}^{p}\int Z_{j}\delta\left(  \int\left\langle \left\langle
\mathbf{d}f^{j},\delta\Gamma\right\rangle +\int\left\langle \left(
\mathbf{d}h\cdot\omega^{\#}\left(  \mathbf{d}f^{j}\right)  \right)  \left(
\Gamma\right)  ,\delta X\right\rangle \right\rangle \right)  ^{\tau_{K}},
\end{align*}
where the last equality follows from Proposition~\ref{lemma restringir intervalo temporal}. Therefore, since $\mathbf{d}f^{j}$ is a deterministic one-form we can conclude that $\left(  \int\left\langle \left\langle \mathbf{d}f^{j},\delta
\Gamma\right\rangle +\int\left\langle \left(  \mathbf{d}h\cdot\omega
^{\#}\left(  \mathbf{d}f^{j}\right)  \right)  \left(  \Gamma\right)  ,\delta
X\right\rangle \right\rangle \right)  ^{\tau_{K}}=0$, which justifies why~(\ref{eq 8}) also holds if we
replace $\beta\in\Omega\left(  M\right)  $ by an arbitrary integrable one-form $\eta$ over $\Gamma$.

Suppose now that $\Gamma$ satisfies the stochastic Hamilton equations up to
$\tau_{K}$ and let $\Sigma:\left(  -s_{0},s_{0}\right)  \times\mathbb{R}%
_{+}\times\Omega\rightarrow M$ be a pathwise variation like in the statement of the theorem. We want to show that%
\[
\left[  \left.  \frac{d}{ds}\right\vert _{s=0}S\left(  \Sigma^{s}\right)
\right]  _{\tau_{K}}=0\text{ \ a.s..}%
\]
Due to~(\ref{eq 9}), we have that
\begin{equation*}
\left[  \left.  \frac{d}{ds}\right\vert _{s=0}S\left(  \Sigma^{s}\right)
\right]  _{\tau_{K}}   =-\left(  \int\left\langle \eta,\delta\Gamma
^{\tau_{K}}\right\rangle +\int\left\langle \mathbf{d}h\left(  \Gamma^{\tau
_{K}}\right)  \left(  \omega^{\#}\left(  \eta\right)  \right)  ,\delta
X\right\rangle \right)  _{\tau_{K}}
+\left\langle \theta\left(  \Gamma^{\tau_{K}}\right)  ,Y\right\rangle
_{\tau_{K}}-\left\langle \theta\left(  \Gamma^{\tau_{K}}\right)
,Y\right\rangle _{t=0}.
\end{equation*}
Since $\Sigma_{0}^{s}=m_{0}$ and $\Sigma_{\tau_{K}}^{s}=\Gamma_{\tau_{K}}$
a.s. for any $s\in\left(  -s_{0},s_{0}\right)  $, then $Y_{0}=Y_{\tau_{K}}=0$
a.s. and both $\left\langle \theta\left(  \Gamma^{\tau_{K}}\right)
,Y\right\rangle _{\tau_{K}}$ and $\left\langle \theta\left(  \Gamma^{\tau_{K}%
}\right)  ,Y\right\rangle _{t=0}$ vanish. Moreover,%
\begin{align}
 \left(  \int\left\langle \eta,\delta\Gamma^{\tau_{K}}\right\rangle
+\int\left\langle \mathbf{d}h\left(  \Gamma^{\tau_{K}}\right)  \left(
\omega^{\#}\left(  \eta\right)  \right)  ,\delta X\right\rangle \right)
_{\tau_{K}} &=\left(  \left(  \int\left\langle \eta,\delta\Gamma^{\tau_{K}}\right\rangle
+\int\left\langle \mathbf{d}h\left(  \Gamma^{\tau_{K}}\right)  \left(
\omega^{\#}\left(  \eta\right)  \right)  ,\delta X\right\rangle \right)
^{\tau_{K}}\right)  _{\tau_{K}}\nonumber\\
&=\left(  \left(  \int\left\langle \eta,\delta\Gamma\right\rangle
+\int\left\langle \mathbf{d}h\left(  \Gamma\right)  \left(  \omega^{\#}\left(
\eta\right)  \right)  ,\delta X\right\rangle \right)  ^{\tau_{K}}\right)
_{\tau_{K}} \label{eq 13}%
\end{align}
which is zero because of (\ref{eq 8}). In the last equality we have used Proposition~\ref{lemma restringir intervalo temporal}.

Conversely, suppose that $\left[  \left.  \frac{d}{ds}\right\vert
_{s=0}S\left(  \Sigma^{s}\right)  \right]  _{\tau_{K}}=0$ a.s. for arbitrary bounded
pathwise variations tending to $\Gamma^{\tau_K}$ uniformly, like in the statement. We want to show that
(\ref{eq 8}) holds. Since our pathwise variations satisfy that $Y_{0}=Y_{\tau_{K}}=0$ a.s., we obtain that
\begin{equation}
\left[  \left.  \frac{d}{ds}\right\vert _{s=0}S\left(  \Sigma^{s}\right)
\right]  _{\tau_{K}}=-\left(  \int\left\langle \eta,\delta\Gamma^{\tau_{K}%
}\right\rangle +\int\left\langle \mathbf{d}h\left(  \Gamma^{\tau_{K}}\right)
\left(  \omega^{\#}\left(  \eta\right)  \right)  ,\delta X\right\rangle
\right)  _{\tau_{K}}=0 \label{eq 10}%
\end{equation}
where $\eta$ is an arbitrary bounded one form over $\Gamma$. Suppose now that
$\eta$ is a semimartingale. Then $\mathbf{1}_{\left[  0,t\right]  }
\eta:\mathbb{R}_{+}\times\Omega\rightarrow T^{\ast}M$ is again bounded and
expressions 
\[
\int\left\langle \mathbf{1}_{\left[  0,t\right]  }\eta,\delta\Gamma^{\tau_{K}%
}\right\rangle \text{ \ \ and \ \ }\int\left\langle \mathbf{d}h\left(
\Gamma^{\tau_{K}}\right)  \left(  \omega^{\#}\left(  \mathbf{1}_{\left[
0,t\right]  }\eta\right)  \right)  ,\delta X\right\rangle
\]
are well-defined by Proposition~\ref{with indicator integral} because both
$\Gamma^{\tau_{K}}$ and $X$ are continuos semimartingales. We already saw in
(\ref{eq 13}) that (\ref{eq 10}) is equivalent to%
\[
\left(  \int\left\langle \eta,\delta\Gamma\right\rangle +\int\left\langle
\mathbf{d}h\left(  \Gamma\right)  \left(  \omega^{\#}\left(  \eta\right)
\right)  ,\delta X\right\rangle \right)  _{\tau_{K}}=0.
\]
Replacing $\eta$ by $\mathbf{1}_{\left[  0,t\right]  }\eta$ in (\ref{eq 10})
and using again the Proposition~\ref{with indicator integral}, we write
\begin{align*}
0  &  =\left(  \int\left\langle \mathbf{1}_{\left[  0,t\right]  }\eta
,\delta\Gamma\right\rangle +\int\left\langle \mathbf{d}h\left(  \Gamma\right)
\left(  \omega^{\#}\left(  \mathbf{1}_{\left[  0,t\right]  }\eta\right)
\right)  ,\delta X\right\rangle \right)  _{\tau_{K}} =\left(  \left(  \int\left\langle \eta,\delta\Gamma\right\rangle
+\int\left\langle \mathbf{d}h\left(  \Gamma\right)  \left(  \omega^{\#}\left(
\eta\right)  \right)  ,\delta X\right\rangle \right)  ^{t}\right)  _{\tau_{K}
}\\
&  =\left(  \int\left\langle \eta,\delta\Gamma\right\rangle +\int\left\langle
\mathbf{d}h\left(  \Gamma\right)  \left(  \omega^{\#}\left(  \eta\right)
\right)  ,\delta X\right\rangle \right)  _{t\wedge\tau_{K}}=\left(  \left(  \int\left\langle \eta,\delta\Gamma\right\rangle
+\int\left\langle \mathbf{d}h\left(  \Gamma\right)  \left(  \omega^{\#}\left(
\eta\right)  \right)  ,\delta X\right\rangle \right)  ^{\tau_{K}}\right)  _{t}.
\end{align*}
Since $t$ is arbitrary this implies that the process $\left(  \int\left\langle \eta,\delta
\Gamma\right\rangle +\int\left\langle \mathbf{d}h\left(  \Gamma\right)
\left(  \omega^{\#}\left(  \eta\right)  \right)  ,\delta X\right\rangle
\right)  ^{\tau_{K}}$ is identically zero, as required.
\quad $\blacksquare$

\section{Proofs and auxiliary results}

\subsection{Proof of Proposition~\ref{prop differentiation processes}}
\label{proof of proposition crucial}

\noindent Before proving the proposition, we recall a technical lemma dealing with
the convergence of sequences in a metric space.

\begin{lemma}
\label{lemma serie convergence}
Let $\left(  E,d\right)  $ be a metric space.
Let $\left\{  x_{n}\right\}  _{n\in\mathbb{N}}$ be a sequence of functions
$x_{n}:\left(  0,\delta\right)  \rightarrow E$ where $\left(  0,\delta\right)
\subset\mathbb{R}$ is an open interval of the real line. Suppose that $x_{n} $
converges uniformly on $\left(  0,\delta\right)  $ to a function $x.$
Additionally, suppose that for any $n$, the limits $\lim\limits_{s\rightarrow0}x_{n}\left(
s\right)  =x_{n}^{\ast}\in E$ exist and so does $\lim\limits_{n\rightarrow\infty
}x_{n}^{\ast}$. Then
\[
\lim_{s\rightarrow0}x\left(  s\right)  =\lim_{n\rightarrow\infty}x_{n}^{\ast
}\text{.}
\]
\end{lemma}

\noindent\textbf{Proof.\ \ }
Let $\varepsilon>0$ be an arbitrary real number. We have%
\[
d\left(  x\left(  s\right)  ,\lim_{n\rightarrow\infty}x_{n}^{\ast}\right)
\leq d\left(  x\left(  s\right)  ,x_{k}\left(  s\right)  \right)  +d\left(
x_{k}\left(  s\right)  ,x_{k}^{\ast}\right)  +d\left(  x_{k}^{\ast}%
,\lim_{n\rightarrow\infty}x_{n}^{\ast}\right)  .
\]
From the definition of limit and since $x_{k}\left(  s\right)  $ converges
uniformly to $x$ on $\left(  0,\delta\right)  ,$ we can choose $k_{0}$ such
that  $d\left(  x_{k}^{\ast},\lim_{n\rightarrow\infty}
x_{n}^{\ast}\right)  <\frac{\varepsilon}{3}$ and $d\left(  x\left(
s\right)  ,x_{k}\left(  s\right)  \right)  <\frac{\varepsilon}{3}$, simultaneously for any
$k\geq k_{0}$. Additionally, since $\lim_{s\rightarrow0}x_{k}\left(  s\right)
=x_{k}^{\ast}$ we choose $s_{0}$ small enough such that $d\left(  x_{k}\left(
s\right)  ,x_{k}^{\ast}\right)<\frac{\varepsilon}{3}  $, for any $s<s_{0}$. Thus,%
\[
d\left(  x\left(  s\right)  ,\lim_{n\rightarrow\infty}x_{n}^{\ast}\right)
<\varepsilon
\]
for any $s<s_{0}.$ Since $\varepsilon>0$ is arbitrary, we conclude that
$\lim\limits_{s\rightarrow0}x\left(  s\right)  =\lim\limits_{n\rightarrow\infty}x_{n}^{\ast
}$.
\quad $\blacktriangledown$

\medskip

\noindent\textbf{Proof of Proposition~\ref{prop differentiation processes}.\ \ } First of
all, the second equality in~(\ref{derivative of functional}) is a straightforward consequence
of~\cite[page 93]{emery}. Now,  let
$\left\{ U_{k}\right\}  _{k\in\mathbb{N}}$ be a countable open covering of $M$ by
coordinate patches. By~\cite[Lemma 3.5]{emery} there exists a sequence $\left\{ 
\tau_{m}\right\} _{m\in\mathbb{N}}$ of 
stopping times such that
$\tau_{0}=0,$ $\tau_{m}\leq\tau_{m+1},$ $\sup_{m}\tau_{m}=\infty,$ a.s., and that,
on each of the sets
\begin{equation}
\label{local interval}
\left[  \tau_{m} ,\tau_{m+1}\right]\cap \{ \tau_{m} <\tau_{m+1}\}
:=\left\{  \left(  t,\omega\right)  \in\mathbb{R}_{+}\times\Omega\mid\tau
_{m+1}\left(  \omega\right)  >\tau_{m}\left(  \omega\right)  \text{ and }%
t\in\left[  \tau_{m}\left(  \omega\right)  ,\tau_{m+1}\left(  \omega\right)
\right]  \right\}
\end{equation}
the semimartingale $\Gamma$ takes its values in  one of the elements  of the
family
$\left\{ U_{k}\right\}  _{k\in\mathbb{N}}$.  

Second, the statement of the proposition is formulated in terms of Stratonovich
integrals. However, the proof will be carried out in the context of It\^o integration
since we will use several times the notion of uniform convergence on compacts
in probability (ucp) which behaves well only with respect to this integral. Regarding
this point we recall that by the very definition of the Stratonovich
integral of a $1$-form $\alpha$ along a semimartingale $\Gamma$ we have that
\begin{equation}
\label{translation Ito Stratonovich}
\int\left\langle \varphi_{s}^{\ast}\alpha,\delta\Gamma\right\rangle
=\int\left\langle d_{2}\left(  \varphi_{s}^{\ast}\alpha\right)  ,d\Gamma
\right\rangle \text{ \ and \ }\int\left\langle \pounds_{Y}\alpha,\delta
\Gamma\right\rangle =\int\left\langle d_{2}\left(  \pounds_{Y}\alpha\right)
,d\Gamma\right\rangle.
\end{equation}

The proof of the proposition follows directly from Lemma~\ref{lemma serie convergence} by
applying it to the sequence of functions given by
\begin{equation*}
x_{n}\left(  s\right)  :=\left(  \int\left\langle \frac{1}{s}\left[
d_{2}\left(  \varphi_{s}^{\ast}\alpha\right)  -d_{2}\left(  \alpha\right)
\right]  ,d\Gamma\right\rangle \right)  ^{\tau_{n}}.
\end{equation*}
This sequence lies in the space $\mathbb{D} $ of c\`agl\`ad  processes endowed with the
topology of the ucp convergence. We recall that this space is metric~\cite[page
57]{protter} and hence we are in the conditions of Lemma~\ref{lemma serie convergence}.
In the following points we verify that the rest of the hypotheses of this result are
satisfied.

\medskip

\noindent {\bf (i) The sequence of functions $\{x _n(s)\}_{n \in  \mathbb{N}}$ converges
uniformly to}
\begin{equation*}
x (s):=\int\left\langle \frac{1}{s}\left[
d_{2}\left(  \varphi_{s}^{\ast}\alpha\right)  -d_{2}\left(  \alpha\right)
\right]  ,d\Gamma\right\rangle.
\end{equation*}  
The pointwise convergence is a consequence of part {\bf (i)} in Proposition~\ref{tau
series}. Moreover, in the proof of that result we saw that if $d: \mathbb{D} \times 
\mathbb{D} \rightarrow \mathbb{R}_+ $ is a distance function  function associated to the
ucp convergence, then for any $t \in \mathbb{R}_+ $ and any $s \in (0, \epsilon)$,
$d(x _n (s), x (s))\leq P(\{ \tau _n< t\})$. Since the right hand side of this inequality
does not depend on $s$ and $P(\{ \tau _n< t\})\rightarrow 0$ as $n \rightarrow \infty$, the
uniform convergence follows.

\medskip

\noindent {\bf (ii) } 
\begin{equation*}
\lim\limits_{\overset{ucp}{s \rightarrow 0}} x _n(s)=
\left(  \int\left\langle d_{2}\left(  \pounds_{Y}\alpha\right)  ,d\Gamma
\right\rangle \right)  ^{\tau_{n}}=: x _n^\ast. 
\end{equation*}
By the construction of the covering $\left\{ U_{k}\right\}  _{k\in\mathbb{N}}$  and of the
stopping times
$\left\{ 
\tau_{m}\right\} _{m\in\mathbb{N}}$, there exists a $k (m) \in \mathbb{N} $  such that  the
semimartingale  $\Gamma$  takes its values in $U_{k (m)}$ when evaluated in the stochastic
interval $( \tau_n, \tau_{n+1}] \subset [ \tau_n, \tau_{n+1}]\cap \{ \tau_n<
\tau_{n+1}\}
$. Now, since $d_2 $ is a linear operator and 
$  \frac{1}{s}\left( \left(  \varphi_{s}^{\ast}
\alpha\right)- \alpha\right)(m) \overset{s \rightarrow 0}{\longrightarrow}  
  \pounds_{Y}\alpha (m)  $,  for any $m \in M $, we have that $ 
\frac{1}{s}\left(  d_{2}\left(  \varphi_{s}^{\ast}
\alpha\right)-   d_{2}\alpha\right)(m) \overset{s \rightarrow 0}{\longrightarrow}  
d_{2}\left(  \pounds_{Y}\alpha\right) (m)  $. Moreover, a straightforward application of
Taylor's theorem shows that 
$ 
\frac{1}{s}\left(  d_{2}\left(  \varphi_{s}^{\ast}
\alpha\right)-   d_{2}\alpha\right)|_{U_{k (m)}} \overset{s \rightarrow
0}{\longrightarrow}   d_{2}\left(  \pounds_{Y}\alpha\right)|_{U_{k (m)}}   
$
uniformly, using a Euclidean norm in $\tau^\ast U_{k (m)}$ (we recall that $U _{k (m)}$ is
a coordinate patch). This fact immediately implies that 
${\bf 1}_{( \tau_n, \tau_{n+1}]}\frac{1}{s}\left(  d_{2}\left(  \varphi_{s}^{\ast}
\alpha\right)-   d_{2}\alpha\right)(\Gamma) \overset{s \rightarrow
0}{\longrightarrow}   
{\bf 1}_{( \tau_n, \tau_{n+1}]}d_{2}\left(  \pounds_{Y}\alpha\right)(\Gamma)
$ in ucp. As by construction the It\^o integral behaves well when we apply it
to a ucp convergent sequence of processes 
we have that
\begin{equation}
\label{with tau inside we will see}
\lim\limits_{\overset{ucp}{s \rightarrow 0}} \int
{\bf 1}_{( \tau_n, \tau_{n+1}]}\left\langle
\frac{1}{s}\left(  d_{2}\left(  \varphi_{s}^{\ast}
\alpha\right)-   d_{2}\alpha\right)(\Gamma), d \Gamma\right\rangle=
\int
{\bf 1}_{( \tau_n, \tau_{n+1}]}\left\langle
d_{2}\left(  \pounds_{Y}\alpha\right)(\Gamma), d \Gamma\right\rangle.
\end{equation}
Consequently,
\begin{align*}
\lim\limits_{\overset{ucp}{s \rightarrow 0}}&\left(\int\left\langle \frac{1}{s}\left[ 
d_{2}\left(  \varphi_{s}^{\ast }\alpha\right)  -d_{2}\left(  \alpha\right)  \right] 
,d\Gamma\right\rangle
\right)  ^{\tau_{n}}\\
 &=\lim\limits_{\overset{ucp}{s \rightarrow 0}}\sum_{m=0}^{n-1}\left[\left( 
\int\left\langle
\frac{1}{s}\left[  d_{2}\left(
\varphi_{s}^{\ast}\alpha\right)  -d_{2}\left(  \alpha\right)  \right]
,d\Gamma\right\rangle \right)  ^{\tau_{m+1}}-\left(  \int\left\langle \frac
{1}{s}\left[  d_{2}\left(  \varphi_{s}^{\ast}\alpha\right)  -d_{2}\left(
\alpha\right)  \right]  ,d\Gamma\right\rangle \right)  ^{\tau_{m}}\right]\\
 &=\lim\limits_{\overset{ucp}{s \rightarrow 0}}\sum_{m=0}^{n-1}\int\mathbf{1}_{\left( 
\tau_{m},\tau_{m+1}\right] }\left\langle \frac{1}{s}\left(  d_{2}\left( 
\varphi_{s}^{\ast}\alpha\right) -d_{2}\alpha\right)  ,d\Gamma\right\rangle
=
\sum_{m=0}^{n-1}\int\mathbf{1}_{\left(  \tau_{m},\tau_{m+1}\right]
}\left\langle d_{2}\left(  \pounds_{Y}\alpha\right)  ,d\Gamma\right\rangle \\
	&=\left(  \int\left\langle d_{2}\left(  \pounds_{Y}\alpha\right)  ,d\Gamma
\right\rangle \right)  ^{\tau_{n}},
\end{align*}
where in the second equality we have used Proposition~\ref{lemma restringir intervalo
temporal} and the third one follows from~(\ref{with tau inside we will see}).  

\medskip

\noindent {\bf (iii)} 
\begin{equation*}
\lim\limits_{n \rightarrow  \infty}x _n^\ast =\int\left\langle d_{2}\left(  \pounds_{Y}\alpha\right)
,d\Gamma\right\rangle.
\end{equation*}
It is a straightforward consequence of  Proposition~\ref{tau series}.

The  equation~(\ref{derivative of functional})  follows from Lemma~\ref{lemma
serie convergence} applied to the sequences $\{x _n\}_{n \in \mathbb{N}}$ and  $\{x
_n^\ast \}_{n
\in \mathbb{N}}$, and using the statements in {\bf (i)}, {\bf (ii)}, and {\bf (iii)}. \quad
$\blacksquare$

\subsection{Proof of Proposition~\ref{derivada pathwise variation}}
\label{appendix proposition derivada pathwise}

We will start the proof by introducing three preparatory results.

\begin{lemma}
\label{producto ucp convergencia}Let $\left\{  X_{n}\right\}  _{n\in
\mathbb{N}}$ and $\left\{  Y_{n}\right\}  _{n\in\mathbb{N}}$ be two sequences
of real valued processes converging in $ucp$ to a couple of processes $X$ and $Y$
respectively. Suppose that, for any $t\in\mathbb{R}_{+}$, the random variables $\sup_{n\in
\mathbb{N}}\sup_{0\leq s\leq t}\left\vert \left(  X_{n}\right)  _{s}
\right\vert$  and $\sup_{0\leq s\leq t}\left\vert Y_{s}\right\vert$ are bounded (their images lie in a compact set of $\mathbb{R}$). Then, the sequence $X_{n}Y_{n}$ converges in
$ucp$ to $XY$ as $n\rightarrow\infty$.
\end{lemma}

\noindent\textbf{Proof.\ \ }
We need to prove that for any $\varepsilon>0$ and any $t\in\mathbb{R}_{+}$,%
\[
P\left(  \left\{  \sup_{0\leq s\leq t}\left\vert \left(  X_{n}Y_{n}\right)
_{s}-\left(  XY\right)  _{s}\right\vert \leq\varepsilon\right\}  \right)
\underset{n\rightarrow\infty}{\longrightarrow}1.
\]
First of all, note that
\[
\sup_{0\leq s\leq t}\left\vert \left(  X_{n}Y_{n}\right)  _{s}-\left(
XY\right)  _{s}\right\vert \leq\sup_{0\leq s\leq t}\left\vert X_{n}\right\vert
\left\vert Y_{n}-Y\right\vert +\sup_{0\leq s\leq t}\left\vert Y\right\vert
\left\vert X_{n}-X\right\vert.
\]
Hence, we have
\begin{align*}
\left\{  \sup_{0\leq s\leq t}\left\vert \left(  X_{n}Y_{n}\right)
_{s}-\left(  XY\right)  _{s}\right\vert \leq\varepsilon\right\}   &
\supseteq\left\{  \sup_{0\leq s\leq t}\left\vert X_{n}\right\vert \left\vert
Y_{n}-Y\right\vert +\sup_{0\leq s\leq t}\left\vert Y\right\vert \left\vert
X_{n}-X\right\vert \leq\varepsilon\right\} \\
&  \supseteq\left\{  \sup_{0\leq s\leq t}\left\vert X_{n}\right\vert
\left\vert Y_{n}-Y\right\vert \leq\frac{\varepsilon}{2}\right\}  \cap\left\{
\sup_{0\leq s\leq t}\left\vert Y\right\vert \left\vert X_{n}-X\right\vert
\leq\frac{\varepsilon}{2}\right\}  .
\end{align*}
Denote
\begin{equation*}
A_{n}   :=\left\{  \sup_{0\leq s\leq t}\left\vert X_{n}\right\vert
\left\vert Y_{n}-Y\right\vert \leq\frac{\varepsilon}{2}\right\}  , \quad\text{and}\quad
B_{n}   :=\left\{  \sup_{0\leq s\leq t}\left\vert Y\right\vert \left\vert
X_{n}-X\right\vert \leq\frac{\varepsilon}{2}\right\}  ,
\end{equation*}
and let $c$ be a constant such that $\sup_{n\in\mathbb{N}}\sup_{0\leq s\leq
t}\left\vert \left(  X_{n}\right)  _{s}\right\vert <c$ and $\sup_{0\leq s\leq
t}\left\vert Y_{s}\right\vert <c$, available by the boundedness hypothesis. Then,%
\begin{align*}
1  &  \geq P\left(  A_{n}\right)  \geq P\left(  \left\{  \sup_{0\leq s\leq
t}\left\vert Y_{n}-Y\right\vert \leq\frac{\varepsilon}{2c}\right\}  \right)
\underset{n\rightarrow\infty}{\longrightarrow}1,\\
1  &  \geq P\left(  B_{n}\right)  \geq P\left(  \left\{  \sup_{0\leq s\leq
t}\left\vert X_{n}-X\right\vert \leq\frac{\varepsilon}{2c}\right\}  \right)
\underset{n\rightarrow\infty}{\longrightarrow}1.
\end{align*}
Thus, $P\left(  A_{n}\right)  \rightarrow1$ and $P\left(  B_{n}\right)
\rightarrow1$ as $n\rightarrow\infty$. But as $P\left(  A_{n}\cap
B_{n}\right)  =P\left(  A_{n}\right)  +P\left(  B_{n}\right)  -P\left(
A_{n}\cup B_{n}\right)  $, we conclude that
\[
P\left(  A_{n}\cap B_{n}\right)  \underset{n\rightarrow\infty}{\longrightarrow
}1.
\]
Since $A_{n}\cap B_{n}\subseteq\left\{  \sup_{0\leq s\leq t}\left\vert \left(
X_{n}Y_{n}\right)  _{s}-\left(  XY\right)  _{s}\right\vert \leq\varepsilon
\right\}  $, we obtain%
\[
P\left(  \left\{  \sup_{0\leq s\leq t}\left\vert \left(  X_{n}Y_{n}\right)
_{s}-\left(  XY\right)  _{s}\right\vert \leq\varepsilon\right\}  \right)
\underset{n\rightarrow\infty}{\longrightarrow}1. \qquad \blacktriangledown
\]

\begin{lemma}
\label{ucp implica probabilidad}Let $\left\{  X_{n}\right\}  _{n\in\mathbb{N}%
}$ be a sequence of real processes converging in $ucp$ to a process $X$.
Let $\tau$ be a stopping time such that $\tau<\infty$ a.s.. Then, the sequence
of random variables $\left\{  \left(  X_{n}\right)  _{\tau}\right\}
_{n\in\mathbb{N}}$ converge in probability to $\left(  X\right)  _{\tau}$.
\end{lemma}

\noindent\textbf{Proof.\ \ }
First of all we show that since $\tau<\infty$ a.s., then $P\left(  \left\{  \tau>t\right\}  \right)  $ converges
to zero as $t\rightarrow\infty$. By contradiction, suppose that this is not the case. Then,
denoting $A_{n}:=\left\{  \tau>n\right\}  $, we have that $A_{n+1}\subseteq
A_{n}$, so $P\left(  A_{n}\right)  $ forms a non-increasing sequence of real
numbers in the interval $\left[  0,1\right]  $. Since this sequence is bounded
below, it must have a limit. This limit corresponds to the probability of
the event $\left\{  \tau=\infty\right\}  $. If it is strictly positive then
there is a contradiction with the fact that $\tau<\infty$ a.s.. So $P\left(
\left\{  \tau>t\right\}  \right)  $ tends to zero as $t\rightarrow\infty$.

We now prove the statement of the lemma. Take some $\varepsilon>0$ and an auxiliary $t\in\mathbb{R}_{+}$. The set
$\left\{  \left\vert \left(  X_{n}\right)  _{\tau}-X_{\tau}\right\vert
>\varepsilon\right\}  $ can be decomposed as the disjoint union of the following two
events,
\[
\left(  \left\{  \left\vert \left(  X_{n}\right)  _{\tau}-X_{\tau}\right\vert
>\varepsilon\right\}  \cap\left\{  \tau\leq t\right\}  \right)  \bigcup\left(
\left\{  \left\vert \left(  X_{n}\right)  _{\tau}-X_{\tau}\right\vert
>\varepsilon\right\}  \cap\left\{  \tau>t\right\}  \right)  .
\]
The first one is contained in the set $\left\{  \sup_{0\leq s\leq t}\left\vert
\left(  X_{n}\right)  _{s}-X_{s}\right\vert >\varepsilon\right\}  $ whose
probability, by hypothesis, converges to zero as $n\rightarrow\infty$. Regarding
the second one,%
\[
P\left(  \left\{  \left\vert \left(  X_{n}\right)  _{\tau}-X_{\tau}\right\vert
>\varepsilon\right\}  \cap\left\{  \tau>t\right\}  \right)  \leq P\left(
\left\{  \tau>t\right\}  \right)  .
\]
But $P\left(  \left\{  \tau>t\right\}  \right)  $ can be made arbitrarily small by taking the auxiliary $t$ big enough. In conclusion, for any
$\varepsilon>0$,
\[
P\left(  \left\{  \left\vert \left(  X_{n}\right)  _{\tau}-X_{\tau}\right\vert
>\varepsilon\right\}  \right)  \underset{n\rightarrow\infty}{\longrightarrow
}0
\]
in probability.
\quad $\blacktriangledown$

\begin{lemma}
\label{ucp implica ucp parado}Let $\left\{  X_{n}\right\}  _{n\in\mathbb{N}}$
be a sequence of real processes converging in $ucp$ to a real process $X$ and
$\tau$ a stopping time. Then, the stopped sequence $\left\{  X_{n}^{\tau
}\right\}  _{n\in\mathbb{N}}$ converges in $ucp$ to $X^{\tau}$ as well.
\end{lemma}

\noindent\textbf{Proof.\ \ }
We just need to observe that, for any $t\in\mathbb{R}_{+}$,%
\[
\sup_{0\leq s\leq t}\left\vert \left(  X_{n}^{\tau}\right)  _{s}-X_{s}^{\tau
}\right\vert =\sup_{0\leq s\leq t}\left\vert \left(  X_{n}\right)
_{\tau\wedge s}-X_{\tau\wedge s}\right\vert \leq\sup_{0\leq s\leq t}\left\vert
\left(  X_{n}\right)  _{s}-X_{s}\right\vert
\]
and, consequently, for any $\varepsilon>0$,%
\[
\left\{  \sup_{0\leq s\leq t}\left\vert \left(  X_{n}\right)  _{s}%
-X_{s}\right\vert \leq\varepsilon\right\}  \subseteq\left\{  \sup_{0\leq s\leq
t}\left\vert \left(  X_{n}^{\tau}\right)  _{s}-X_{s}^{\tau}\right\vert
\leq\varepsilon\right\}  .
\]
Hence, since by hypothesis $P\left(  \left\{  \sup_{0\leq s\leq t}\left\vert \left(  X_{n}\right)
_{s}-X_{s}\right\vert \leq\varepsilon\right\}  \right)  $ converges to $1$ as
$n\rightarrow\infty$, then so does $P\left(  \left\{  \sup_{0\leq s\leq
t}\left\vert \left(  X_{n}^{\tau}\right)  _{s}-X_{s}^{\tau}\right\vert
\leq\varepsilon\right\}  \right)  $.
\quad $\blacktriangledown$

\medskip

We now proceed with the proof of the proposition.  We will start by using Whitney's Embedding Theorem and the remarks in~\cite[\S 7.7]{emery} to visualize $M$ as an embedded submanifold of $\mathbb{R}^p $ , for some $p\in\mathbb{N}$, and to write down  our Stratonovich integrals as \emph{real} Stratonovich integrals. Indeed,   there exists a family of functions $\left\{
h^{1},...,h^{p}\right\}  \subset C^{\infty}\left(  \mathbb{R}^{p}\right)  $
such that, in the embedded picture, the one form $\alpha$ can be written as $\alpha=\sum_{j=1}^{p}%
Z_{j}\mathbf{d}h^{j}$, where $Z_{j}\in C^{\infty}\left(  \mathbb{R}%
^{p}\right)  $ for $j\in\left\{  1,...,p\right\}  $. Therefore, using the
properties of the Stratonovich integral (see \cite[Proposition 7.4]{emery}),%
\begin{equation}
\frac{1}{s}\left[  \int\left\langle \alpha,\delta\Sigma^{s}\right\rangle
-\int\left\langle \alpha,\delta\Gamma^{\tau_{K}}\right\rangle \right]
=\sum_{j=1}^{p}\frac{1}{s}\left[  \int Z_{j}\left(  \Sigma^{s}\right)
\delta\left(  h^{j}\left(  \Sigma^{s}\right)  \right)  -\int Z_{j}\left(
\Gamma^{\tau_K}\right)  \delta\left(  h^{j}\left(  \Gamma^{\tau_{K}}\right)  \right)
\right]  . \label{eq plus and minus}
\end{equation}
Adding and subtracting the term $\sum_{j=1}^{p}\int Z_{j}\left(  \Sigma
^{s}\right)  \delta h^{j}\left(  \Gamma^{\tau_{K}}\right)  $ in the right hand side of~(\ref{eq plus and minus}), we have%
\begin{align}
\frac{1}{s}\left[  \int\left\langle \alpha,\delta\Sigma^{s}\right\rangle
-\int\left\langle \alpha,\delta\Gamma^{\tau_{K}}\right\rangle \right]   &
=\sum_{j=1}^{p}\underset{(1)}{\underbrace{\frac{1}{s}\left[  \int Z_{j}\left(
\Sigma^{s}\right)  \delta h^{j}\left(  \Sigma^{s}\right)  -\int Z_{j}\left(
\Sigma^{s}\right)  \delta h^{j}\left(  \Gamma^{\tau_{K}}\right)  \right]  }%
}\nonumber\\
&  +\sum_{j=1}^{p}\underset{(2)}{\underbrace{\frac{1}{s}\left[  \int\left(
Z_{j}\left(  \Sigma^{s}\right)  -Z_{j}\left(  \Gamma^{\tau_{K}}\right)
\right)  \delta h^{j}\left(  \Gamma^{\tau_{K}}\right)  \right]  }}.
\label{eq 3}%
\end{align}
We are going to study the terms (1) and (2) separately. We start by considering
\[\sigma_{n}=\left\{  0=T_{0}^{n}\leq T_{1}^{n}\leq \ldots\leq T_{k_{n}}^{n}
<\infty\right\},\] 
a sequence of random partitions that tends to the identity (in the sense of~\cite[page 64]{protter}).

\medskip

\noindent {\bf The expression (1):} We want to study the $ucp$ convergence of $\frac{1}{s}\left[  \int
Z_{j}\left(  \Sigma^{s}\right)  \delta h^{j}\left(  \Sigma^{s}\right)  -\int
Z_{j}\left(  \Sigma^{s}\right)  \delta h^{j}\left(  \Gamma^{\tau_{K}}\right)
\right]  $ as $s\rightarrow0$. Define
\begin{align*}
x_{n}\left(  s\right)   &  :=\frac{1}{s}\left(  \sum_{i=0}^{k_{n}-1}\frac
{1}{2}\left(  Z_{j}\left(  \Sigma^{s}\right)  _{T_{i+1}^{n}}+Z_{j}\left(
\Sigma^{s}\right)  _{T_{i}^{n}}\right)  \left(  h^{j}\left(  \Sigma
^{s}\right)  ^{T_{i+1}^{n}}-h^{j}\left(  \Sigma^{s}\right)  ^{T_{i}^{n}%
}\right)  \right. \\
&  \left.  -\sum_{i=0}^{k_{n}-1}\frac{1}{2}\left(  Z_{j}\left(  \Sigma
^{s}\right)  _{T_{i+1}^{n}}+Z_{j}\left(  \Sigma^{s}\right)  _{T_{i}^{n}%
}\right)  \left(  h^{j}\left(  \Gamma^{\tau_{K}}\right)  ^{T_{i+1}^{n}}%
-h^{j}\left(  \Gamma^{\tau_{K}}\right)  ^{T_{i}^{n}}\right)  \right) \\
&  =\sum_{i=0}^{k_{n}-1}\frac{1}{2}\left(  Z_{j}\left(  \Sigma^{s}\right)
_{T_{i+1}^{n}}+Z_{j}\left(  \Sigma^{s}\right)  _{T_{i}^{n}}\right)  \left(
\frac{h^{j}\left(  \Sigma^{s}\right)  ^{T_{i+1}^{n}}-h^{j}\left(  \Gamma
^{\tau_{K}}\right)  ^{T_{i+1}^{n}}}{s}-\frac{h^{j}\left(  \Sigma^{s}\right)
^{T_{i}^{n}}-h^{j}\left(  \Gamma^{\tau_{K}}\right)  ^{T_{i}^{n}}}{s}\right)  .
\end{align*}
which corresponds to the discretization of the Stratonovich integrals
$\frac{1}{s}\left[  \int Z_{j}\left(  \Sigma^{s}\right)  \delta h^{j}\left(
\Sigma^{s}\right)  \right.  \allowbreak-\left.  \int Z_{j}\left(  \Sigma
^{s}\right)  \delta h^{j}\left(  \Gamma^{\tau_{K}}\right)  \right]  $ using the random partitions
of $\sigma_{n}$. Indeed, by \cite[Corollary 1, page 291]{protter},
\[
x_{n}\left(  s\right)  \underset{n\rightarrow\infty}{\underset{ucp}
{\longrightarrow}}\frac{1}{s}\left[  \int Z_{j}\left(  \Sigma^{s}\right)
\delta h^{j}\left(  \Sigma^{s}\right)  -\int Z_{j}\left(  \Sigma^{s}\right)
\delta h^{j}\left(  \Gamma^{\tau_{K}}\right)  \right]  .
\]
On the other hand, as $T_{i}^{n}<\infty$ a.s. for any $i\in\left\{
1,...,k_{n}\right\}  $, part {\bf (i)} in Definition~\ref{def pathwise variation} and Lemma~\ref{ucp implica probabilidad} imply that
\[
\frac{1}{2}\left(  Z_{j}\left(  \Sigma^{s}\right)  _{T_{i+1}^{n}}+Z_{j}\left(
\Sigma^{s}\right)  _{T_{i}^{n}}\right)  \underset{s\rightarrow0}%
{\underset{ucp}{\longrightarrow}}\frac{1}{2}\left(  Z_{j}\left(  \Gamma
^{\tau_{K}}\right)  _{T_{i+1}^{n}}+Z_{j}\left(  \Gamma^{\tau_{K}}\right)
_{T_{i}^{n}}\right)
\]
The convergence above is in
probability but, for convenience, we prefer to regard these random variables as trivial
processes. Furthermore, part {\bf (ii)} in Definition~\ref{def pathwise variation} and Lemma~\ref{ucp implica ucp parado}
\begin{align*}
&  \frac{h^{j}\left(  \Sigma^{s}\right)  ^{T_{i+1}^{n}}-h^{j}\left(
\Gamma^{\tau_{K}}\right)  ^{T_{i+1}^{n}}}{s}=\left(  \frac{h^{j}\left(
\Sigma^{s}\right)  -h^{j}\left(  \Gamma^{\tau_{K}}\right)  }{s}\right)
^{T_{i+1}^{n}}\underset{s\rightarrow0}{\underset{ucp}{\longrightarrow}%
}Y\left[  h^{j}\right]  ^{T_{i+1}^{n}},\\
&  \frac{h^{j}\left(  \Sigma^{s}\right)  ^{T_{i}^{n}}-h^{j}\left(
\Gamma^{\tau_{K}}\right)  ^{T_{i}^{n}}}{s}=\left(  \frac{h^{j}\left(
\Sigma^{s}\right)  -h^{j}\left(  \Gamma^{\tau_{K}}\right)  }{s}\right)
^{T_{i}^{n}}\underset{s\rightarrow0}{\underset{ucp}{\longrightarrow}}Y\left[
h^{j}\right]  ^{T_{i}^{n}}%
\end{align*}
by Definition (\ref{def pathwise variation}) item 3 and Lemma
(\ref{ucp implica ucp parado}). Now, since by hypothesis $\Sigma$ and $Y $  are bounded then so are 
$\frac{1}{2}\left(  Z_{j}\left(
\Sigma^{s}\right)  _{T_{i+1}^{n}}+Z_{j}\left(  \Sigma^{s}\right)  _{T_{i}^{n}
}\right)  $ and $Y\left[  h^{j}\right]  ={\bf i}_{Y}\mathbf{d}h^{j}$ ($\mathbf{d}h^{j}$ is only evaluated on the compact
$K$ since $Y $ is a vector field over $\Gamma^{\tau_K} $) and hence by Lemma~\ref{producto ucp convergencia}
\[
x_{n}\left(  s\right)  \underset{s\rightarrow0}{\underset{ucp}{\longrightarrow
}}\sum_{i=0}^{k_{n-1}}\frac{1}{2}\left(  Z_{j}\left(  \Gamma^{\tau_{K}%
}\right)  _{T_{i+1}^{n}}+Z_{j}\left(  \Gamma^{\tau_{K}}\right)  _{T_{i}^{n}%
}\right)  \left(  Y\left[  h^{j}\right]  ^{T_{i+1}^{n}}-Y\left[  h^{j}\right]
^{T_{i}^{n}}\right)  =:x_{n}^{\ast}%
\]
In addition, by \cite[Corollary 1,
page 291]{protter},%
\[
x_{n}^{\ast}\underset{n\rightarrow\infty}{\underset{ucp}{\longrightarrow}}\int
Z_{j}\left(  \Gamma^{\tau_{K}}\right)  \delta\left(  Y\left[  h^{j}\right]
\right)  .
\]
Hence, by Lemma~\ref{lemma serie convergence} we conclude that
\begin{equation}
\label{we will substitute 1}
\frac{1}{s}\left[  \int Z_{j}\left(  \Sigma^{s}\right)  \delta h^{j}\left(
\Sigma^{s}\right)  -\int Z_{j}\left(  \Sigma^{s}\right)  \delta h^{j}\left(
\Gamma^{\tau_{K}}\right)  \right]  \underset{s\rightarrow0}{\underset
{ucp}{\longrightarrow}}\int Z_{j}\left(  \Gamma^{\tau_{K}}\right)
\delta\left(  Y\left[  h^{j}\right]  \right)  .
\end{equation}

\noindent {\bf The expression (2):} We want to study now the $ucp$ convergence of $\frac{1}{s}\int\left(
Z_{j}\left(  \Sigma^{s}\right)  -Z_{j}\left(  \Gamma^{\tau_{K}}\right)
\right)  \delta h^{j}\left(  \Gamma^{\tau_{K}}\right)  $ as $s\rightarrow0$. As in the previous paragraphs, we define
\begin{align*}
y_{n}\left(  s\right)   &  :=\frac{1}{s}\left(  \sum_{i=0}^{k_{n}-1}\frac
{1}{2}\left(  Z_{j}\left(  \Sigma^{s}\right)  _{T_{i+1}^{n}}+Z_{j}\left(
\Sigma^{s}\right)  _{T_{i}^{n}}\right)  \left(  h^{j}\left(  \Gamma^{\tau_{K}%
}\right)  ^{T_{i+1}^{n}}-h^{j}\left(  \Gamma^{\tau_{K}}\right)  ^{T_{i}^{n}%
}\right)  \right. \\
&  -\left.  \sum_{i=0}^{k_{n}-1}\frac{1}{2}\left(  Z_{j}\left(  \Gamma
^{\tau_{K}}\right)  _{T_{i+1}^{n}}+Z_{j}\left(  \Gamma^{\tau_{K}}\right)
_{T_{i}^{n}}\right)  \left(  h^{j}\left(  \Gamma^{\tau_{K}}\right)
^{T_{i+1}^{n}}-h^{j}\left(  \Gamma^{\tau_{K}}\right)  ^{T_{i}^{n}}\right)
\right) \\
&  =\sum_{i=0}^{k_{n}-1}\frac{1}{2}\left(  \frac{Z_{j}\left(  \Sigma
^{s}\right)  _{T_{i+1}^{n}}-Z_{j}\left(  \Gamma^{\tau_{K}}\right)
_{T_{i+1}^{n}}}{s}+\frac{Z_{j}\left(  \Sigma^{s}\right)  _{T_{i}^{n}}%
-Z_{j}\left(  \Gamma^{\tau_{K}}\right)  _{T_{i}^{n}}}{s}\right)  \left(
h^{j}\left(  \Gamma^{\tau_{K}}\right)  ^{T_{i+1}^{n}}-h^{j}\left(
\Gamma^{\tau_{K}}\right)  ^{T_{i}^{n}}\right)
\end{align*}
as a discretization of the Stratonovich integral $\frac{1}{s}\int\left(
Z_{j}\left(  \Sigma^{s}\right)  -Z_{j}\left(  \Gamma^{\tau_{K}}\right)
\right)  \delta h^{j}\left(  \Gamma^{\tau_{K}}\right)  $ using
$\sigma_{n}$. Then, by construction,
\[
y_{n}\left(  s\right)  \underset{n\rightarrow\infty}{\underset{ucp}%
{\longrightarrow}}\frac{1}{s}\int\left(  Z_{j}\left(  \Sigma^{s}\right)
-Z_{j}\left(  \Gamma^{\tau_{K}}\right)  \right)  \delta h^{j}\left(
\Gamma^{\tau_{K}}\right)  .
\]
On the other hand, invoking Definition~\ref{def pathwise variation} and
Lemma~\ref{ucp implica probabilidad} we have that
\begin{align*}
\frac{Z_{j}\left(  \Sigma^{s}\right)  _{T_{i+1}^{n}}-Z_{j}\left(  \Gamma
^{\tau_{K}}\right)  _{T_{i+1}^{n}}}{s}  &  =\left(  \frac{Z_{j}\left(
\Sigma^{s}\right)  -Z_{j}\left(  \Gamma^{\tau_{K}}\right)  }{s}\right)
_{T_{i+1}^{n}}\underset{s\rightarrow0}{\underset{ucp}{\longrightarrow}%
}Y\left[  Z_{j}\right]  _{T_{i+1}^{n}}\\
\frac{Z_{j}\left(  \Sigma^{s}\right)  _{T_{i}^{n}}-Z_{j}\left(  \Gamma
^{\tau_{K}}\right)  _{T_{i}^{n}}}{s}  &  =\left(  \frac{Z_{j}\left(
\Sigma^{s}\right)  -Z_{j}\left(  \Gamma^{\tau_{K}}\right)  }{s}\right)
_{T_{i}^{n}}\underset{s\rightarrow0}{\underset{ucp}{\longrightarrow}}Y\left[
Z_{j}\right]  _{T_{i}^{n}}.
\end{align*}
We now use again the boundedness of $\Sigma $ and $Y $ to guarantee the boundedness of   $Y\left[  Z_{j}\right]  _{T_{i+1}^{n}}=\left( {\bf  i}_{Y}\mathbf{d}
Z_{j}\right)  _{T_{i+1}^{n}}$ and $Y\left[  Z_{j}\right]  _{T_{i}^{n}}=\left(
{\bf i}_{Y}\mathbf{d}Z_{j}\right)  _{T_{i}^{n}}$ (notice that $\mathbf{d}Z_{j}$ is only evaluated on the compact set $K$
because $Y$ is a vector field over
$\Gamma^{\tau_{K}}\subseteq K$). Therefore,  by Lemma~\ref{producto ucp convergencia},
\[
x_{n}\left(  s\right)  \underset{s\rightarrow0}{\underset{ucp}{\longrightarrow
}}\sum_{i=0}^{k_{n}-1}\frac{1}{2}\left(  Y\left[  Z_{j}\right]  _{T_{i+1}^{n}%
}+Y\left[  Z_{j}\right]  _{T_{i}^{n}}\right)  \left(  h^{j}\left(
\Gamma^{\tau_{K}}\right)  ^{T_{i+1}^{n}}-h^{j}\left(  \Gamma^{\tau_{K}%
}\right)  ^{T_{i}^{n}}\right)  :=x_{n}^{\ast}.
\]
Additionally, the sequence $\{x_{n}^{\ast}\}_{n \in \mathbb{N}}$ obviously converge
in $ucp$ to $\int Y\left[  Z_{j}\right]  \delta\left(  h^{j}\left(
\Gamma^{\tau_{K}}\right)  \right)  $ as $n\rightarrow\infty$. Hence, by Lemma
\ref{lemma serie convergence}, we conclude that
\begin{equation}
\label{we will substitute 2}
\frac{1}{s}\left[  \int\left(  Z_{j}\left(  \Sigma^{s}\right)  -Z_{j}\left(
\Gamma^{\tau_{K}}\right)  \right)  \delta h^{j}\left(  \Gamma^{\tau_{K}%
}\right)  \right]  \underset{s\rightarrow0}{\underset{ucp}{\longrightarrow}%
}\int Y\left[  Z_{j}\right]  \delta\left(  h^{j}\left(  \Gamma^{\tau_{K}%
}\right)  \right)  .
\end{equation}

To sum up, if we substitute~(\ref{we will substitute 1}) and~(\ref{we will substitute 2}) in~(\ref{eq 3}) we obtain that
\[
\frac{1}{s}\left[  \int\left\langle \alpha,\delta\Sigma^{s}\right\rangle
-\int\left\langle \alpha,\delta\Gamma^{\tau_{K}}\right\rangle \right]
\underset{s\rightarrow0}{\underset{ucp}{\longrightarrow}}\sum_{j=1}^{p}\int
Z_{j}\left(  \Gamma^{\tau_{K}}\right)  \delta\left(  Y\left[  h^{j}\right]
\right)  +\int Y\left[  Z_{j}\right]  \delta\left(  h^{j}\left(  \Gamma
^{\tau_{K}}\right)  \right)  .
\]
Using the integration by parts formula,%
\begin{align*}
\int Z_{j}\left(  \Gamma^{\tau_{K}}\right)  \delta\left(  Y\left[
h^{j}\right]  \right)   &  =Z_{j}\left(  \Gamma^{\tau_{K}}\right)  Y\left[
h^{j}\right]  -\left(  Z_{j}\left(  \Gamma^{\tau_{K}}\right)  Y\left[
h^{j}\right]  \right)  _{t=0}-\int Y\left[  h^{j}\right]  \delta\left(
Z_{j}\left(  \Gamma^{\tau_{K}}\right)  \right) \\
&  =\left\langle \alpha\left(  \Gamma^{\tau_{K}}\right)  ,Y\right\rangle
-\left\langle \alpha\left(  \Gamma^{\tau_{K}}\right)  ,Y\right\rangle
_{t=0}-\int Y\left[  h^{j}\right]  \delta\left(  Z_{j}\left(  \Gamma^{\tau
_{K}}\right)  \right)
\end{align*}
and, consequently,%
\[%
\begin{array}
[c]{c}%
\frac{1}{s}\left[  \int\left\langle \alpha,\delta\Sigma^{s}\right\rangle
-\int\left\langle \alpha,\delta\Gamma^{\tau_{K}}\right\rangle \right]
\underset{s\rightarrow0}{\underset{ucp}{\longrightarrow}}\int Y\left[
Z_{j}\right]  \delta\left(  h^{j}\left(  \Gamma^{\tau_{K}}\right)  \right)
-\int Y\left[  h^{j}\right]  \delta\left(  Z_{j}\left(  \Gamma^{\tau_{K}%
}\right)  \right) \\
+\left\langle \alpha\left(  \Gamma^{\tau_{K}}\right)  ,Y\right\rangle
-\left\langle \alpha\left(  \Gamma^{\tau_{K}}\right)  ,Y\right\rangle _{t=0}.
\end{array}
\]
In order to conclude the proof, we claim that
\begin{equation}
\int Y\left[  Z_{j}\right]  \delta\left(  h^{j}\left(  \Gamma^{\tau_{K}
}\right)  \right)  -\int Y\left[  h^{j}\right]  \delta\left(  Z_{j}\left(
\Gamma^{\tau_{K}}\right)  \right)  =\int\left\langle {\bf i}_{Y}\mathbf{d}%
\alpha,\delta\Gamma^{\tau_{K}}\right\rangle . \label{eq 4}%
\end{equation}
Indeed,
\[
\mathbf{d}\alpha=\mathbf{d}\left(  \sum_{j=1}^{p}Z_{j}\mathbf{d}h^{j}\right)
=\sum_{j=1}^{p}\mathbf{d}Z_{j}\wedge\mathbf{d}h^{j}, \quad\text{and}\quad
{\bf i}_{Y}\mathbf{d}\alpha=\sum_{j=1}^{p}\left(  Y\left[  Z_{j}\right]
\mathbf{d}h^{j}-Y\left[  h^{j}\right]  \mathbf{d}Z_{j}\right)
\]
which proofs~(\ref{eq 4}), as required.
\quad $\blacksquare$

\subsection{Auxiliary results about integrals and stopping times} 
In the following paragraphs 
we collect three results that are used in the paper in relation with the interplay
between stopping times and integration limits.

\begin{proposition}
\label{lemma restringir intervalo temporal}
Let $X$ be a continuous
semimartingale defined on $\left[  0,\zeta_{X}\right)  $ and $\Gamma$ a continuous
semimartingale. Let $\tau,~\xi$ be two stopping times such that $\tau\leq
\xi<\zeta_{X}.$ Then,%
\[
\left(  X\cdot \Gamma\right)  ^{\tau}=\left(  \mathbf{1}_{\left[  0,\tau\right]
}X\right)  \cdot \Gamma=\left(  X\cdot \Gamma^{\tau}\right)  \quad \mbox{ and } \quad
\left(  X\cdot \Gamma\right)  ^{\xi}-\left(  X\cdot \Gamma\right)  ^{\tau}=\left(
\mathbf{1}_{\left(  \tau,\xi\right]  }X\right)  \cdot \Gamma\]
An equivalent result holds when dealing with the Stratonovich
integral, namely
\begin{equation*}
\left(  \int X\delta \Gamma\right)  ^{\tau}    =\int X\delta
\Gamma^{\tau}=\left(  \int X ^\tau\delta \Gamma\right)  ^{\tau} .
\end{equation*}

\end{proposition}

\noindent\textbf{Proof.\ \ }
By~\cite[Theorem 12, page 60]{protter} we have that
$
\mathbf{1}_{\left[  0,\tau\right]  }X\cdot \Gamma=\left(  X\cdot \Gamma\right)  ^{\tau
}=\left(  X\cdot \Gamma^{\tau}\right)$.
Therefore,%
\[
\left(  X\cdot \Gamma\right)  ^{\xi}-\left(  X\cdot \Gamma\right)  ^{\tau}%
=\mathbf{1}_{\left[  0,\xi\right]  }X\cdot \Gamma-\mathbf{1}_{\left[
0,\tau\right]  }X\cdot \Gamma=\left[  \left(  \mathbf{1}_{\left[  0,\xi\right]
}-\mathbf{1}_{\left[  0,\tau\right]  }\right)  X\right]  \cdot \Gamma=\left(
\mathbf{1}_{\left(  \tau,\xi\right]  }X\right)  \cdot \Gamma.
\]
As to the Stratonovich integral, since $X$ and $\Gamma $
are semimartingales, we can write~\cite[Theorem 23, page 68]{protter} that 
\begin{equation*}
\left(  \int X\delta \Gamma\right)  ^{\tau}=\left(  X\cdot \Gamma\right)  ^{\tau
}+\frac{1}{2}\left[  X,\Gamma\right]  ^{\tau}=\left(  X\cdot \Gamma^{\tau}\right) 
+\frac{1}{2}\left[  X,\Gamma^{\tau}\right] =\int X\delta \Gamma^{\tau}. 
\end{equation*}
Finally, observe that for any process, $\left(  X^{\tau}\right)  ^{\tau
}=X^{\tau}$. On the other hand, taking into account that $\mathbf{1}_{\left[
0,\tau\right]  }X=\mathbf{1}_{\left[  0,\tau\right]  }X^{\tau}$ and $\left[
\Gamma,X\right]  =\left[  X,\Gamma\right]  $, we have%
\begin{align*}
\left(  \int X\delta \Gamma\right)  ^{\tau}  &  =\mathbf{1}_{\left[  0,\tau\right]
}X\cdot \Gamma+\frac{1}{2}\left[  X,\Gamma\right]  ^{\tau}=\mathbf{1}_{\left[
0,\tau\right]  }X^{\tau}\cdot \Gamma+\left(  \frac{1}{2}\left[  X,\Gamma\right]  ^{\tau
}\right)  ^{\tau}\\
&  =\left(  X^{\tau}\cdot \Gamma\right)  ^{\tau}+\left(  \frac{1}{2}\left[
X^{\tau},\Gamma\right]  \right)  ^{\tau}=\left(  X^{\tau}\cdot \Gamma+\frac{1}{2}\left[
X^{\tau},\Gamma\right]  \right)  ^{\tau}=\left(  \int X^{\tau}\delta \Gamma\right)
^{\tau}. \quad \blacksquare
\end{align*}

\begin{proposition}
\label{tau series}
Let $X:\mathbb{R}_{+}\times\Omega\rightarrow\mathbb{R}$ be a real valued
process. Let $\left\{  \tau_{n}\right\}  _{n\in\mathbb{N}}$ be a sequence of
stopping times such that a.s. $\tau_{0}=0$, $\tau_{n}\leq\tau_{n+1}$, for all
$n\in\mathbb{N}$, and $\sup_{n\in\mathbb{N}}\tau_{n}=\infty$. Then,%
\[
X=\underset{n\rightarrow\infty}{\lim_{ucp}}X^{\tau_{n}}.
\]
In particular, if $\Gamma:\mathbb{R}_{+}\times\Omega\rightarrow M$ is a continuous
$M$-valued semimartingale and $\eta\in\Omega_{2}\left(  M\right)  $ then, 
\[
\int\left\langle \eta,d\Gamma\right\rangle =\underset{k\rightarrow\infty}%
{\lim_{ucp}}\left(  \int\left\langle \eta,d\Gamma\right\rangle \right)
^{\tau_{k}}=\underset{k\rightarrow\infty}{\lim_{ucp}}\sum_{n=0}^{k-1}%
\int\mathbf{1}_{\left(  \tau_{n},\tau_{n+1}\right]  }\left\langle \eta
,d\Gamma\right\rangle .
\]
\end{proposition}

\noindent\textbf{Proof.\ \ } Let $\varepsilon>0$ and $t\in\mathbb{R}_{+}$. Then for any $s\in\left[
0,t\right]  $ one has%
\[
\left\{  \left\vert X^{\tau_{n}}-X\right\vert _{s}>\varepsilon\right\}
\subseteq\left\{  \tau_{n}<s\right\}  \subseteq\left\{  \tau_{n}<t\right\}  .
\]
Hence for any $t\in\mathbb{R}_{+}$%
\[
P\left(  \left\{  \left\vert X^{\tau_{n}}-X\right\vert _{s}>\varepsilon
\right\}  \right)  \leq P\left(  \left\{  \tau_{n}<t\right\}  \right)  .
\]
The result follows because $P\left(  \left\{  \tau_{n}<t\right\}  \right)
\rightarrow0$ as $n\rightarrow\infty$ since $\tau_{n}\rightarrow\infty$ a.s.,
and hence in probability.

Let now $\Gamma$ be a $M$-valued continuous semimartingale and $\eta\in
\Omega_{2}\left(  M\right)  $. Notice first that $\left(
{\displaystyle\int}
\left\langle \eta,d\Gamma\right\rangle \right)  ^{\tau_{0}}=0$ because
$\tau_{0}=0$. Consequently, by Proposition~\ref{lemma restringir intervalo temporal} we can write%
\[
\left(  \int\left\langle \eta,d\Gamma\right\rangle \right)  ^{\tau_{k}}%
=\sum_{n=0}^{k-1}\left(  \int\left\langle \eta,d\Gamma\right\rangle \right)
^{\tau_{n+1}}-\left(  \int\left\langle \eta,d\Gamma\right\rangle \right)
^{\tau_{n}}=\sum_{n=0}^{k-1}\int\mathbf{1}_{\left(  \tau_{n},\tau
_{n+1}\right]  }\left\langle \eta,d\Gamma\right\rangle
\]
and the result follows. \quad $\blacksquare$

\begin{proposition}
\label{with indicator integral}
Let $X$ and $Y$ be two real semimartingales.
Suppose that $X$ is continuous and $X_{0}=0$. Then, for any $t\in
\mathbb{R}_{+}$, the Stratonovich integral $\int\left(  \mathbf{1}_{\left[
0,t\right]  }Y\right)  \delta X$ is well defined and equal to $\left(  \int Y\delta
X\right)  ^{t}$.
\end{proposition}

\noindent\textbf{Proof.\ \ }
If $\int\left(  \mathbf{1}_{\left[  0,t\right]  }Y\right)  \delta X$ was well
defined, it should be equal to $\int\left(  \mathbf{1}_{\left[  0,t\right]
}Y\right)  dX+\frac{1}{2}\left[  \mathbf{1}_{\left[  0,t\right]  }Y,X\right]
$. Since $\int\left(  \mathbf{1}_{\left[  0,t\right]  }Y\right)  dX$ is
well defined, the only thing that we need to check is that $\left[  \mathbf{1}_{\left[
0,t\right]  }Y,X\right]  $ exists. On the other hand, recall that
(\cite[Theorem 12 page 60 and Theorem 23 page 68]{protter})%
\[
\left(  \int Y\delta X\right)  ^{t}=\int\left(  \mathbf{1}_{\left[
0,t\right]  }Y\right)  dX+\frac{1}{2}\left[  Y,X\right]  ^{t}=\int\left(
\mathbf{1}_{\left[  0,t\right]  }Y\right)  dX+\frac{1}{2}\left[
Y^{t},X\right]  .
\]
Hence, what we are actually going to proceed by showing that $\left[
\mathbf{1}_{\left[  0,t\right]  }Y,X\right]  $ is equal to $\left[
Y^{t},X\right]  $.

Let $\sigma_{n}=\left\{  0=T_{0}^{n}\leq T_{1}^{n}\leq...\leq T_{k_{n}}
^{n}<\infty\right\}  $ be a sequence of random partitions tending to the identity (in the sense of~\cite[page 64]{protter}).  Given two real processes $X$ and $Y$, their quadratic
variation, if it exists, can be defined as the limit in $ucp$ when $n \rightarrow  \infty $ of the following sums
\[
\left[  Y,X\right]  =\underset{n\rightarrow\infty}{\lim_{ucp}}\sum
_{i=0}^{k_{n}-1}\left(  Y^{T_{i+1}^{n}}-Y^{T_{i}^{n}}\right)  \left(
X^{T_{i+1}^{n}}-X^{T_{i}^{n}}\right)  .
\]
Let now
\begin{align*}
H_{n}  &  :=\sum_{i=0}^{k_{n}-1}\left(  \left(  Y^{t}\right)  ^{T_{i+1}^{n}%
}-\left(  Y^{t}\right)  ^{T_{i}^{n}}\right)  \left(  X^{T_{i+1}^{n}}%
-X^{T_{i}^{n}}\right)  ,\\
G_{n}  &  :=\sum_{i=0}^{k_{n}-1}\left(  \left(  \mathbf{1}_{\left[
0,t\right]  }Y\right)  ^{T_{i+1}^{n}}-\left(  \mathbf{1}_{\left[  0,t\right]
}Y^{t}\right)  ^{T_{i}^{n}}\right)  \left(  X^{T_{i+1}^{n}}-X^{T_{i}^{n}%
}\right)  .
\end{align*}
It is clear that the sequence $\{H_{n}\}_{n \in  \mathbb{N}}$ converges uniformly on compacts in probability to
$\left[  Y^{t},X\right]  $. We are going to prove that there exists such a convergence for
the sequence of processes $\left\{  G_{n}\right\}  _{n\in\mathbb{N}}$ by showing  that the elements $\left(  G_{n}\right)  _{s}$ coincide with $\left(
H_{n}\right)  _{s}$,  for any $s\in\mathbb{R}_{+}$, up to a set whose probability tends to zero as $n\rightarrow\infty$. We will consider two cases:

\medskip

\noindent {\bf 1. The case  $s\leq t$.}
Given a specific $i\in\left\{  0,...,k_{n}-1\right\}  $, and recalling that by construction
$T_{i}^{n}\leq T_{i+1}^{n}$ a.s., it is clear that $\left(  \left(  Y^{t}\right)  ^{T_{i+1}
^{n}}-\left(  Y^{t}\right)  ^{T_{i}^{n}}\right)  _{s}=Y_{T_{i+1}^{n}\wedge
s}-Y_{T_{i}^{n}\wedge s}$ is different from $0$ only for those $\omega
\in\Omega$ in $\left\{  T_{i}^{n}<s\right\}  $ in which case it takes the value 
\begin{equation}
Y_{T_{i+1}^{n}\wedge s}-Y_{T_{i}^{n}}. \label{eq 11}%
\end{equation}
On the other hand, $\left(  \left(  \mathbf{1}_{\left[  0,t\right]  }Y\right)
^{T_{i+1}^{n}}-\left(  \mathbf{1}_{\left[  0,t\right]  }Y^{t}\right)
^{T_{i}^{n}}\right)  _{s}$ is again different from $0$ only in the set
$\left\{  T_{i}^{n}<s\right\}  $ and there it is equal to (\ref{eq 11}).
Therefore, $\left(  G_{n}\right)  _{s}=\left(  H_{n}\right)  _{s}$ whenever $s\leq
t$.

\medskip

\noindent {\bf 2. The case $s>t$.}
In this case, $\left(  \left(  Y^{t}\right)  ^{T_{i+1}^{n}}-\left(
Y^{t}\right)  ^{T_{i}^{n}}\right)  _{s}=Y_{t\wedge T_{i+1}^{n}}-Y_{t\wedge
T_{i}^{n}}$ which is different from $0$ only in the set $\left\{  T_{i}%
^{n}<t\right\}  $, where it takes the value
\begin{equation}
Y_{t\wedge T_{i+1}^{n}}-Y_{T_{i}^{n}} \label{eq 12}
\end{equation}
However, in this case $\left(  \left(  \mathbf{1}_{\left[  0,t\right]
}Y\right)  ^{T_{i+1}^{n}}-\left(  \mathbf{1}_{\left[  0,t\right]  }%
Y^{t}\right)  ^{T_{i}^{n}}\right)  _{s}=\mathbf{1}_{\left\{  T_{i+1}^{n}\leq
t\right\}  }Y_{t\wedge T_{i+1}^{n}}-\mathbf{1}_{\left\{  T_{i}^{n}\leq
t\right\}  }Y_{t\wedge T_{i}^{n}}$, which is equal to (\ref{eq 12}) in the set
$\left\{  T_{i+1}^{n}\leq t\right\}  $ (which contains $\left\{  T_{i}
^{n}<t\right\}  $ since $T_{i}^{n}\leq T_{i+1}^{n}$), but differs from
(\ref{eq 12}) in
\[
A_{i}^{n}\left(  t\right)  :=\left\{  T_{i}^{n}\leq t<T_{i+1}^{n}\right\}
\]
where it takes the value $-Y_{T_{i}^{n}}$. For any other $\omega\in\Omega$ not in these
sets, $\left(  \left(  \mathbf{1}_{\left[  0,t\right]  }Y\right)
^{T_{i+1}^{n}}-\left(  \mathbf{1}_{\left[  0,t\right]  }Y^{t}\right)
^{T_{i}^{n}}\right)  _{s}\left(  \omega\right)  =0$. Therefore, whenever $s>t$,
$\left(  G_{n}\right)  _{s}$ and $\left(  H_{n}\right)  _{s}$ are different only for the
$\omega\in A_{i}^{n}\left(  t\right)  $. Observe that, since $t$ is
fixed, only one of the sets $\left\{  A_{i}^{n}\left(  t\right)  \right\}
_{i\in \{0, \ldots,  k_{n}-1\}}$ is non-empty and, on it,
\[
\left(  H_{n}\right)  _{s}-\left(  G_{n}\right)  _{s}=Y_{t}\left(
X_{t}-X_{T_{i}^{n}}\right)  .
\]

\medskip

To sum up, the analysis that we just carried out shows that for any  $u\in\mathbb{R}_{+}$
\[
\sup_{0\leq s\leq u}\left\vert \left(  H_{n}\right)  _{s}-\left(
G_{n}\right)  _{s}\right\vert =\mathbf{1}_{A_{i}^{n}\left(  t\right)
}\left\vert Y_{t}\right\vert \left\vert \left(  X_{t}-X_{T_{i}^{n}}\right)
\right\vert
\]
for some $i\in\left\{  0,...,k_{n}-1\right\}  $. If $X$ is continuous, this
expression tells us that $\sup_{0\leq s\leq u}\left\vert \left(  H_{n}\right)
_{s}-\left(  G_{n}\right)  _{s}\right\vert \rightarrow0$ a.s. as
$n\rightarrow\infty$ which, in turn, implies that $\sup_{0\leq s\leq
u}\left\vert \left(  H_{n}\right)  _{s}-\left(  G_{n}\right)  _{s}\right\vert
$ converges to $0$ in probability as well. That is, for any $\varepsilon>0$,%
\[
P\left(  \left\{  \sup_{0\leq s\leq u}\left\vert \left(  H_{n}\right)
_{s}-\left(  G_{n}\right)  _{s}\right\vert >\varepsilon\right\}  \right)
\rightarrow0,\text{ \ as }n\rightarrow\infty,
\]
which is the same as saying that $H_{n}-G_{n}$ converges to $0$ in $ucp$.
Thus, since $G_{n}=H_{n}-\left(  H_{n}-G_{n}\right)  $ and the limit in $ucp$
as $n\rightarrow\infty$ exist for the both sequences $\left\{  H_{n}\right\}
_{n\in\mathbb{N}}$ and $\left\{  H_{n}-G_{n}\right\}  _{n\in\mathbb{N}}$, so
does the limit of $\left\{  G_{n}\right\}  _{n\in\mathbb{N}}$  which, by
definition, is the quadratic variation $\left[  \mathbf{1}_{\left[
0,t\right]  }Y,X\right]  $. Moreover, as $\left(  H_{n}-G_{n}\right)
\rightarrow0$ in $ucp$ as $n\rightarrow\infty$,
\[
\left[  Y^{t},X\right]  =\underset{n\rightarrow\infty}{\lim_{ucp}}
H_{n}=\underset{n\rightarrow\infty}{\lim_{ucp}}G_{n}=\left[  \mathbf{1}
_{\left[  0,t\right]  }Y,X\right],
\]
which concludes the proof.
\quad $\blacksquare$

\section{Appendices}

\subsection{Preliminaries on semimartingales and integration}

In the following paragraphs we state a few standard definitions and results on
manifold valued semimartingales and integration. Semimartingales are the
natural setup for stochastic differential equations and, in particular, for
the equations that we handle in this paper. For proofs and additional details
the reader is encouraged to check, for instance, with~\cite{chung williams,
durrett, emery, ikeda watanabe, legall integration, protter}, and references therein.

\medskip

\noindent\textbf{Semimartingales.} The first element in our setup for
stochastic processes is a probability space $(\Omega,\mathcal{F}, P)$ together
with a filtration $\{ \mathcal{F} _{t}\mid t \geq0\} $ of $\mathcal{F} $ such
that $\mathcal{F} _{0}$ contains all the negligible events (complete
filtration) and the map $t \longmapsto\mathcal{F}_{t} $ is right-continuous,
that is, $\mathcal{F} _{t}=\bigcap_{\epsilon>0} \mathcal{F} _{t+ \epsilon}$.

A real-valued {\bfseries\itshape  martingale} $\Gamma: \mathbb{R}_{+}=[0, \infty)
\times\Omega\rightarrow\mathbb{R} $ is a stochastic process such that for
every pair $t,s \in\mathbb{R}_{+} $ such that $s \leq t $, we have:

\begin{description}
\item[(i)] $\Gamma$ is $\mathcal{F} _{t} $-adapted, that is, $\Gamma_{t} $ is
$\mathcal{F}_{t}$-measurable.

\item[(ii)] $\Gamma_{s}=E[\Gamma_{t}\mid\mathcal{F} _{s}]$.

\item[(iii)] $\Gamma_{t}$ is integrable: $E[|\Gamma_{t} |] <+ \infty$.
\end{description}

For any $p \in[1, \infty) $, $\Gamma$ is called a $L ^{p} $%
{\bfseries\itshape -martingale} whenever $\Gamma$ is a martingale and
$\Gamma_{t} \in L ^{p}(\Omega) $, for each $t$. If $\sup_{t \in\mathbb{R}_{+}}
\mathrm{E}[|\Gamma_{t}|^{p}]< \infty$, we say that $\Gamma$ is $L ^{p}
${\bfseries\itshape
-bounded}. The process $\Gamma$ is {\bfseries\itshape  locally bounded } if
for any time $t\geq0$, $\sup\{| \Gamma_{s}(\omega)|\mid s\leq t\}< \infty$, almost surely.
Every continuous process is locally bounded. Recall that a process is said to
be {\bfseries\itshape  continuous} when its paths are continuous. Most
processes considered in this paper will be of this kind. Given two continuous
processes $X $ and $Y $ we will write $X=Y $ when they are a modification of
each other or when they are indistinguishable since these two concepts
coincide for continuous processes.

A random variable $\tau:\Omega\rightarrow\lbrack0,+\infty]$ is called a
{\bfseries\itshape stopping time} with respect to the filtration
$\{\mathcal{F}_{t}\mid t\geq0\}$ if for every $t\geq0$ the set $\{\omega
\mid\tau(\omega)\leq t\}$ belongs to $\mathcal{F}_{t}$. Given a stopping time
$\tau,$ we define%
\[
\mathcal{F}_{\tau}=\left\{  \Lambda\in\mathcal{F}~|~\Lambda\cap\left\{
\tau\leq t\right\}  \in\mathcal{F}_{t}\text{ for any }t\in\mathbb{R}%
_{+}\right\}  .
\]
Given an adapted process $\Gamma,$ it can be shown that $\Gamma_{\tau}$ is
$\mathcal{F}_{\tau}$-measurable. Furthermore, the {\bfseries\itshape stopped
process} $\Gamma^{\tau}$ is defined as
\[
\Gamma_{t}^{\tau}:=\Gamma_{t\wedge\tau}:=\Gamma_{t}\boldsymbol{1}_{\{t\leq
\tau\}}+\Gamma_{\tau}\boldsymbol{1}_{\{t>\tau\}}.
\]

A {\bfseries\itshape  continuous local martingale} is a continuous adapted
process $\Gamma$ such that for any $n \in\mathbb{N} $, $\Gamma^{\tau_{n}}
\boldsymbol{1}_{\{\tau_{n}>0\}}$ is a martingale, where $\tau_{n}
$ is the stopping time $\tau_{n}:=\inf\{t \geq0\mid| \Gamma_{t}|= n \} $.

We say that the stochastic process $\Gamma: \mathbb{R}_{+} \times
\Omega\rightarrow\mathbb{R} $ has {\bfseries\itshape  finite variation}
whenever it is adapted and has bounded variation on compact subintervals of
$\mathbb{R}_{+} $. This means that for each fixed $\omega\in\Omega$, the path
$t \longmapsto\Gamma_{t}(\omega)$ has bounded variation on compact
subintervals of $\mathbb{R}_{+} $, that is, the supremum $\sup\left\{
\sum_{i=1}^{p} | \Gamma_{t _{i}}( \omega)-\Gamma_{t _{i-1}}( \omega)|\right\}
$ over all the partitions $0=t _{0}< t _{1}< \cdots< t _{p}= t $ of the
interval $[0,t] $ is finite.

A {\bfseries\itshape continuous semimartingale} is the sum of a continuous
local martingale and a process with finite variation. It can be proved that a
given semimartingale has a unique decomposition of the form $\Gamma=
\Gamma_{0}+V+ \Lambda$, with $\Gamma_{0}$ the initial value of $\Gamma$, $V$ a
finite variation process, and $\Lambda$ a local continuous semimartingale.
Both $V$ and $\Lambda$ are null at zero.

\medskip

\noindent\textbf{The It\^o integral with respect to a continuous
semimartingale.} Let $\Gamma: \mathbb{R}_{+} \times\Omega\rightarrow\mathbb{R}
$ be a continuous local martingale. It can be shown that there exists a unique
increasing process with finite variation $[ \Gamma, \Gamma] _{t}$ such that
$\Gamma_{t} ^{2}- [ \Gamma, \Gamma] _{t}$ is a local continuous martingale. We
will refer to $[ \Gamma, \Gamma] _{t}$ as the {\bfseries\itshape  quadratic
variation} of $\Gamma$. Given $\Gamma= \Gamma_{0}+V+ \Lambda, \Gamma^{\prime}=
\Gamma_{0}^{\prime}+V^{\prime}+ \Lambda^{\prime}$ two continuous local
martingales we define their {\bfseries\itshape  joint quadratic variation} or {\bfi 
quadratic covariation} as
\[
[\Gamma, \Gamma^{\prime}]_{t}=\frac{1}{2} \left(  [\Lambda+ \Lambda^{\prime},
\Lambda+ \Lambda^{\prime}]_{t}-[\Lambda, \Lambda]_{t}-[\Lambda^{\prime},
\Lambda^{\prime}]_{t} \right)  .
\]

Let $\left\{  X_{n}\right\}  _{n\in\mathbb{N}}$ be a sequence of processes. We
will say that $\left\{  X_{n}\right\}  _{n\in\mathbb{N}}$ {\bfseries\itshape
converges uniformly on compacts in probability } (abbreviated {\bfi  ucp}) to a process
$X$ if for any
$\varepsilon>0$ and any $t\in\mathbb{R}_{+},$
\[
P\left(  \left\{  \sup_{0\leq s\leq t}\left\vert X_{n}-X\right\vert
_{s}\right\}  >\varepsilon\right)  \longrightarrow0,
\]
as $n\rightarrow\infty$.

Following \cite{protter}, we denote by $\mathbb{L}$ the space of processes
$X:\mathbb{R}_{+}\times\Omega\rightarrow\mathbb{R}$ whose paths are
left-continuous and have right limits. These are usually called
{\bfi  c\`{a}gl\`{a}d} processes, which are initials in French for
\textit{left-continuous with right limits}. We say that a process
$X\in\mathbb{L}$ is {\bfseries\itshape elementary} whenever it can be expressed
as
\[
X=X_{0}\boldsymbol{1}_{\left\{  0\right\}  }+\sum_{i=1}^{p-1}X_{i}%
\boldsymbol{1}_{(\tau_{i},\tau_{i+1}]},
\]
where $0\leq\tau_{1}<\cdots<\tau_{p-1}<\tau_{p}$ are stopping times, and $X_{0}$
and $X_{i}$ are $\mathcal{F}_{0}$ and $\mathcal{F}_{\tau_{i}}
$-measurable random variables, respectively such that $|X_{0}|<\infty$ and
$|X_{i}|<\infty$ a.s. for all $i \in \{1, \ldots,p-1\}$.
$\boldsymbol{1}_{(\tau_{i} ,\tau_{i+1}]}$ is the characteristic function
of the set $(\tau_{i},\tau_{i+1}]=\left\{ 
\left( t,\omega\right)  \in R_{+}\times\Omega~|~t\in(\tau_{i}\left(  \omega\right)
,\tau_{i+1}\left(  \omega\right)  ]\right\}  $ and $\boldsymbol{1}_{\left\{
0\right\}  }$ of $\left\{  \left(  t,\omega\right)  \in R_{+}\times
\Omega~|~t=0\right\}  .$ It can be shown (see \cite[Theorem 10, page 57]{protter})
that the set of elementary processes is dense in $\mathbb{L}$ in the $ucp$ topology.

Let $\Gamma$ be a semimartingale such that $\Gamma_{0}=0$ and $X $ elementary. We define
{\bfseries\itshape It\^{o}'s stochastic integral} of $X$ with respect to
$\Gamma$ as given by
\begin{equation}
X\cdot\Gamma:=\int Xd\Gamma:=\sum_{i=1}^{p-1}X_{i}(\Gamma^{\tau_{i+1}}%
-\Gamma^{\tau_{i}}). \label{integral Ito definite}%
\end{equation}
In the sequel we will exchangeably use the symbols $X\cdot\Gamma$ and $\int
Xd\Gamma$ to denote the It\^{o} stochastic integral. It is a deep result that,
if $\Gamma$ is a semimartingale, the It\^{o} stochastic integral is a
continuous map from $\mathbb{L}$ into the space of processes whose paths are
right-continuous and have left limits ({\bfi  c\`adl\`ag}), usually denoted by
$\mathbb{D}$, equipped also with the $ucp$ topology. Therefore we can extend the It\^{o}
integral to the whole $\mathbb{L}$. In particular, we can integrate any
continuous adapted processes with respect to any semimartingale.

Given any stopping time $\tau$ we define
\[
\int_{0}^{\tau}Xd\Gamma:=(X\cdot\Gamma)_{\tau}.
\]
It can be shown that $(\boldsymbol{1}_{[0,\tau]}X)\cdot\Gamma=(X\cdot
\Gamma)^{\tau}=X\cdot\Gamma^{\tau}$. If there exists a stopping times
$\zeta_{\Gamma}$ such that the semimartingale $\Gamma$ is defined only on the
stochastic intervals $[0,\zeta_{\Gamma})$, then we may define the It\^{o}
integral of $X$ with respect to $\Gamma$ on any interval $[0,\tau]$ such that
$\tau<\zeta_{\Gamma}$ by means of $X\cdot\Gamma^{\tau}$.

\medskip

\noindent\textbf{The Stratonovich integral and stochastic calculus.} Given
$\Gamma$ and $X$ two semimartingales we define the {\bfseries\itshape
Stratonovich integral} of $X$ along $\Gamma$ as
\[
\int_{0}^{t}X\delta\Gamma=\int_{0}^{t}Xd\Gamma+\frac{1}{2}[X,\Gamma]_{t}.
\]
Let $X^{1},\ldots,X^{p}$ be $p$ continuous semimartingales and $f\in
C^{2}(\mathbb{R}^{p})$. The celebrated {\bfseries\itshape It\^{o} formula}
states that
\begin{align*}
f(X_{t}^{1},\ldots,X_{t}^{p})  &  =f(X_{0}^{1},\ldots,X_{0}^{p})+\sum
_{i=1}^{p}\int_{0}^{t}\frac{\partial f}{\partial x^{i}}(X_{s}^{1},\ldots
,X_{s}^{p})dX_{s}^{i}\\
&  +\frac{1}{2}\sum_{i,j=1}^{p}\int_{0}^{t}\frac{\partial^{2}f}{\partial
x^{i}\partial x^{j}}(X_{s}^{1},\ldots,X_{s}^{p})d[X^{i},X^{j}]_{s}%
\end{align*}
The analogue of this equality for the Stratonovich integral is
\[
f(X_{t}^{1},\ldots,X_{t}^{p})=f(X_{0}^{1},\ldots,X_{0}^{p})+\sum_{i=1}^{p}%
\int_{0}^{t}\frac{\partial f}{\partial x^{i}}(X_{s}^{1},\ldots,X_{s}%
^{p})\delta X_{s}^{i}.
\]
An important particular case of these relations are the {\bfseries\itshape
integration by parts} formulas
\begin{align*}
\int_{0}^{t}Xd\Gamma &  =(X\Gamma)_{t}-(X\Gamma)_{0}-\int_{0}^{t}\Gamma
dX-\frac{1}{2}[X,\Gamma]_{t},\\
\int_{0}^{t}X\delta\Gamma &  =(X\Gamma)_{t}-(X\Gamma)_{0}-\int_{0}^{t}%
\Gamma\delta X.
\end{align*}

\medskip

\noindent\textbf{Stochastic differential equations.} Let $\Gamma=(\Gamma^{1},
\ldots, \Gamma^{p}) $ be $p$ semimartingales with $\Gamma_{0}=0 $ and $f :
\mathbb{R} ^{q} \times\mathbb{R}^{p} \rightarrow\mathbb{R}^{q} $ a smooth
function. A {\bfseries\itshape  solution} of the {\bfseries\itshape  It\^o
stochastic differential equation}
\begin{equation}
\label{sde statement}d X ^{i}=\sum_{j=1}^{p} f ^{i} _{j}(X, \Gamma) d
\Gamma^{j}%
\end{equation}
with initial condition the random vector $X _{0}=(X _{0}^{1}, \ldots, X
_{0}^{q})$ is a stochastic process $X_{t}=(X _{t}^{1}, \ldots, X _{t}^{q})$
such that $X^{i} _{t} -X^{i} _{0}=\sum_{j=1}^{p} \int_{0} ^{t} f ^{i} _{j}(X,
\Gamma) d \Gamma^{j} $. It can be shown~\cite[page 310]{protter} that for any
$x \in\mathbb{R} ^{q} $ there exists a stopping time $\zeta: \mathbb{R} ^{q}
\times\Omega\rightarrow\mathbb{R}_{+} $ and a time-continuous solution $X(t,
\omega, x )$ of~(\ref{sde statement}) with initial condition $x$ and defined
in the time interval $[0, \zeta(x, \omega))$. Additionally, $\limsup_{t
\rightarrow\zeta(x, \omega)}\| X _{t}(\omega)\|= \infty$ a.s. on $\{\zeta<
\infty\}$ and $X$ is smooth on $x$ in the open set $\{x\mid\zeta(x ,\omega)>
t\}$. Finally, the solution $X $ is a semimartingale.

\subsection{Second order vectors and forms}

In the paragraphs that follow we review the basic tools on second order geometry needed in
the definition of the stochastic integral of a form along a manifold valued semimartingale.
The reader interested in the proofs of the statements cited in this section is
encouraged to check with~\cite{emery}, and references therein.  

Let $M$ be a finite dimensional, second-countable, locally compact Hausdorff (and hence
paracompact) manifold. Given $m \in  M $, a {\bfi  tangent vector} at $m$ {\bfi  of order
two} with no constant term is a differential operator
$L:\mathcal{C}^{\infty}\left(  M\right)  \longrightarrow\mathbb{R}$
that satisfies
\[
L\left[  f^{3}\right]  \left(  m\right)  =3f\left(  m\right)  L\left[
f^{2}\right]  \left(  m\right)  -3f^{2}\left(  m\right)  L\left[  f\right]
\left(  m\right)  .
\]
The vector space of tangent vectors of order two at $m$ is denoted as
$\tau_{m}M$. The manifold $\tau M:=\bigcup_{m\in  M}\tau_{m}M$ is referred to as the
{\bfi  second order tangent bundle} of $M$. Notice that the (first order) tangent bundle $TM
$ of $M$ is contained in $\tau M $.  A {\bfi  vector field of order two} is a smooth section of
the bundle
$\tau M\rightarrow M$. We denote the set of vector fields order two by $\mathfrak{X} _2 (M)
$. If
$Y,Z\in\mathfrak{X}(  M)  $ then the product
$Z Y\in  \mathfrak{X} _2 (M)$. Conversely,
every second order vector field $L \in \mathfrak{X} _2 (M)$ can be written as a finite sum
of fields of the form
$ZY$ and $W$, with $Z,Y,W\in\mathfrak{X}(  M)$.

The {\bfi  forms of order two} $\Omega _2 (M) $ are the smooth sections of the cotangent
bundle of order two
$\tau^\ast M :=\bigcup_{m\in  M}\tau^\ast _{m}M$. For any $f, g,h \in  C^\infty(M)$ and 
$L \in  \mathfrak{X} _2 (M) $ we define
$d_2f \in  \Omega _2 (M) $ by $d_2f(L):=L[f]$,  and
$d_{2}f\cdot d_{2}g\in  \Omega _2 (M) $ as
\[
d_{2}f\cdot d_{2}g[L]:=\frac{1}{2}\left(  L\left[  fg\right]  -fL\left[ 
g\right] -gL\left[  f\right]  \right).
\]
It is easy to show that for any $Y,Z, W\in  \mathfrak{X}(M) $,
\[
d_{2}f\cdot d_{2}g\left[  Z  Y\right]    =\frac{1}{2}\left(  Z\left[
f\right]  Y\left[  g\right]  +Z\left[  g\right]  Y\left[  f\right]  \right)
 \quad \mbox{ and } \quad
d_{2}f\cdot d_{2}g\left[  W\right]    =0.
\]
More generally, let $\alpha_m,\beta_m\in T_m ^\ast M$ and
choose $f,g\in\mathcal{C}^{\infty}\left(  M\right)  $  two functions
 such that $\mathbf{d} f(m)=\alpha_m$ and
$ \mathbf{d}g(m)=\beta_m.$ It is easy to check that $(\mathbf{d}f\cdot \mathbf{d}g) (m)$ does
not depend on the particular choice of $f$ and $g$ above and hence we can write
$\alpha_m\cdot\beta_m$ to denote  $(\mathbf{d}f\cdot \mathbf{d}g) (m)$. If  $\alpha,\beta\in
\Omega (M)$  then we can define 
$\alpha\cdot\beta \in \Omega_2(M)$ as $\left(\alpha\cdot\beta\right)(m):= \alpha(m) \cdot
\beta (m) $. This product is commutative and $\mathcal{C}^{\infty}\left(  M\right) 
$-bilinear.
It can be shown that
every second order form can be locally written as a finite sum of forms of the type
$\mathbf{d}f \cdot \mathbf{d}g$ and $d_2h$.

The $d_2 $  operator can also be defined on forms by using a result  (Theorem 7.1
in~\cite{emery}) that claims that there exists a unique linear operator
$d_2 :
\Omega (M)
\rightarrow
\Omega _2 (M)$  characterized by
\[
d_{2}\left(  \mathbf{d}f\right)    =d_{2}f \quad \mbox{ and } \quad
d_{2}\left(  f\alpha\right)    =\mathbf{d}f\cdot\alpha+fd_{2}\alpha.
\]

\subsection{Stochastic integrals of forms along a semimartingale}

Let $M$ be a manifold. A continuous $M$-valued stochastic process $X: \mathbb{R}_+
\times\Omega \rightarrow  M$ is called a {\bfi  continuous $M$-valued semimartingale} if for
each smooth function
$f \in  C^\infty(M)$, the real valued process $f \circ X$ is a
(real-valued) continuous semimartingale. We say that $X$ is {\bfi  locally
bounded} if the sets $\{X _s( \omega)\mid 0\leq s\leq t\}$ are relatively compact in $M$ for each
$t \in \mathbb{R}_+$, a.s.

Let $X$ be a $M$-valued semimartingale and $\theta: \mathbb{R}_+ \times  \Omega \rightarrow
\tau^\ast M$ be a c\`agl\`ad locally bounded process over $X$, that is,
$\pi\circ\theta=X$, where $\pi: \tau ^\ast  M \rightarrow  M $ is the canonical projection.
It can be shown (see~\cite[Theorem 6.24]{emery}) that there exists a unique linear map
$\theta\longmapsto
\int\left\langle \theta,dX\right\rangle $ that associates to each such $\theta$
a continuous real valued semimartingale and that is fully characterized by the following
properties:  for any $f\in\mathcal{C}^{\infty}\left(  M\right)  $ and any locally bounded
c\`agl\`ad real-valued
process
$K$,
\begin{equation}
\label{properties Ito integral manifold}
\int\left\langle d_{2}f \circ X,dX\right\rangle   =f\left(  X\right)
-f\left(  X_{0}\right), \quad\mbox{and} \quad
\int\left\langle K\theta,dX\right\rangle   =\int Kd\left(
\int\left\langle \theta,dX\right\rangle \right).
\end{equation}
The stochastic process 
$\int\left\langle \theta,dX\right\rangle $ will be called the {\bfi 
It\^o integral} of $\theta$ along $X $. If $\alpha \in \Omega _2 (M) $, we
will write in the sequel  the It\^o integral of $\alpha$ along $X$, that is, $\int \langle
\alpha \circ  X,  d X\rangle
$ as $\int \langle \alpha ,  d X\rangle $.

The integral of a $(0,2)$-tensor $b$ on $M$ along $X$ is the image of the unique linear mapping
$b\longmapsto \int b(dX, dX)$ onto the space of real continuous processes with finite variation
that for all $f,g, \in C^\infty(M) $ satisfies
\begin{equation}
\label{integration tensors}
\int (fb)(dX, dX)=\int (f \circ X) d \left(\int b(d X, d X) \right) \quad \mbox{ and } \quad\int
\left( df\otimes dg\right)(dX, dX)=[f \circ X, g \circ  X].
\end{equation}

If $\alpha\in\Omega\left(
M\right)  $ and $X$ is a semimartingale on $M,$ the real semimartingale
$\int\left\langle d_{2}\alpha,dX\right\rangle $ is called the
{\bfi  Stratonovich integral} of $\alpha$ along $X$ and is denoted by
$
\int\left\langle \alpha,\delta X\right\rangle$. This definition can be generalized by taking $\beta  $  a $T ^\ast M $ valued semimartingale over $X$ and by defining the Stratonovich integral as the unique real valued semimartingale that satisfies the properties 
\begin{equation}
\label{properties strato}
\int\left\langle \mathbf{d}f,\delta X\right\rangle   =f\left(  X\right)  -f\left(
X_{0}\right),\quad\mbox{and}
\quad
\int\left\langle Z \beta,\delta X\right\rangle   =\int Z\left(  X\right)
\delta\left(  \int\left\langle \beta,\delta X\right\rangle \right),
\end{equation}
for any
$f\in\mathcal{C}^{\infty}\left(  M\right)$ and any continuous real valued semimartingale $Z $.
Finally, it can be shown that (see~\cite[Proposition 6.31]{emery}) for any $f,g \in
C^\infty(M) $,
\begin{equation}
\label{integral with quadratic variation}
\int\left\langle \mathbf{d}f\cdot \mathbf{d}g,dX\right\rangle =\frac{1}{2}\left[
f\left(  X\right)  ,g\left(  X\right)  \right].
\end{equation}

\subsection{Stochastic differential equations on manifolds}
\label{Stochastic differential equations on manifolds}

The reader interested in the details of the material presented in this section is
encouraged to check with the chapter 7 in~\cite{emery}. 

Let
$M$ and
$N$ be two manifolds. A {\bfi  Stratonovich operator} from $M$ to $N$ is a family $\{
e(x,y)\}_{x \in M, y \in  N }$ such that $e(x,y) :T _xM \rightarrow T
_yN$ is a linear mapping that depends smoothly on its two entries. Let $e ^\ast (x,y): T
^\ast _y N \rightarrow T _x ^\ast  M $ be the adjoint of $e(x,y)$. 

Let $X$  be a $M$-valued
semimartingale. We say that a
$N$-valued semimartingale is a solution of the the {\bfi  Stratonovich stochastic
differential equation}
\begin{equation}
\label{sde stocastic stratonovich}
\delta Y= e(X,Y) \delta X
\end{equation}
if for any $\alpha \in  \Omega (N) $, the following equality between Stratonovich integrals
holds:
\[
\int   \langle \alpha, \delta Y\rangle =\int  \langle  e ^\ast (X,Y) \alpha,
\delta X\rangle.
\]
It can be shown~\cite[Theorem 7.21]{emery} that given a semimartingale $X$ in $M$, a
$\mathcal{F} _0$ measurable random variable $Y _0$, and a Stratonovich operator $e $ from $M
$ to
$N$, there are a  stopping time $\zeta>0 $ and
a solution
$Y 
$ of~(\ref{sde stocastic stratonovich}) with initial condition $Y _0$ defined on the set
$\{(t, \omega)\in 
\mathbb{R}_+ \times  \Omega\mid t \in  [0, \zeta (\omega))\}$ that has the following
maximality and uniqueness property: if   $\zeta' $ is another
stopping time such that  $\zeta'< \zeta $ and $Y' $ is another solution defined on  $\{(t,
\omega)\in 
\mathbb{R}_+ \times  \Omega\mid t \in  [0, \zeta' (\omega))\}$, then $Y'$ and $Y $ coincide
in this set. If $\zeta $ is finite then $Y $ explodes at time $\zeta$, that is, the path $Y _t $
with $t \in  [0, \zeta) $ is not contained in any compact subset of $N$.

The stochastic differential equations from the It\^o integration point of view require the
notion of {\bfi  Schwartz operator} whose construction we briefly review. The reader
interested in the details of this construction is encouraged to check with~\cite{emery}.
Note first that we can associate to any element
$L \in\mathfrak{X}_2(M)$  a symmetric tensor $
\widehat{L}\in\mathfrak{X} (M) \otimes \mathfrak{X} (M) $. Second, given $x \in M $ and $y
\in N $, a linear mapping from $\tau _x M $  into $\tau _y N $ is called a {\bfi  Schwartz
morphism} whenever $f(T _xM)\subset T _yN $ and $\widehat{f(L)}= \left(f|_{T _x M}\otimes
f|_{T _x M}\right)(\widehat{L}) $, for any $L \in  \tau _x M $. Third, let
$M$ and
$N$ be two manifolds; a {\bfi  Schwartz operator} from $M$ to $N$ is a family $\{
f(x,y)\}_{x \in M, y \in  N }$ such that $f(x,y) :\tau _xM \rightarrow \tau
_yN$ is a Schwartz operator that depends smoothly on its two entries. Let $f ^\ast (x,y):
\tau ^\ast _y N \rightarrow \tau _x ^\ast  M $ be the adjoint of $f(x,y)$. 
Finally, let $X$  be a $M$-valued
semimartingale. We say that a
$N$-valued semimartingale is a solution of the the {\bfi  It\^o stochastic
differential equation}
\begin{equation}
\label{sde stocastic ito}
d Y= f(X,Y) d X
\end{equation}
if for any $\alpha \in  \Omega_2 (N) $, the following equality between It\^o integrals
holds:
\[
\int  \langle \alpha, d Y\rangle =\int  \langle  f ^\ast (X,Y) \alpha,
d X\rangle.
\]
There exists an existence and uniqueness result for the solutions of these stochastic
differential equations analogous to the one for Stratonovich differential equations. 

Given a Stratonovich operator $e$ from $M$ to $N$, there exists a unique Schwartz operator
$f:\tau M\times N\rightarrow\tau N$  defined as follows. Let $ \gamma (t)=(x (t), y (t))\in 
M \times N$ be
 a smooth curve
that verifies $e(x (t) , y (t))(\dot x (t))= \dot y (t) $,
for all $t$. We define  $f (x (t), y (t))\left(L_{\ddot x (t)}\right) :=\left(L_{\ddot y
(t)}\right)$, where the second order differential operators $\left(L_{\ddot x (t)}\right) 
\in\tau_{x\left(  t\right) }M$ and $\left(L_{\ddot y (t)}\right)  \in\tau_{y\left(  t\right)
}N$ are defined  as $\left(L_{\ddot x (t)}\right)\left[
h\right]  :=\frac{d^{2}}{dt^2}h\left(  x\left(  t\right)  \right)  $ and $\left(L_{\ddot y
(t)}\right) \left[ g\right]  :=\frac{d^{2}}{dt^2}g\left(  y\left(  t\right)  \right)  $, for
any $h \in C^\infty(M) $ and 
$g\in C^{\infty}\left(  N\right)$. This relation completely determines $f$ since the
vectors of the form $L_{\ddot x (t)}$ span $\tau_{x\left(  t\right)
}M$. Moreover, the It\^o and Stratonovich equations $\delta Y=
e (X,Y) \delta X $  and $d Y= f(X,Y) d X $ are equivalent, that is, they have the same
solutions.

\medskip

\noindent\textbf{Acknowledgments} 
We thank Michel Emery and Jean-Claude Zambrini for carefully going through the paper and for their valuable comments and
suggestions.  We also thank Nawaf Bou-Rabee and Jerry Marsden for stimulating discussions on stochastic variational integrators. The authors acknowledge partial support from the French Agence National de la Recherche, contract number JC05-41465. J.-A. L.-C. acknowledges
support from the Spanish Ministerio de Educaci\'on y Ciencia grant number
BES-2004-4914. He also acknowledges partial support from MEC grant
BFM2006-10531 and Gobierno de Arag\'on grant DGA-grupos consolidados
225-206. J.-P. O. has been partially supported by a ``Bonus
Qualit\'e Recherche" contract from the Universit\'e de Franche-Comt\'e.


\addcontentsline{toc}{section}{Bibliography}
\begin{thebibliography}{MM99}
\footnotesize
\bibitem[AM78]{fom}  Abraham, R., and Marsden, J.E. 
[1978] \textit{Foundations of Mechanics}. Second  edition,
Addison-Wesley.

\bibitem[A03]{ludwig arnold}Arnold, L. [2003] {\it Random Dynamical Systems}. Springer Monographs
in Mathematics. Springer Verlag.

\bibitem[Ar89]{arnold}  Arnold, V.I. 
[1989] \textit{Mathematical Methods of Classical Mechanics}. Second 
edition. Volume {\bf 60} of
\textit{Graduate Texts in Mathematics}, Springer Verlag. 

\bibitem[B81]{bismut 81}Bismut, J.-M.[1981] {\it M\'ecanique Al\'eatoire}. Lecture Notes in
Mathematics, volume {\bf 866}. Springer-Verlag.

\bibitem[BRO07]{bou 1}Bou-Rabee, N. and Owhadi, H. [2007] Stochastic Variational Integrators. {\it arXiv:0708.2187}.

\bibitem[BRO07a]{bou 2}Bou-Rabee, N. and Owhadi, H. [2007]  Stochastic Variational Partitioned Runge-Kutta Integrators for Constrained Systems. {\it arXiv:0709.2222}.

\bibitem[BJ76]{box jenkins 76}Box, G. E. P. and Jenkins, G. M.[1976] {\it Time Series
Analysis: Forecasting and Control}. Holden-Day.

\bibitem[CH06]{chorin hald}Chorin, A.J. and Hald, O.H.[2006] {\it Stochastic
Tools in Mathematics and Science}. Surveys and Tutorials in the Applied
Mathematical Sciences, volume 1. Springer Verlag.

\bibitem[CW90]{chung williams} Chung, K. L. and Williams, R. J. [1990] {\it Introduction to
Stochastic Integration}. Second edition. Probability and its Applications. Birkh\"auser
Verlag.

\bibitem[CD06]{cresson darses 06}Cresson, J. and Darses, S. [2006]
Plongement stochastique des syst\`emes lagrangiens. {\it   C. R. Math. Acad. Sci. Paris}, 
{\bf 342}, 333-336. 

\bibitem[Du96]{durrett}Durrett, R. [1996] {\it Stochastic Calculus. A Practical Introduction.}
Probability and  Stochastics Series. CRC Press.

\bibitem[E89]{emery} \'Emery, M.  [1989] {\it Stochastic Calculus in
Manifolds}. Springer-Verlag.

\bibitem[E90]{emery transfer}\'Emery, M. [1990] On two transfer principles in
stochastic differential geometry. {\it S\'eminaire de Probabilit\'es}, XXIV,
1988/89, 407-441, Lecture Notes in Math., {\bf 1426}, Springer Verlag.
Correction: {\it S\'eminaire de Probabilit\'es}, XXVI, 633, Lecture Notes in
Math., {\bf 1526}, Springer Verlag,  1992.

\bibitem[G66]{gihman}Gihman, I. I. [1966] Stability of solutions of stochastic differential equations. (Russian)  Limit Theorems Statist. Inference (Russian) pp. 14--45 Izdat. "Fan", Tashkent. English translation: Selected Transl. Statist. and Probability, {\bf 12}, Amer. Math. Soc., pp. 125-154. MR {\bf 40}, number 944.

\bibitem[Ha80]{hasminskii} Hasminskii, R. Z. [1980] {\it Stochastic Stability of Differential Equations.} Translated from the Russian by D. Louvish. Monographs and Textbooks on Mechanics of Solids and Fluids: Mechanics and Analysis, {\bf 7}. Sijthoff and Noordhoff, Alphen aan den Rijn---Germantown, Md. 

\bibitem[H02]{Hsu-Brownian motion}Hsu, E. P. [2002] \textit{Stochastic
Analysis on Manifolds}. Graduate Studies in Mathematics, \textbf{38}. American
Mathematical Society.

\bibitem[IW89]{ikeda watanabe}Ikeda, N. and Watanabe, S. [1989] {\it Stochastic Differential
Equations and Diffusion Processes}. Second edition. North-Holland Mathematical Library, 24.
North-Holland Publishing Co.

\bibitem[I01]{imkeller}Imkeller, P. and  Lederer, C. [2001] Some formulas for Lyapunov exponents and rotation numbers in two dimensions and the stability of the harmonic oscillator and the inverted pendulum. {\it Dyn. Syst.}, {\bf  16}(1), 29-61.

\bibitem[K81]{Kunita}Kunita, H. [1981] Some extensions of It\^{o}'s formula.
{\it Seminaire de probabilit\'{e}s de Strasbourg} XV, 118-141. 
Lecture Notes in Mathematics, {\bf 850}, Springer Verlag. 

\bibitem[LeG97]{legall integration} Le Gall, J.-F. [1997] {\it Mouvement Brownian et Calcul
Stochastique.} Notes de Cours DEA 1996-97. Available at 
\mbox{http://www.dma.ens.fr/~legall/}.

\bibitem[Ne67]{nelson 67} Nelson, E. [1967] {\it Dynamical Theories of Brownian Motion.}
Princeton University Press.

\bibitem[LL76]{landau} Landau, L.D., and Lifshitz, E.M. [1976] \textit{Mechanics}. 
Volume 1 of \textit{Course of Theoretical Physics}. Third Edition. Pergamon Press.

\bibitem[L06]{Xue-Mei}Li, X.-M. An averaging principle for integrable
stochastic Hamiltonian systems. {\it Preprint available at}
http://www.lboro.ac.uk/departments/ma/research/preprints/papers06/06-27.pdf.

\bibitem[M81]{meyer 81}Meyer, P.-A. [1981]
G\'eom\'etrie stochastique sans larmes. Seminar on
Probability, XV (Univ. Strasbourg, Strasbourg, 1979/1980),  44--102,  Lecture
Notes in Math., {\bf 850}, Springer Verlag. 

\bibitem[M82]{meyer 82}Meyer, P.-A. [1982] G\'eom\'etrie diff\'erentielle stochastique.
II.  Seminar on Probability, XVI,
Supplement,  165--207, Lecture Notes in Math., {\bf 921}, Springer-Verlag.

\bibitem[Ok03]{oksendal}{\O}ksendal, B.[2003]{\it Stochastic Differential
Equations}. Sixth Edition. Universitext. Springer-Verlag.

\bibitem[O06]{inverted pendulum}Ovseyevich, A. I. [2006] The stability of an
inverted pendulum when there are rapid random oscillations of the suspension
point. \textit{Journal of Applied Mathematics and Mechanics}, \textbf{70}, 762-768.

\bibitem[O83]{neill 83}O'Neill, B.[1983] {\it Semi-Riemannian Geometry. With
Applications to Relativity} Pure and Applied Mathematics, volume {\bf 103}.
Academic Press. 

\bibitem[LO07]{lo}L\'azaro-Cam\'{\i}, J.-A. and Ortega, J.-P.  [2007]
Reduction of stochastic Hamiltonian systems. {\it In preparation}.

\bibitem[M99]{misawa} Misawa, T. [1999] Conserved quantities and symmetries related to
stochastic dynamical systems. {\it Ann. Inst. Statist. Math.}, {\bf  51}(4),
779--802.

\bibitem[OP04]{leibniz}Ortega, J.-P. and Planas-Bielsa, V. [2004] Dynamics on Leibniz
manifolds. {\it J. Geom. Phys.}, {\bf 52}(1), 1-27.

\bibitem[P90]{protter}Protter, P. [2005] {\it Stochastic Integration and Differential
Equations. A New Approach}. Applications of Mathematics, volume {\bf 21}. Second Edition. Springer-Verlag.

\bibitem[Sch82]{schwartz 82} Schwartz, L. [1982] G\'eom\'etrie diff\'erentielle du 2\`eme ordre,
semi-martingales et \'equations diff\'erentielles stochastiques sur une vari\'et\'e
diff\'erentielle. Seminar on Probability, XVI, Supplement, 1--148, Lecture Notes in Math.,
{\bf 921}, Springer Verlag.

\bibitem[TZ97]{thieullen zambrini 97}Thieullen, M. and  Zambrini, J. C. [1997]  Probability
and quantum symmetries. I. The theorem  of Noether in Schršdinger's Euclidean quantum
mechanics. {\it Ann. Inst. H. PoincarŽ Phys. ThŽor.}, {\bf  67}(3), 297--338.

\bibitem[TZ97a]{thieullen zambrini 97a}Thieullen, M. and  Zambrini, J. C. [1997] Symmetries
in the stochastic calculus of variations. {\it Probab. Theory Related Fields}, {\bf 107}(3),
401--427.

\bibitem[W80]{Watanabe}Watanabe S. [1980] Differential and variation for flow
of diffeomorphisms defined by stochastic differential equations on manifolds
(in Japanese). S\={u}kaiken K\={o}kyuroku, \textbf{391}.

\bibitem[Y81]{yasue 81}Yasue, K. [1981] Stochastic calculus of variations. {\it Journal of
Functional Analysis}, {\bf 41}, 327-340.

\bibitem[ZY82]{zambrini yasue 82}Zambrini, J.-C. and  Yasue, K.[1982] Semi-classical quantum
mechanics and stochastic calculus of variations. {\it Annals of Physics}, {\bf 143},
54-83.

\bibitem[ZM84]{zheng meyer} Zheng, W. A. and Meyer, P.-A. [1984] Quelques r\'esultats de
"m\'ecanique stochastique".  Seminar on
probability, XVIII, 223--244, Lecture Notes in Math., {\bf 1059}, Springer-Verlag. 

\end{thebibliography}
\end{document}